\documentclass[11pt,a4paper,leqno,noamsfonts]{amsart}
 
\usepackage{titlesec}

\titleformat{\part}
  {\centering\normalfont\Large\bfseries} 
  {\textbf{Part \thepart.}} 
  {1em} 
  {} 

\makeatletter
\renewcommand{\l@part}{\@tocline{0}{0pt}{0pc}{}{\bfseries}} 
\makeatother

\renewcommand{\thesection}{\arabic{section}}  
\renewcommand{\thesubsection}{\thesection.\arabic{subsection}}  
\renewcommand{\thesubsubsection}{\thesubsection.\arabic{subsubsection}}  

\titleformat{\section}
  {\normalfont\large\centering\scshape} 
  {\thesection.} 
  {0.5em}
  {}

\titleformat{\subsection}
  {\normalfont\bfseries}
  {\thesubsection.} 
  {0.5em}
  {}

\titleformat{\subsubsection}
  {\normalfont\bfseries}
  {\thesubsubsection.} 
  {0.5em}
  {}

\makeatletter
\renewcommand{\l@section}{\@tocline{1}{0pt}{0em}{1em}{}} 
\makeatother
 
\usepackage[english]{babel}
\usepackage[dvipsnames]{xcolor}
\usepackage{graphicx,pifont}
\usepackage{titlesec}
\usepackage[utopia]{mathdesign}
\usepackage[a4paper,top=2.8cm,bottom=2.8cm,
           left=2.8cm,right=2.8cm,marginparwidth=60pt]{geometry}
\usepackage[utf8]{inputenc}
\usepackage{braket,caption,comment,mathtools,stmaryrd,enumitem} \usepackage[normalem]{ulem}
\usepackage[usestackEOL]{stackengine}
\usepackage{multirow,booktabs,microtype,relsize,soul}
\usepackage[foot]{amsaddr}
\usepackage[colorlinks,bookmarks]{hyperref} %
      \hypersetup{colorlinks,%
            citecolor=britishracinggreen,%
            filecolor=black,%
            linkcolor=cobalt,%
            urlcolor=cornellred}
      \numberwithin{equation}{section}
\usepackage[capitalise]{cleveref}
\Crefname{conjecture}{Conjecture}{Conjectures}
\newenvironment{resume}{%
  \par\medskip
  \noindent{\scshape Résumé.}\enspace\ignorespaces
}{%
  \par
}
\input{macros.tex}
\allowdisplaybreaks

\title[Invariants of nested Hilbert and Quot schemes on surfaces]{Invariants of nested Hilbert and Quot schemes on surfaces}
\author{Nadir Fasola, Michele Graffeo, Danilo Lewański, Andrea T. Ricolfi}

\begin{document}

\maketitle

\begin{abstract}
Let $(S,p)$ be a smooth pointed surface. In the first part of this paper we study motivic invariants of punctual nested Hilbert  schemes attached to $(S,p)$ using the Hilbert--Samuel stratification. We compute two infinite families of motivic classes of punctual nested Hilbert schemes, corresponding to nestings of the form $(2,n)$ and $(3,n)$. As a consequence, we are able to give a lower bound for the number of irreducible components of $S_p^{[2,n]}$ and $S_p^{[3,n]}$.

In the second part of this paper we characterise completely the generating series of Euler characteristics of all nested Hilbert and Quot schemes. 
This is achieved via a novel technique, involving differential operators modelled on the enumerative problem, which we introduce.
From this analysis, we deduce that in the Hilbert scheme case the generating series is the product of a rational function by the celebrated Euler's product formula counting integer partitions. In higher rank, we derive functional equations relating the nested Quot scheme generating series to the rank one series, corresponding to nested Hilbert schemes.

\begin{resume}
Soit $(S,p)$ une surface lisse pointée. Dans la première partie de cet article, nous étudions des invariants motiviques des schémas de Hilbert ponctuels emboîtés associés à $(S,p)$, en utilisant la stratification de Hilbert--Samuel. Nous calculons deux familles infinies de classes motiviques de schémas de Hilbert ponctuels emboîtés, correspondant à des emboîtements de la forme $(2,n)$ et $(3,n)$. En conséquence, nous obtenons une borne inférieure pour le nombre de composantes irréductibles de $S_p^{[2,n]}$ et de $S_p^{[3,n]}$.

Dans la seconde partie de cet article, nous caractérisons complètement les séries génératrices des caractéristiques d’Euler de tous les schémas de Hilbert et touse les schémas Quot emboîtés.
Cela est obtenu au moyen d’une nouvelle technique, que nous introduisons, faisant intervenir des opérateurs différentiels modélisés sur le problème énumératif considéré.
À partir de cette analyse, nous déduisons que, dans le cas des schémas de Hilbert, la série génératrice est le produit d’une fonction rationnelle par la célèbre formule du produit d’Euler comptant les partitions d’entiers. En rang supérieur, nous obtenons des équations fonctionnelles reliant les séries génératrices des schémas Quot emboîtés à la série de rang un, correspondant aux schémas de Hilbert emboîtés.
\end{resume}
\end{abstract}

\renewcommand{\contentsname}{}
\tableofcontents

\section{Introduction}
Let $S$ be a smooth complex quasiprojective surface, $p \in S$ a closed point. This paper studies motivic invariants of the \emph{punctual nested Hilbert scheme}
\[
S_p^{[n_1,\ldots,n_\ell]},
\]
the fine moduli space of flags $Z_1 \subset \cdots \subset Z_\ell$ of finite subschemes of $S$ of fixed length $\chi(\OO_{Z_i}) = n_i$, such that $Z_\ell$, and hence every $Z_i$, is entirely supported at $p$. The higher rank version of this moduli space, the \emph{nested Quot scheme}, is also analysed (cf.~\Cref{sec:euler-quot}).

The paper is divided in two parts. The first part, consisting of Sections \ref{sec:useful-strata}--\ref{sec:nested_3n}, computes the motive of all punctual nested Hilbert schemes of the form $S_p^{[i,n]}$, for $i=2,3$ (see \Cref{sec:K0-results} for an outline). The second part, consisting of Sections \ref{subsec:nestedpart}--\ref{sec:euler-quot}, computes various generating series of Euler characteristics of punctual nested Hilbert schemes and Quot schemes (see, respectively, Sections \ref{sec:intro-hilb} and \ref{sec:nested-quot-intro} for an outline). We give more details on our main results in the next two sections.

\subsection{\texorpdfstring{Results in $K_0(\Var_{\BC})$}{}} 
\label{sec:K0-results}
The `universal motivic invariant' (or, simply, the \emph{motive}) of the punctual Hilbert scheme $S_p^{[n]}$, namely its class in the Grothendieck ring of varieties $K_0(\Var_{\BC})$, was computed by G\"{o}ttsche \cite{Gottsche-motivic}. The resulting motivic generating function
\[
\mathsf{Hilb}^\bullet(t) = \sum_{n\geqslant 0} \,[S_p^{[n]}]t^n = 1+t+(\BL+1)t^2+\cdots
\]
can be written in the explicit form
\begin{equation}\label{eq:gottsche}
    \mathsf{Hilb}^\bullet(t) = \prod_{j\geqslant 1}\,\frac{1}{1-\BL^{j-1}t^j},
\end{equation}
where $\BL$ is the \emph{Lefschetz motive}, the class of $\BA^1$ in $K_0(\Var_{\BC})$. 

We consider, for $i\geqslant 0$, the generating functions
\begin{equation}
\label{eqn:nested-series}
\mathsf{Hilb}^{[i,\bullet]}(t) = \sum_{n\geqslant i}\,[S_p^{[i,n]}]t^n\,\in\,K_0(\Var_{\BC})\llbracket t \rrbracket.
\end{equation}
Note that the generating functions so far introduced do not depend on the pair $(S,p)$, as they all equal the series attached to the pair $(\BA^2,0)$.

Our first main result is the determination of the series \eqref{eqn:nested-series} in the first nontrivial cases, corresponding to $i=2,3$.\footnote{The case $i=1$ is covered by G\"{o}ttsche's formula, see  \cite[Thm.~5.1]{Gottsche-motivic} or, more precisely, the proof of \cite[Thm.~5.3]{Gottsche-motivic}.}

\begin{thm}[Corollaries \ref{cor:k=2}, \ref{cor:k=3}]
\label{thm:intro-k=2and3}
There are identities
\begin{align*}
\mathsf{Hilb}^{[2,\bullet]}(t)&=[\BP^1]\mathsf{Hilb}^\bullet(t)+[\BP^1]\frac{t(\BL-1)-1}{1-\BL t}, \\
\mathsf{Hilb}^{[3,\bullet]}(t)&=[\BP^2]\mathsf{Hilb}^\bullet(t)-\frac{\mathsf h(t)}{(1 - \BL t) (1 - \BL^2 t^2)},
\end{align*}
where $\mathsf h(t) \in \BZ[\BL,t] \subset K_0(\Var_{\BC})[t]$ is the polynomial
\[
\mathsf h(t) = [\BP^2] - (\BL^3-1) t - (\BL^3-1)[\BP^1] t^2 + \BL^2(\BL^3 - 1) t^3 - \BL^2 t^4 .
\]
\end{thm}

The proof goes by a series of stratification arguments according to the pairs of compatible Hilbert--Samuel functions of a nesting $Z_1 \subset Z_2$, and exploits crucially a number of Zariski locally trivial fibrations (cf.~\Cref{subsec:HS-strata}) whose targets are the homogeneous loci in the \emph{unnested} punctual Hilbert scheme $S_p^{[n]}$. This allows us to explicitly determine the motives of $S_p^{[i,n]}$ for $i=2,3$ and arbitrary $n$, see \Cref{thm:motivepunctual} and \Cref{motive-punctual-3n}.

Along the way, we are able to confirm that, for $n\geqslant 4$, the number of irreducible components of maximal dimension in $S_p^{[2,n]}$ is equal to $\floor{n/2}$, see \Cref{cor:irrcomp2n}. This is weaker than the statement proved in \cite[Cor.~7.5]{BULOIS}, which also includes purity, but we provide a new proof. Similarly, we prove in \Cref{cor:irrcomp3n} that the number of irreducible components of maximal dimension in $S_p^{[3,n]}$ is equal to   {  $\floor{\frac{n^2+4}{12}} $} for $n\geqslant 5$.

\subsection{Euler characteristic results} \label{subsec:intronestedhilb}
The second part of this paper is devoted to the computation of the topological Euler characteristic of general nested Hilbert schemes and Quot schemes. 

\subsubsection{Nested Hilbert schemes}
\label{sec:intro-hilb}
Our starting point is the basic observation that the positive integer 
\[
\chi(S_p^{[n_1,\ldots,n_\ell]})
\]
agrees with the number of $\ell$-tuples $(\lambda_1,\ldots,\lambda_\ell)$ of Ferrers diagrams\footnote{Ferrers diagram are just the `visual' counterpart of integer partitions, see \Cref{sec:nested-Ferrers}.} such that $\lambda_i \subset \BN^2$ has size $\lvert \lambda_i\rvert = n_i$ for all $i=1,\ldots,\ell$, and $\lambda_1 \subset \cdots \subset \lambda_\ell$.
Our main object of study is the generating function\footnote{The `$\mathsf F$' in the notation, here and throughout, stands for `flag'.} $\FFZ_{\bk}(q)\in\BZ\llbracket q\rrbracket$ defined, for an arbitrary $s$-tuple $\bk=(k_1,\ldots,k_s) \in\BZ^s_{\geqslant 0}$, by the formula
\begin{equation}
\label{eqn_FZk}
\FFZ_{\bk}(q)=\sum_{n\geqslant 0} \chi(S_p^{[n,n+k_1,n+k_1+k_2,\ldots,n+\sum_{i=1}^sk_i]})q^n.
\end{equation}
When $s=0$, we recover the series 
\[
\mathsf Z(q) = \chi \mathsf{Hilb}^\bullet (q) = \sum_{n \geqslant 0}\chi(S_p^{[n]})q^n = 1+q+2q^2+3q^3+5q^4+\cdots 
\]
which by the celebrated formula proved by Euler \cite[Chapter 16]{Euler_partitions} can be expressed as the infinite product\footnote{Note that this expression can be recovered from \Cref{eq:gottsche} setting $\BL=1$.}
\begin{equation}
\label{classical-Z}
\mathsf Z(q) = \prod_{j\geqslant 1}\frac{1}{1-q^j}.
\end{equation}
The most general structural result we obtain is the following.

\begin{thm}[{\Cref{thm:ZDnested}}]
\label{thm:intro:Z_k}
Let $s>0$ be a positive integer. Given a sequence of nonnegative integers $\bk=(k_1,\ldots,k_s)\in \BZ^s_{\geqslant 0}$ summing up to $K$, the generating series $\FFZ_{\bk}(q)$ satisfies the relation
\begin{equation}
\label{eqn:Zk/Z}
    \frac{\FFZ_{\bk}(q)}{\ZZ(q)} = \frac{\PP_{\bk}(q)}{\prod_{j=1}^{K} (1 - q^j)},
\end{equation}
where $\PP_{\bk}(q)$ is a polynomial in $q$ of degree bounded by $\frac{5}{4}K^2 - 2K + 1$. In particular, the ratio $\FFZ_{\bk}/\ZZ$ is a rational function in $q$ with only roots of unity as poles.
\end{thm}
The proof goes, roughly, as follows. We only explain here the case $s=1$, as it already contains the key insights. See \Cref{subsec:g.f.tantinest} for the general case. In the case $s=1$, after setting $D=k_1$, we observe that in order to determine the coefficients of 
\begin{equation}
\label{ZD-series-intro}
\FFZ_D(q) = \sum_{n\geqslant 0}\chi(S_p^{[n,n+D]})q^n,
\end{equation}
we have to figure out, given a \emph{skew Ferrers diagram} $\lambda$ of size $D$, the number of ways to add it to a Ferrers diagram of size $n$ such that the result is a (non-skew) Ferrers diagram of size $n + D$. Let $A_{n-D}(\lambda)$ be this number, and set
\[
\FFZ_\lambda(q) = \sum_{n\geqslant D} A_{n-D}(\lambda) q^n.    
\]
To $\lambda$, we associate (cf.~\Cref{eq:defTlambda}) the differential operator 
\[
T_{\lambda} = 
\sum_{j \geq 0} \prod_{k=1}^M \left[ \left( \frac{\mathrm{d}}{\mathrm{d}y_{j+\sum_{i=k+1}^M \ell_i,v_k}} - \frac{\mathrm{d}}{\mathrm{d}y_{j + \sum_{i=k+1}^M \ell_i,v_k+1}} \delta_{k < M}\right) \prod_{p=0}^{\ell_k-2} \frac{\mathrm{d}}{\mathrm{d}y_{j+\sum_{i=k}^{M} \ell_i - p - 1, 0}}\right],
\]
depending on formal variables $\vec{y} = (y_{j,h}\,|\,(j,h) \in \BZ \times \BN)$, where $M=M_\lambda \in \BZ_{\geq 1}$ is an integer depending on the shape of $\lambda$ (cf.~\Cref{not:SK-paths}). Furthermore, to $D$ we associate a $\vec{y}$-decorated version of the partition function \eqref{classical-Z}, namely (cf.~\Cref{def:F-partition-function})
\[
F^{[D]}(\vec{y}; q) = 
    \prod_{j\geq 1} 
    \left( 
    y_{j,0} + y_{j,1}q^j + y_{j,1}y_{j,2} q^{2j} + \dots +
    \left(\prod_{k=1}^D y_{j,k}\right) \sum_{m_j \geq D} q^{jm_j}
    \right).
\]
We prove in \Cref{thm:skewferret} the key relation
\[
\FFZ_\lambda(q) = T_{\lambda}.F^{[D]}(\vec{y}; q)\Bigg|_{\vec{y} = 1}.
\]
The series $\FFZ_D$, cf.~\Cref{ZD-series-intro}, is then reconstructed (cf.~\Cref{cor:TD-applied-to-FD}) as
\begin{equation}
\label{eqn:ZD-intro}
\FFZ_D(q) = T_D.F^{[D]}(\vec{y}; q) \Bigg{|}_{\vec{y} = 1},
\end{equation}
where $T_D = \sum_{\lvert \lambda \rvert = D} T_\lambda$, the sum running over all skew Ferrers diagrams of size $D$. From here, to get the $s=1$ specialisation of \Cref{thm:intro:Z_k}, one is reduced to proving its `$\lambda$-analogue', stating that the ratio $\FFZ_\lambda / \ZZ$ has the same structure of the right hand side of \Cref{eqn:Zk/Z}: this is done in \Cref{prop:Zlambda}.

In \Cref{sec:in-action} we include explicit calculations for $D=2,3$ in order to show concretely how the process works. This yields, for instance, the identity
\[
\frac{\FFZ_3(q)}{\ZZ(q)} = \frac{3 - q - q^2}{(1-q)(1-q^2)(1-q^3)}.
\]
See also Sections \ref{sec:some-polynomials} and \ref{sec:some-ratios} for some explicit ratios $\FFZ_D/\ZZ$ and $\FFZ_{\bk}/\ZZ$. Furthermore, we determine the values at 0 and 1 of the polynomials $\mathsf P_\lambda$, $\mathsf P_D$ and $\mathsf P_{\bk}$ in \Cref{lem:Plambda0110}, \Cref{lemma:specialisation} and \Cref{rmk:specialisation} respectively. 

It is an interesting question whether there exists a motivic refinement of the techniques leading to the proof of \Cref{thm:intro:Z_k}. We leave this for future research.

\subsubsection{Nested Quot schemes}
\label{sec:nested-quot-intro}

In \Cref{sec:euler-quot} we establish a `higher rank' version of the constructions and results obtained in 
Sections~\ref{sec:g.f.}--\ref{sec:proof-B}. Geometrically, this corresponds to determining the Euler characteristic $\chi_r^{[\bn]}$ of the \emph{punctual nested Quot scheme}
\[
\Quot_r^{[\bn]}(S)_p
\]
attached to a smooth pointed surface $(S,p)$, where $\bn = (n_1,\ldots,n_\ell) \in \BZ^{\ell}_{\geq 0}$ is a nondecreasing sequence of integers. This space parametrises flags of subsequent quotients
\[
\begin{tikzcd}
\OO_{S}^{\oplus r} \arrow[two heads]{r} & F_\ell\arrow[two heads]{r} & F_{\ell-1}\arrow[two heads]{r} & \cdots \arrow[two heads]{r} & F_1
\end{tikzcd}
\]
where $F_i$ is a 0-dimensional coherent sheaf on $S$ such that $\chi(F_i) = n_i$ for $i=1,\ldots,\ell$ and $F_\ell$ is entirely supported at $p$.

In the unnested setup (i.e.~when $\ell=1)$, we observe, cf.~\Cref{lemma:sZ-geometric}, that the series
\begin{equation}
\label{seriesQ(q,s)-intro}
\QQ(q,s)=\sum_{r\ge 0}\sum_{n \geq 0}\chi_r^{[n]}q^ns^r\in\BZ\llbracket q,s\rrbracket
\end{equation}
satisfies the relation 
\begin{equation}
\label{eqn:Q(qs)-intro}
\QQ(q,s) = \frac{1}{1-\mathsf Z(q)s}.
\end{equation}
We view this as a first indication that the higher rank theory is \emph{fully determined} by the rank $1$ theory. This heuristic phenomenon has been observed in several places, see for instance the motivic computations \cite{MR_nested_Quot, CR_framed_motivic, CRR_higher_rank} and the recent results in  Donaldson--Thomas theory \cite{FMR_higher_rank, fasola2023tetrahedron,FT_1}. We prove that this heuristic holds true in the nested case as well for surfaces, by providing explicit closed formulas, as we now explain.

{First of all, we introduce the higher rank versions of the series given in \Cref{sec:intro-hilb}, namely 
\begin{align*} 
    \FFQ_{r,\bk}(q)&=\sum_{n\ge 0}\chi_r^{[n,n+k_1,n+k_1+k_2,\ldots,n+\sum_{i=1}^sk_i]}q^n \in\BZ\llbracket q\rrbracket,& \text{for } r\in\BZ_{\ge 0},  \bk=(k_1,\ldots,k_\ell) \in \BZ_{\geq 0}^\ell,\\
    \FFQ_{\bk}(q,s)&=\sum_{r\ge 0}\FFQ_{r,\bk}(q)  s^r \in\BZ\llbracket q,s\rrbracket, & \text{for }  \bk=(k_1,\ldots,k_\ell) \in \BZ_{\geq 0}^\ell,\\ 
    \FFQ(q,s,\boldsymbol{\vec v} ) &=\sum_{\ell\ge 1}\sum_{\bk\in\BZ^{\ell}_{\ge0}}\FFQ_{\bk}(q,s) \boldsymbol{\vec v}^{\bk}\in\BZ\llbracket q,s,\boldsymbol{\vec v} \rrbracket,\\
\FFZ(q,\boldsymbol{\vec v})  &=\sum_{\ell\ge 1}\sum_{\bk\in\BZ^{\ell}_{\ge0}}\FFZ_{\bk}(q)\boldsymbol{\vec v}^{\bk}\in\BZ\llbracket q,\boldsymbol{\vec v}\rrbracket.
\end{align*}
where $\boldsymbol{\vec{v}}=(v_i\,|\,i \in \BZ_{\geq 1})$ is a countable collection of variables, and $\boldsymbol{\vec v}^{\bk}=\prod_{i=1}^\ell v_i^{k_i}$ is a monomial. The series $\FFZ_{\bk}$ was defined in \Cref{eqn_FZk} and clearly agrees with $\FFQ_{1,\bk}$. 

We prove the following relations.

\begin{thm}
[{\Cref{thm:Z_Dr}, \Cref{cor:FQ/Z^r}}, \Cref{cor:FQ/Z^r:multinested}]
\label{thm:intro-Z_Dr}
There is an identity of formal power series
\[
\mathsf{FQ}(q,s,\boldsymbol{\vec v}) =
\frac{1}{ 1 - \mathsf{FZ}(q,\boldsymbol{\vec v})s}\in\BZ\llbracket q,s,\boldsymbol{\vec v}\rrbracket.
\]
Moreover, for any choice of $r\in\BZ_{\ge0}$ and   vector ${\bk} = (k_1, \dots, k_\ell)\in\BZ^\ell_{\ge0}$, we have an identity
\begin{equation}
\frac{\FFQ_{r,\bk}(q)}{\ZZ(q)^r}  
= \frac{P_{\bk}(q)}{\prod_{j=1}^K (1 - q^j)^{m_j}}
\end{equation}
in $\BZ\llbracket q \rrbracket$, where $K = \sum_{i=1}^s k_i$, $m_j = m_j(\bk)$ are nonnegative integers, and $P_{\bk}(q)$ is a polynomial that depends on the vector $\bk$. In particular, the ratio $\FFQ_{r,\bk}(q) / \ZZ(q)^r$ is a rational function in $q$. In the special case of length one nestings, for any choice of $D \in \BZ_{\geq 0}$ we have an identity
\[
\frac{\FFQ_{r,D}(q)}{\ZZ(q)^r}=\frac{\mathsf P_{r,D}(q)}{\prod_{j=1}^D(1-q^j)^{\min\left(r, \floor{\frac{D}{j}}\right)}}, \]
where $\mathsf P_{r,D}(q)\in\BZ[q]$ is a polynomial.
\end{thm}

Furthermore, we prove that the specialisation
\[
\FFQ(q,s,v) = \FFQ(q,s,\boldsymbol{\vec v})\big|_{{v_1=v},\,v_{>1} = 0}
\]
is entirely determined by the application of certain differential operators to the series $\QQ(q,s)$, cf.~Equation~\eqref{seriesQ(q,s)-intro}, as expressed by the following result.

\begin{thm}[{\Cref{thm:Z_Dr3}}]
\label{thm:Z_Dr3-intro}
There is an identity of formal power series
\begin{equation}
\FFQ(q,s,v)=\exp{\left(\FFZ(q,v)-1\right)}\Big{|}_{v^k \mapsto v^ks^k\frac{\dd^k}{\dd s^k}}.\QQ(q,s)
\in\BZ\llbracket q,v,s\rrbracket,
\end{equation}
where $\FFZ(q,v) = \FFZ(q,\boldsymbol{\vec v})\big|_{{v_1=v},\,v_{>1} = 0}$ and the substitution of variables is meant after the expansion of the exponential as a power series in $v$.
\end{thm}

By means of a general power structure argument, in \Cref{sec:power-structure} we explain how to pass from invariants of \emph{punctual} Quot schemes to invariants of global Quot schemes on an arbitrary smooth surface $S$. As an example, if $S=\dP_6$ is the 6-th del Pezzo surface, we are able to calculate the rank 6 invariant
\[
\chi\left( \Quot_6^{[6,6+6]}(\dP_6)\right)=120806108165466.
\]
}

\subsection{Geometry of nested Hilbert schemes and related works}
The Hilbert scheme of points $S^{[n]}$ on a smooth surface $S$ is a smooth $2n$-dimensional variety, by a celebrated result of Fogarty \cite{Fogarty_Hilb}. See Nakajima's book \cite{Nakajima} for an overview of this vaste subject.
The geometry of \emph{nested} Hilbert schemes on surfaces is still largely unexplored, for instance the quest for (elementary) components, and the understanding of their scheme structure. See \cite{ALESSIONESTED,Rasul-irr-nested} for results on irreducibility for several types of nestings and \cite[Prop.~3.7]{ALESSIONESTED} for reducibility examples. It is proved in \cite{UPDATES} that there exist $n_1<\cdots < n_5$ such that $S^{[1,n_1,\ldots,n_5]}$ has generically nonreduced components. See Cheah's work \cite{MR1616606} for a characterisation of smoothness (in all dimensions), and \cite{Lissite-quot} for a similar characterisation in higher rank. Motivic aspects of the theory of Hilbert and Quot schemes, besides G\"{o}ttsche's seminal work \cite{Gottsche-motivic}, can be found in \cite{dCM_motives,GLMHilb,mozgovoy2019motivic,ricolfi2019motive,MR_nested_Quot,MR-hyperquot,MR18,double-nested-1}. For K-theoretic and geometric aspects of nested Hilbert schemes, and their connection to moduli of flags of sheaves and framed sheaves, the reader can consult \cite{BULOIS, Negut_flags,BFT_flags,cazzaniga2020framed,Mon_double_nested} and the references therein. On the cohomological side, for classical Hilbert schemes of surfaces we refer to \cite{Gott2,ESHilb,Nakajima}, whereas for combinatorial aspects on 3-folds (in the language of \emph{plane partitions}) we refer to \cite{Atkin,MacMahon-statistics} and finally we refer to the recent work \cite{GMMR2} for a combinatorial advance in all dimensions.

\subsection{Conventions}
\label{sec:conventions}
We work over the field $\BC$ of complex numbers throughout. For a scheme $Y$ and locally closed subschemes $U_i \subset Y$, we abuse notation writing `$Y = \coprod_i U_i$' to mean that $(U_i)_i$ form a \emph{stratification} of $Y$, in the sense that the canonical morphism $\coprod_i U_i \to Y$ is bijective. As we deal with \emph{punctual} Hilbert (and Quot) schemes only, the reader is free to replace the (fixed, but arbitrary) pointed smooth surface $(S,p)$ with the local model $(\BA^2,0)$. In fact, our analysis is carried out in the local case, exploiting without mention the identification $S_p^{[\bn]} = (\BA^2)^{[\bn]}_0$. If $\mathsf c$ is any `condition', by $\delta_{\mathsf c} \in \set{0,1}$ we denote the associated Kronecker delta, taking the value $1$ if $\mathsf c$ is satisfied, and $0$ otherwise.

\subsection*{Acknowledgements} The authors thank Paolo Lella, Yuze Luan, Sergej Monavari, Riccardo Moschetti, and Alessio Sammartano for useful discussions.  N.~F., D.~L. and A.~R. are partially supported by the PRIN project 2022BTA242 “Geometry of algebraic structures: moduli, invariants, deformations”. M.~G., D.~L. and A.~R. are members of the GNSAGA - INdAM.

We are particularly grateful to the anonymous referee for the careful review and for suggesting many useful improvements, including \Cref{rmk:alternative} and the final version of \Cref{thm:intro-Z_Dr}.

\part*{Motivic calculations}
We prove in this part of the paper \Cref{thm:intro-k=2and3} from the Introduction. In other words, we will compute the generating functions \eqref{eqn:nested-series} of the motives $S_p^{[i,n]}$ for $i=2,3$. This will be achieved as a combination of Corollaries \ref{cor:k=2} and \ref{cor:k=3}. 

In the next section we recall the formal definition of the moduli space $S^{[\bn]}$ we are interested in, along with the basic material (Hilbert--Samuel functions, apolarity and the Grothendieck ring of varieties) that will be used in the main computations.

Before diving into the computation of the motives of $S_p^{[i,n]}$ for $i=2,3$, we write down in the \Cref{sec:useful-strata} an explicit auxiliary stratification of the (unnested) punctual Hilbert scheme $S_p^{[n]}$, according to the value $\bh_I(2) \in \set{1,2,3}$ of the Hilbert--Samuel function of a non-curvilinear ideal $I \subset \BC[x,y]$ of colength $n$.

\section{Background material}

\subsection{Nested Hilbert schemes}
Let $X$ be a complex quasiprojective scheme. Fix an integer $\ell > 0$ and an $\ell$-tuple of nonnegative integers $\bn = (n_1,\ldots,n_\ell)$ such that $n_1 \leqslant \cdots \leqslant n_\ell$. The $\bn$-nested Hilbert functor attached to $X$ is the functor $\Sch_{\BC}^{\opp} \to \Sets$ which assigns to a $\BC$-scheme $B$ the set of flags
\[
\begin{tikzcd}
    \CZ_1 \arrow[hook]{r} & \CZ_2 \arrow[hook]{r} & \cdots  \arrow[hook]{r} & \CZ_\ell
\end{tikzcd}
\]
of $B$-flat closed subschemes $\CZ_i \into X \times B$ such that $\CZ_i \to B$ is finite of relative length $n_i$ for $i=1,\ldots,\ell$. This functor is represented by a quasiprojective scheme $X^{[\bn]}$, called the $\bn$-\emph{nested Hilbert scheme}, see e.g.~\cite[Theorem 4.5.1]{Sernesi-Deformations} and \cite{Kleppe}. The case $\ell=1$ recovers the classical Hilbert scheme $X^{[n]}$. 

\begin{notation}
Given an an $\bn$-\emph{nesting} $Z_1\subset \cdots \subset Z_\ell$ of finite subschemes of $X$, the corresponding $\BC$-valued point of $X^{[\bn]}$ will be denoted by $[Z_1\subset \cdots \subset Z_\ell]$, or by $[\mathscr I_1 \supset \cdots \supset \mathscr I_\ell]$, if $\mathscr I_j \subset \OO_X$ is the ideal sheaf cutting out $Z_j \subset X$.
\end{notation}

\subsubsection{Punctual locus}
By construction, there is a closed immersion
\begin{equation}
\label{eqn:iota-inclusion}
\begin{tikzcd}
X^{[\bn]} \arrow[hook]{r}{\iota} & \displaystyle\prod_{i=1}^\ell X^{[n_i]},
\end{tikzcd}
\end{equation}
cut out by the nesting condition. Fix a closed point $p \in X$. The \emph{punctual nested Hilbert scheme} attached to $(X,\bn,p)$ is the projective closed subscheme $X^{[\bn]}_p \subset X^{[\bn]}$ whose points correspond to $\ell$-tuples $(Z_1,\ldots,Z_\ell)$ such that $Z_\ell$ (and hence every $Z_i$) is entirely supported at $p$ --- such schemes are known as as \emph{fat points}. The scheme $X^{[\bn]}_p$ is defined as the scheme-theoretic fibre of the composition
\[
\begin{tikzcd}
X^{[\bn]} \arrow{r}{\pr_\ell} &  X^{[n_\ell]} \arrow{r}{} & \Sym^{n_\ell}(X)
\end{tikzcd}
\]
over the cycle $n_\ell\cdot p$, where the first map is the map \eqref{eqn:iota-inclusion} followed by the $\ell$-th projection and $X^{[n_\ell]} \to \Sym^{n_\ell}(X)$ is the Hilbert--Chow morphism (see \cite{Rydh1} for its construction in the most general setup, or \cite[Thm.~2.16]{Bertin}, expanding on Grothendieck's original argument \cite[Sec.~6]{Grothendieck_Quot}; see also \cite{Fantechi-Ricolfi-motivic,Fantechi-Ricolfi-structural}). It can be proved (along the same lines of the argument in \cite[Sec.~2.1]{ricolfi2019motive}) that, if $X$ is smooth of dimension $d$ at $p$, then the punctual nested Hilbert scheme attached to $(X,\bn,p)$ is (non-canonically) isomorphic to the punctual nested Hilbert scheme attached to $(\BA^d,\bn,0)$, where $0 \in \BA^d$ is the origin.

\subsubsection{Curvilinear locus}
\label{sec:puttane-curvilinee}
Let $X$ be a smooth $d$-dimensional variety. The easiest example of a fat point in $X$ is provided by a \emph{curvilinear subscheme}, namely a fat point $Z \subset X$ such that $\OO_Z(Z) \cong \BC[t]/t^{n+1}$, for some $n \geqslant 0$. For each $p \in X$, the curvilinear schemes supported at $p$ form an open subscheme 
\begin{equation}
\label{eqn:curvilinear-locus}
\begin{tikzcd}
\mathsf C_p^{[n]} \arrow[hook]{r} & X^{[n]}_p.
\end{tikzcd}
\end{equation}
The \emph{curvilinear locus} in the nested setup is the open subscheme defined as the fibre product
\[
\begin{tikzcd}
\mathsf C_p^{[\bn]} \MySymb{dr}\arrow{d}\arrow[hook]{r} & X^{[\bn]}_p\arrow{d}{\pr_\ell} \\
\mathsf C_p^{[n_\ell]} \arrow[hook]{r} & X^{[n_\ell]}_p
\end{tikzcd}
\]
whose points correspond to $\bn$-nestings $Z_1\subset \cdots \subset Z_\ell$ such that $Z_\ell$ (and hence every $Z_i$) is curvilinear.

\begin{remark}
\label{rmk:curvilinear-loci}
We record here, for future use, the following observations.
\begin{itemize}
\item [\mylabel{curv-1}{(1)}] If $\dim X=2$, the curvilinear locus \eqref{eqn:curvilinear-locus} is dense in $X_p^{[n]}$, and admits a Zariski locally trivial fibration with target $\BP^1$ and fibre $\BA^{n-2}$, see \cite[Prop.~IV.1.1]{Bria1}.
\item [\mylabel{curv-2}{(2)}] For each $\bn = (n_1,\ldots,n_\ell)$, the $\ell$-th projection
\[
\begin{tikzcd}
\mathsf C_p^{[\bn]} \arrow{r}{\pr_\ell} & \mathsf C_p^{[n_\ell]}    
\end{tikzcd}
\]
is a bijective morphism. More generally, projecting onto the $k$-th factor yields an affine bundle $\pr_k\colon \mathsf C_p^{[\bn]} \to \mathsf C_p^{[n_k]}$ with fibre $\BA^{n_\ell-n_k}$.
\item [\mylabel{curv-3}{(3)}] Just as the punctual loci, the curvilinear loci do not depend on the pair $(X,p)$, but only on the ambient dimension $d=\dim X$. Since we will fix $X=\BA^2$ very soon, we will just denote them by $\mathsf C^{[\bn]}$, viewed inside $(\BA^2)^{[\bn]}_0$.
\end{itemize}
\end{remark}

\subsubsection{Tangent space}\label{subsub:tgnested}
Fix a closed point $[Z] \in X^{[\bn]}$ corresponding to an $\ell$-tuple $(Z_1,\ldots,Z_\ell)$. The differential of \eqref{eqn:iota-inclusion} at $[Z]$ is an injection of $\BC$-vector spaces
\[
\begin{tikzcd}
\mathsf T_{[Z]}X^{[\bn]} \arrow[hook]{rr}{\dd \iota|_{[Z]}} && \displaystyle\prod_{i=1}^\ell\Hom_{\OO_X}(\mathscr I_i,\OO_{Z_i}),
\end{tikzcd}
\]
where $\mathscr I_i \subset \OO_X$ denotes the ideal sheaf of $Z_i \into X$. This map realises the tangent space $\mathsf T_{[Z]}X^{[\bn]}$ as the locus of $\ell$-tuples of tangent vectors
\[
\delta = (\delta_1,\ldots,\delta_\ell) \in \prod_{i=1}^\ell\Hom_{\OO_X}(\mathscr I_i,\OO_{Z_i})
\]
such that all squares in the diagram
\begin{equation}
\label{diag:tangents}
\begin{tikzcd}[row sep=large,column sep=large]
\mathscr I_\ell\arrow[hook]{r}\arrow{d}{\delta_\ell} & \mathscr I_{\ell-1}\arrow[hook]{r}\arrow{d}{\delta_{\ell-1}} & \cdots \arrow[hook]{r} &  \mathscr I_2\arrow[hook]{r}\arrow{d}{\delta_2} & \mathscr I_1\arrow{d}{\delta_1} \\
\OO_{Z_\ell}\arrow[two heads]{r} & \OO_{Z_{\ell-1}}\arrow[two heads]{r} & \cdots \arrow[two heads]{r} & \OO_{Z_2}\arrow[two heads]{r} & \OO_{Z_1}
\end{tikzcd}
\end{equation}
commute, see e.g.~\cite[Prop.~4.5.3(i)]{Sernesi-Deformations}.

\subsubsection{The  Bia{\l{}}ynicki-Birula decomposition of the punctual nested Hilbert scheme}\label{subsec:BBdeco}
Let $X$ be a smooth quasiprojective variety and let $p\in X$ be a closed point with ideal sheaf $\Fm_p\subset \OO_X$. For a fat point $Z\subset X$ supported at $p$, and given an integer $k \geqslant 0$, we set
\[ 
\OI_Z^{\geqslant k} = \OI_Z\cap \Fm_p^k \subset \OI_Z \quad \mbox{ and }\quad \OO_Z^{\geqslant k} = (\Fm_p^k + \OI_Z)/\OI_Z \subset \OO_Z. 
\]

\begin{definition}
\label{def:negativetangents}
Let $[Z]\in X^{[\bn]}_p \subset X^{[\bn]}$ be a point. The \textit{nonnegative part of the tangent space} $\mathsf T_{[Z]}   X^{[\bn]}$ is the vector subspace
\[
\mathsf T_{[Z]}^{\geqslant0}X^{[\bn]}   = \Set{\delta\in\mathsf T_{ [Z]} X^{[\bn]}  | \delta_i(\OI_{Z_i}^{\geqslant k}) \subset \OO_{Z_i}^{\geqslant k}\mbox{ for all }k \in\BN\mbox{ and for }i=1,\ldots,\ell} \subset \mathsf T_{[Z]}X^{[\bn]}. 
\] 
\end{definition}

Notice that nonnegative tangent vectors can be understood as commutative diagrams of the form \eqref{diag:tangents} where $\delta_i\in \mathsf T_{[Z_i]}^{\geqslant 0}X^{[n_i]}$ for all $i=1,\ldots,\ell$.

When $X=\BA^d$, the nonnegative part of the tangent space can be interpreted as the tangent space to the so-called  Bia{\l{}}ynicki-Birula decomposition of the punctual nested Hilbert scheme, whose definition we recall now. We restrict to the case $d=2$. Consider the lift of the diagonal action of the torus $\mathbb G_m=\Spec \BC[s,s^{-1}]$ on $\BA^2$, given by homotheties, to the nested Hilbert scheme $(\BA^2)^{[\bn]}_0$. \textit{The  Bia{\l{}}ynicki-Birula decomposition of the punctual nested Hilbert scheme} is the quasiprojective scheme $(\BA^2)_0^{[\bn],+}$ representing the functor $\Sch_{\BC}^{\opp} \to \Sets$ sending
\[
B\mapsto 
\Set{\overline{\mathbb G}_m \times B \xrightarrow{\varphi} \nested{\bn}|\varphi\mbox{ is }\mathbb G_m\mbox{-equivariant}
}
\]
where, by convention, one sets $\overline{\mathbb G}_m=\Spec \BC[s^{-1}]$, see \cite{ELEMENTARY,UPDATES}.
 
\begin{prop}[{\cite[Thm.~4.11]{ELEMENTARY}}]\label{rem:nonneg-punctual}
Let $[Z]\in (\BA^2)_0^{[\bn],+}$ be a point. Then  
    \[
    \mathsf T_{[Z]}(\BA^2)_0^{[\bn],+}= \mathsf T_{[Z]}^{\geqslant0}   \nested{\bn} .
    \]     
\end{prop}

{ 
\begin{remark}\label{rmk:novo}
It is worth mentioning that \cite[Thm.~4.11]{ELEMENTARY} only involves nestings of two schemes ($\ell=2$) in its original form. However, the proof can be immediately adapted to the general case, see also \cite{UPDATES}.

We also stress that, if $I\subset R $ is a homogeneous ideal of finite colength $n$, then the tangent space $\mathsf T_{[I]}(\BA^2)^{[n]}\cong\Hom_R(I,R/I)$ is naturally $\BZ$-graded. In particular, by \Cref{def:negativetangents}, positive tangent vectors to the nested Hilbert schemes are naturally identified with $\ell$-tuples of positive tangent vectors making the diagram \eqref{diag:tangents} commute. Moreover, by standard arguments, see \cite{multigraded}, in this setting one has
\[
\mathsf T_{[I]}((\BA^2)_0^{[n],+})^{\BG_m}= \mathsf T_{[I]}^{=0}   \nested{n},
\]
and similarly for the nestings of homogeneous ideals.
\end{remark}
}

\begin{remark}
\label{rem:canschem}
Set-theoretically $(\BA^2)_0^{[\bn],+}$ coincides with the punctual nested Hilbert scheme $(\BA^2)_0^{[\bn]}$, in the sense that there is a bijective morphism between the reductions of these two schemes. This fact will be exploited in the paper for the computation of the motive of some punctual nested Hilbert schemes.  

On the other hand, the Hilbert--Samuel strata $H_{\bh}$, recalled in \Cref{subsec:HS-strata} below, are unions of connected components of the punctual Hilbert scheme  endowed with the schematic structure induced by the Bia{\l{}}ynicki-Birula decomposition. This induces a canonical schematic structure on the Hilbert--Samuel strata   and it allows us to define morphisms with target $H_{\bh}$ canonically; this, in turn, provides an upper bound for the dimension of the fibre of the initial ideal morphism \eqref{initial-ideal-map} and, crucially, it allows us to compute the motive of the punctual Hilbert scheme (and of its nested versions) in terms of Hilbert--Samuel strata.
\end{remark}

\subsection{Tools from commutative algebra}
\label{sec:comm-algebra}
In this subsection we recall a couple of useful concepts from commutative algebra: Hilbert--Samuel functions and apolarity.

\subsubsection{Hilbert--Samuel functions}
\label{sec:HS-functions}
The Hilbert--Samuel function is a classical invariant attached to  graded modules. We refer the reader to \cite[Sec.~5.1]{EISENBUD} and \cite[Lect.~20]{HARRIS} for more details.

From now on, we work with $X=\BA^2$, with coordinate ring $R=\BC[x,y]$, and we denote by $\Fm = (x,y) \subset R$ the maximal ideal of the origin $0 \in \BA^2 = \Spec R$. We endow $R$ with the standard grading $\deg(x)=\deg(y)=1$. The $i$-th graded piece of $R$ will be denoted $R_i$. A fat point $Z \subset \BA^2$ supported at the origin corresponds, algebraically, to a local Artinian (finite type) $\BC$-algebra $(A,\Fm_A)$ with residue field $\BC$, namely
\begin{equation}\label{eqn:local-artin}
(A,\Fm_A) = (R/I,\Fm/I).
\end{equation}

\begin{definition}
Let $M=\bigoplus_{i\in\BZ} M_i$ be a finite graded $R$-module. The \emph{Hilbert--Samuel function} $\bh_M$ attached to $M$ is the function
\[
\begin{tikzcd}[row sep = tiny]
\BZ\arrow[r,"\bh_M"]& \BN\\
        i \arrow[r,mapsto] & \dim_{\BC} M_i. 
\end{tikzcd}
\]
Let $I\subset R$ be the ideal of a fat point $Z = \Spec R/I \subset \BA^2$. Form the local Artinian $\BC$-algebra $(A,\Fm_A)$ as in \eqref{eqn:local-artin}. The \emph{Hilbert--Samuel function} of $A$ (that we shall denote $\bh_I$ or by $\bh_Z$ throughout) is defined to be the Hilbert--Samuel function of the graded $R$-module
\begin{equation}\label{eqn:grading-A-module}
\mathsf{gr}_{\Fm_A}(A) =\bigoplus_{i\geqslant 0} \mathfrak m_A^{i}/\mathfrak m_A^{i+1}.
\end{equation}
\end{definition}
 
We record in the next remark a few useful observations, that we exploit to set up some more notation. 

\begin{remark}
\label{rmks-on-HS}
Suppose $A=R/I$ corresponds to a fat point $Z \subset \BA^2$. The function $\bh_I$ has the following properties:
\begin{itemize}
\item [\mylabel{HF1}{(i)}] It vanishes for almost all values (i.e.~it has finite support), and satisfies $\bh_{I}(0)=1$. Therefore we represent $\bh_I$ as a finite string of integers, namely $\bh_I = (1,h_1,\ldots,h_t)$, where $t$ is the largest index such that $h_t \neq 0$.
\item [\mylabel{HF2}{(ii)}] The value $\bh_{I}(1) = \dim_{\BC} \Fm_A/\Fm_A^2$ is nothing but the \emph{embedding dimension} of $Z$, namely the smallest dimension of a scheme containing $Z$ and smooth at the unique point of $Z$. Since we work in dimension 2, we have $\bh_{I}(1) \in \set{0,1,2}$. Note that $Z$ is curvilinear if and only if $\bh_I(1)\leqslant 1$.
\item [\mylabel{HF3}{(iii)}] The length $\chi(\OO_Z)=\dim_{\BC}A$ can be computed as the finite sum $\lvert \bh \rvert = 1+\sum_{i > 0} h_i$.
\item [\mylabel{HF4}{(iv)}] Let $\In I \subset R$ be the initial ideal of $I$ (i.e.~the ideal generated by the initial forms of the polynomials in $I$, cf.~\cite[Sec.~5.1]{EISENBUD}). Then the isomorphism of graded $R$-algebras $\mathsf{gr}_{\Fm_A}(A)\cong R/\In I$ shows that $\bh_{\In I} = \bh_I$.
\end{itemize}
\end{remark}

\subsubsection{Hilbert--Samuel strata and homogeneous loci}
\label{subsec:HS-strata}
We are ready to introduce the \emph{Hilbert--Samuel strata}. We refer the reader to \cite{multigraded} for more details.

\begin{definition}
Given a function $\bh\colon\BZ\rightarrow \BN$ with finite support, the \emph{Hilbert--Samuel stratum} $H_{\bh}$  is the (possibly empty) locally closed subset
\[
H_{\bh}=\Set{[I]\in \nested{\left|\bh\right|}_0 | \bh_I=\bh}\subset \nested{\left|\bh\right|}_0.
\]
We endow it with the reduced induced subscheme structure.
\end{definition}

The \emph{Hilbert--Samuel stratification}
\[
\nested{n}_0 = \coprod_{\lvert \bh \rvert = n} H_{\bh} 
\]
will be our main tool in later computations (cf.~\Cref{rem:canschem}). Its nested version, exploited in later sections, is a straightforward generalisation.

The standard scaling action $\BG_m\times \BA^2\rightarrow \BA^2$ lifts to the Hilbert scheme $\nested{n}_0$. 
Under this action, the Hilbert--Samuel strata are $\BG_m$-invariant locally closed subsets of the punctual Hilbert scheme. Hence, the action of the torus $\BG_m$ restricts to an action on $H_{\bh} $, for all $\bh\colon \BZ\rightarrow \BN$. Set theoretically, the $\BG_m$-fixed locus of the action on $H_{\bh} $ consists of \emph{homogeneous} ideals (with fixed Hilbert--Samuel function $\bh$). This locus defines a closed subset \cite[Prop.~1.5]{multigraded}
\[
\OH_{\bh}  \subset H_{\bh} ,
\]
which{, whenever not specified differently,} we endow with the reduced induced subscheme structure. 

\begin{prop}[{\cite[p.~773]{8POINTS}}]\label{prop:initial morphism}
Fix $\bh = (1, h_1, \ldots , h_t)$. There is a morphism
\begin{equation}\label{initial-ideal-map}
    \begin{tikzcd}
        H_{\bh} \arrow[r,"\pi_{\bh}"] &\OH_{\bh}
    \end{tikzcd}
    \end{equation}
which, on closed points, sends an ideal to its initial ideal.
\end{prop}
 
In the 2-dimensional case, which is of interest to us, the initial ideal morphism recalled in \Cref{prop:initial morphism} behaves particularly well, as explained in the following theorem collecting two celebrated results by Brian\c{c}on and Iarrobino.

\begin{theorem}[{\cite[Thm.~III.3.1]{Bria1}} and {\cite[Thm 1]{Iarropunctual}}]
\label{thm:briarro}
Let $H_{\bh}\subset \nested{\lvert \bh \rvert}_0$ be the stratum associated to the Hilbert--Samuel function $\bh$. Denote by $d,s>0$  the integers such that
     \[
     \bh = ( 1,2,\ldots,d, h_{d},\ldots,h_{d+s-1},0,\ldots), 
     \]
with $h_d<d+1$ and $h_{i}\geqslant h_{i+1}$, for $i\geqslant d$. Then, $H_{\bh}\subset \nested{\lvert \bh \rvert}_0$ is smooth, connected and of dimension

\begin{equation}\label{eq:briarro}
\dim H_{\bh} = \lvert \bh \rvert-d-\sum_{i\geqslant d} \binom{h_{i-1}-h_{i}}{2}.     
 \end{equation}
    Moreover, the homogeneous locus $\OH_{\bh}$ has dimension 
    \[
    \dim \OH_{\bh} = \sum_{i\geqslant d} (h_{i-1}-h_i+1)(h_{i}-h_{i+1}),
    \]  
    and the initial ideal morphism $
    \begin{tikzcd}
        H_{\bh} \arrow[r,"\pi_{\bh}"] &\OH_{\bh}
    \end{tikzcd}$ defines a Zariski locally-trivial fibration  with fibres  affine spaces.
\end{theorem}

\begin{remark}\label{rmk:novo2}
\Cref{thm:briarro} is a powerful tool allowing one to compute the fibres of the map $\pi_{\bh}$ in the unnested case. This in turn gives important information about the motives of the Hilbert--Samuel strata, see \Cref{prop:briarromotivico}. 
On the other hand, also in the nested case there is an analogue of the Hilbert--Samuel stratification and of the initial ideal morphism in \Cref{prop:initial morphism}, see \cite{2step} for an explicit presentation. Since there is no nested version of \Cref{thm:briarro} available, in order to compute the fibres of the initial ideal morphism one has to apply \Cref{rem:nonneg-punctual} and to show that the positive tangents are unobstructed, see \Cref{rmk:novo}. This is what we will do in \Cref{lemma:Y1112}.
\end{remark}

\subsubsection{Apolarity}
\label{sec:apolarity}
A powerful tool for parametrising families of ideals is \emph{Macaulay duality}, also known as \emph{apolarity} \cite{Macaulay}. We now present some results on apolarity that we shall need through the rest of the paper. The reader can consult \cite{8POINTS,Hilb_11,EMSALEM,Geramita,Iarrobook} for more details. Working in characteristic $0$ is crucial in this section.

Let us set
\[
R=\BC[x,y], \quad R^*=\BC[x^*,y^*].
\]
We view $R^*$ as an $R$-module via the action 
\[
\begin{tikzcd}[row sep =tiny]
    R\times R^*\arrow[r]& R^* \\
    (x^{\alpha}y^{\beta},p(x^*,y^*))\arrow[r,mapsto] & \displaystyle\frac{\partial^{\alpha+\beta} }{\partial^{\alpha} x^*  \partial^{\beta} y^*}p,
\end{tikzcd}\]
where $\alpha,\beta \in\BZ_{\geqslant 0}$. This induces, for every  $k\geqslant 0$, a perfect pairing
\[
\begin{tikzcd}
R_k\times R_k^* \arrow{r} & R_0^* = \BC
\end{tikzcd}
\]
and, consequently, a notion of orthogonality. 
\begin{definition}[{\cite{Hilb_11}}]
 An \emph{inverse system} is a graded vector subspace of $R^*$ closed under differentiation. If $S\subset R^*$ is a finite subset consisting only of homogeneous forms, then the inverse system generated by $S$ is the smallest graded subspace $\langle S\rangle\subset R^*$ containing $S$ and closed under differentiation. The \emph{apolar ideal} attached to $S$ is  
\[
S^\perp=\Set{r\in  R  | r\cdot \langle S\rangle =0 }\subset R.
\]
If $I\subset R$ is a homogeneous ideal, its associated inverse system is 
\[
I^\perp =\Set{r^*\in  R^*  | I \cdot r^* =0 }\subset R^*.
\]
\end{definition} 
\begin{remark}
Notice that if $I\subset R$ is a homogeneous ideal, then $I^\perp\subset R^*$ is a graded subspace closed under differentiation. Conversely, every graded vector subspace $T\subset R^*$ closed under differentiation is orthogonal to the homogeneous ideal $T^\perp\subset R$.
Moreover, if $V\subset R^*_j$ is a vector subspace, then
\[
\dim_\BC (V^\perp)_j =\dim_\BC {R^*_j}/ V,
\]
which yields an isomorphism of graded vector spaces $R/I\cong I^\perp$, see \cite[Sec.~2]{8POINTS}.
\end{remark}
 
\begin{prop}
\label{prop:gorapolar}
Let $f \in  R^*$ be a homogeneous form of degree
$d\geqslant 0$. Then the Hilbert--Samuel function of the ideal $\langle f\rangle^\perp$ is
symmetric, i.e.~$\bh_{\langle f\rangle^\perp}(i) = \bh_{\langle f\rangle^\perp}(d - i)$ for all $i=0,\ldots,d$. 

Suppose also that $\bh_{\langle f\rangle^\perp}$ has the form
\[
\bh_{\langle f\rangle^\perp}=(1,h_1,h_2,\ldots,k,1),
\]
then $k=1,2$. If $k=2$, then 
$\dim_{\BC}\Span\left(\frac{\mathrm{d}f}{\mathrm{d}x^*},\frac{\mathrm{d}f}{\mathrm{d}y^*}\right)=2$. If $k=1$, then $h_1=1$ and $f=\ell^d$, for some $\ell \in R^*_1$.  
\end{prop} 

\begin{proof}
The first part of the statement is \cite[Lemma 2.12]{Iarrobook}. The second part of the statement is \cite[Prop.~10]{CARLINI}. The number of essential variables computed in loc.~cit.~is identified with $h_1$ by \cite[Lemma 6]{CARLINI} which in turn equals $k$ by the symmetry described in the first part of the statement.
\end{proof} 

\subsection{Grothendieck ring of varieties}
\label{sec:K0(Var)}

The Grothendieck ring of complex varieties $K_0(\Var_{\BC})$ is the free abelian group generated by isomorphism classes $[X]$ of $\BC$-varieties, modulo the \emph{scissor relations}, namely the relations 
\[
[X] = [Y] + [X \setminus Y]
\]
whenever $Y \subset X$ is a closed subvariety. The group structure agrees, on generators, with the disjoint union of varieties and has neutral element $0 = [\emptyset]$. Any finite type $\BC$-scheme $X$ has a well-defined class $[X]$ in $K_0(\Var_{\BC})$, namely the motive of its reduction. The fibre product defines a ring structure on $K_0(\Var_{\BC})$, with neutral element $1 = [\Spec \BC]$. Every constructible subset $Z \subset X$ of a variety $X$ has a well-defined motivic class $[Z]$, independent upon the decomposition of $Z$ into locally closed subsets. A bijective morphism $X \to Y$ of schemes induces an identity $[X] = [Y]$ in $K_0(\Var_{\BC})$. Another fundamental rule is the \emph{fibration property}, which says that if $X \to Y$ is a Zariski locally trivial fibration, with fibre $F$, then $[X] = [Y][F]$ in $K_0(\Var_{\BC})$.

The \emph{Lefschetz motive} is the class
\[
\BL = [\BA^1]\in K_0(\Var_{\BC}).
\]
It allows one, for instance, to express the class of the projective space $\BP^e$ as $[\BP^e] = \sum_{0\leq i\leq e}\BL^i$.

We will not need this in what follows, but the Lefschetz motive is of remarkable importance in birational geometry. For instance, the quotient $K_0(\Var_{\BC}) / (\BL)$ is isomorphic to the monoid algebra generated by stable birational equivalence classes of varieties \cite{zbMATH02069674}; it was proved in \cite[Thm.~7.5]{Galkin-Shinder} that if $\BL$ is not a 0-divisor then a very general smooth cubic 4-fold is irrational. However, $\BL$ turned out to be a 0-divisor \cite{zbMATH06837514}, and this made it possible to define the notion of $\BL$-equivalence for algebraic varieties, and to relate it, conjecturally, with derived equivalence \cite{zbMATH06941785}.

The ring $K_0(\Var_\BC)$ has the following universal property. Suppose $A$  is a ring and $w(-)$ is an $A$-valued invariant of algebraic varieties, such that 
\begin{itemize}
\item [$\circ$] $w(\pt)=1$,
\item [$\circ$] $w(\emptyset) = 0$,
\item [$\circ$] $w(X\times Y) = w(X)w(Y)$ for every two varieties $X$ and $Y$, and
\item [$\circ$] $w(X) = w(Y) + w(X\setminus Y)$ for every variety $X$ and closed subvariety $Y \subset X$. 
\end{itemize}  
Then $w$ is called a \emph{motivic invariant} of algebraic varieties, and there is precisely one ring homomorphism $w\colon K_0(\Var_\BC) \to A$ sending the class $[X]$ of a variety $X$ to $w(X) \in A$. Ring homomorphisms out of $K_0(\Var_{\BC})$ are called \emph{motivic measures}, or generalised Euler characteristics, or realisations \cite{DenefLoeser1,LooijengaMM}. Therefore the knowledge of the motive of a scheme implies the knowledge of all of its motivic invariants. 

An immediate consequence of \Cref{thm:briarro} is the following equality in the Grothendieck ring of varieties.

\begin{prop}
\label{prop:briarromotivico} 
Let $H_{\bh}\subset \nested{\lvert \bh \rvert}_0$ be the stratum associated to the Hilbert--Samuel function 
     \[
     \bh = ( 1,2,\ldots,d, h_{d},\ldots,h_{t},0,\ldots). 
     \]
Then there is an identity 
\[
  [H_{\bh}]=[\OH_{\bh}] \BL^{\binom{d}{2} -\sum_{i\geqslant d } \frac{(h_{i-1}-h_{i})(h_{i-1}+h_{i} -2 h_{i+1}-1 )}{2} -h_{i+1} }
\] 
in $K_0(\Var_\BC)$.
\end{prop}

\section{A useful stratification of \texorpdfstring{$S_p^{[n]}$}{}}
\label{sec:useful-strata}
The most natural way to compute the motive  of a variety is to stratify it by locally closed subschemes whose  motives are known. In \cite{MOTIVES} the main strategy consisted in stratifying the punctual Hilbert scheme $(\BA^d)^{[n]}_0$ according to the first entry of the Hilbert--Samuel vector, namely the embedding dimension. Here we analyse the punctual Hilbert scheme $S_p^{[n]}=(\BA^2)^{[n]}_0$ of a smooth surface by using a finer stratification: we stratify the locus parametrising ideals $I\subset R=\BC[x,y]$ with $\bh_I(1)=2$ (the complement of the curvilinear locus) according to the value $\bh_I(2)$. 
This yields the decomposition
\[
\nested{n}_0=\mathsf C ^{[n]} \amalg H_1^{[n]} \amalg H_2^{[n]} \amalg H_3^{[n]},
\]
where $\mathsf C^{[n]}\subset \nested{n}_0$ is the curvilinear locus (cf.~\Cref{sec:puttane-curvilinee}), and
\[
H_i^{[n]}=\Set{[I]\in \nested{n}_0\smallsetminus \mathsf C^{[n]} | \bh_I(2)= i}.
\]
The main strategy to compute the motive of these strata consists in analysing the possible Hilbert--Samuel functions occurring in each stratum. 
  
The curvilinear locus satisfies
\begin{equation}
    \label{eq:curv}
    [\mathsf C^{[n]}]=[\BP^1]\BL^{n-2}
\end{equation}
by \Cref{rmk:curvilinear-loci}\ref{curv-1}.
We move to the computation of the motives of the loci parametrising ideals $I$ of embedding dimension $\bh_I(1)=2$.

Before we start, we fix the following notation.

\begin{notation}
If $H \subset \nested{n}_0$ is a $\BG_m$-invariant locally closed subscheme, we denote by $\OH \subset H$ its $\BG_m$-fixed locus. Set-theoretically, it consists of points $[I] \in H$ such that the corresponding ideal $I\subset R$ is \emph{homogeneous}.
\end{notation}

\subsection{The stratum \texorpdfstring{$H_1^{[n]}$}{}}\label{sec:stratum_H1}
A point $[Z]\in H_1^{[n]}$ corresponds to a scheme of length $\chi(\OO_Z)=n\geqslant 4$ and Hilbert--Samuel function function
\[
\bh=(1,2,1^{n-3}).
\]
We distinguish the special case $n=4$ from the case $n\geqslant 5$.

If $n=4$, then any homogeneous ideal $I$ corresponding to a point of $\OH_1^{[n]}$ is of the form $I=\langle f,x^*,y^* \rangle^\perp $, for some degree 2 homogeneous form  $f\in R_2^*$. Therefore, we have $\mathscr H_1 ^{[n]}\cong\BP^2$. On the other hand, if $n\geqslant 5 $, a homogeneous ideal $I$ with $\bh_I\equiv \bh$ is of the form $I=\langle\ell^{n-2},x^*,y^* \rangle^\perp$, for some linear form $\ell\in R_1^*$, by \cite[Prop.~2.16]{Hilb_11}. Therefore, we have $\OH_1^{[n]}\cong\BP^1$. Summing up and applying \Cref{prop:briarromotivico}, we obtain the identities 
\begin{equation}
\label{eq:H1}
    [H_1^{[n]}] =\begin{cases}
[\BP^2]&\mbox{ for }n=4,\\
    [\BP^1]\BL^{n-3}&\mbox{ for }n\geqslant 5.
\end{cases}
\end{equation}

\subsection{The stratum \texorpdfstring{$H_2^{[n]}$}{}}\label{sec:stratum_H2}  
In order to compute the motive of the stratum $H_2^{[n]}$, we stratify again. We start by the observation that only few Hilbert--Samuel functions can occur on $H_2^{[n]}$. Indeed, for every $s \leqslant \left\lfloor\frac{n-1}{2}\right\rfloor$ we get one possible Hilbert--Samuel function in $H_2^{[n]}$, namely
\[
\bh_s = (1,2^s,1^{n-2s-1}).
\]
We leave the dependence on $n$ implicit to avoid making the notation too cumbersome. Therefore we have a stratification
\[
H_2^{[n]}=\coprod_{s=2}^{\left\lfloor \frac{n-1}{2}\right\rfloor} H_{\bh_s}.
\]
In the next sections of the paper we will focus on the punctual nested Hilbert scheme on a smooth surface $S$. We shall need a finer stratification than the one we just presented. More precisely, we will stratify $H_{\bh_s}$ according to the number $i \in \set{1,2}$ of distinct roots of the unique quadric generator of $\In I$,  
i.e.~we define
\begin{equation}
\label{def:H^i_s}
H_{\bh_s}^i=\Set{[I]\in H_{\bh_s} |
\begin{array}{cc}
\mbox{the unique generator of }(\In I)_2 \\
\mbox{has }i\mbox{ distinct roots}
\end{array}
}.
\end{equation}
Notice that $H_{\bh_s}=H_{\bh_s}^1\amalg H_{\bh_s}^2$ is a locally closed stratification.

\begin{prop}
\label{prop:stratquad}
In $K_0(\Var_\BC)$, we have identities
\[
[H_{\bh_s}^1]=\begin{cases}
  [\BP^1]\BL^{n-5}  &\mbox{ if }s=\frac{n-1}{2},\\
    [\BP^1][\BP^1]\BL^{n-5}&\mbox{ if }s=\frac{n-2}{2},\\  
  [\BP^1] \BL^{n-4}   &\mbox{ if }1<s<\frac{n-2}{2}, 
\end{cases}
\qquad\qquad\qquad
[H_{\bh_s}^2]=\begin{cases}
  \BL^{n-3}  &\mbox{ if }s=\frac{n-1}{2}, \\
     [\BP^1]\BL^{n-3}&\mbox{ if }s=\frac{n-2}{2} ,\\  
  [\BP^1]\BL^{n-3}  &\mbox{ if }1<s < \frac{n-2}{2}.  
\end{cases}
\] 
\end{prop}

\begin{proof} 
By \Cref{prop:briarromotivico}, we only need to compute the motives $[\OH_{\bh_s}^i]$, for $i=1,2$. 

We start from the case $s=\frac{n-1}{2}$. In particular $n=2s+1$ is odd. The homogeneous ideals parametrised by  $\OH_{\bh_s}$ have one of the following forms 
\[
I_1=(\ell_1^2)+\Fm^{s+1} \qquad\mbox{ and }\qquad I_2=(\ell_1\ell_2)+\Fm^{s+1},
\]
for some linear forms $\ell_1\not= \ell_2\in \BP R_1$. Precisely, we have $[I_i]\in\OH_{\bh_s}^i$, for $i=1,2$. This observation identifies $\OH_{\bh_s}^1$ with the image of the second Veronese embedding of $\BP R_1$ in $\BP R_2$, i.e. with a conic in $\BP R_2\cong\BP^2$. Hence, we have
\[
[\OH_{\bh_s}^1] = [\BP^1],\qquad [\OH_{\bh_s}^2] =  \BL^2.
\]

We move now to the case $s=\frac{n-2}{2}$. In particular $n=2s+2$ is even. The homogeneous ideals corresponding to points of $\OH_{\bh_s}$ are uniquely determined by their unique generators in degree 2 and $s+1$. This locus defines a $\BP^1$-bundle over $\BP R_2$ whose points correspond to homogeneous ideals of the following possible forms
\[
I_{1}= (\ell_1^2,f)+\Fm^{s+2}  \qquad\mbox{ and }\qquad I_2= (\ell_1\ell_2,g)+\Fm^{s+2} ,
\]
for some $\ell_1,\ell_2\in R_1$, $f\in \BP (R_{s+1}/ (\ell_1^2)_{s+1})\cong\BP^1$  and $g\in \BP(R_{s+1}/ (\ell_1\ell_2)_{s+1})\cong\BP^1$. In particular, we have $[I_i] \in\OH_{\bh_s}^i$. By arguing as above, we get
\[
[\OH_{\bh_s}^1] = [\BP^1]^2,\qquad [\OH_{\bh_s}^2] = [\BP^1]\BL^2.
\]
Finally, we focus on the general case $n-2s-1$. The homogeneous locus $\OH_{\bh_s}$ is a $\BP^1$-bundle over $\BP R_1^*$. Indeed, there is a Zariski locally trivial fibration
\[
\begin{tikzcd}[row sep=.75em]
    \OH_{\bh_s}\arrow[r] & v_{1,n-s-1}(\BP R_1^*) & \hspace{-1cm}\subset ~~\BP R_{n-s-1}^*\\
{[I]} \arrow[mapsto,r] & {[(I^\perp)_{n-s-1}]},
\end{tikzcd}
\]
where $v_{1,n-s-1}$ is the $(n-s-1)$-th Veronese embedding of $\BP^1$, and whose fibre over a point $[\langle\ell^{n-s-1}\rangle]\in\BP R_1^*$ is $\BP(R_2^*/\langle\ell^2\rangle)\cong\BP^1$.

Moreover, $[\OH_{\bh_s}^2]=[\BP R_1^*]\BL$, as any homogeneous ideal $[I]\in\OH_{\bh_s}^2$ is of the form $\langle\ell_1^{n-s-1},\ell_2^{s}\rangle^\perp$, with $(\ell_1,\ell_2)\in(\BP R_1^*)^{\times2}\smallsetminus\Delta$. Indeed, by \Cref{prop:gorapolar}, the respective derivatives of the dual generators of $I$ must both span a 1-dimensional vector space, for otherwise the ideal would not contain a quadric with distinct roots. Thus, we also get $[\OH_{\bh_s}^1]=\BP^1$, which concludes the proof.
\end{proof}

Summing up we obtain the following corollaries.

\begin{corollary} 
In $K_0(\Var_\BC)$, we have identities
\[
[H_{\bh_s}]=\begin{cases}
  [\BP^2]\BL^{n-5}  &\mbox{ if }s=\frac{n-1}{2}, \\
    [\BP^2][\BP^1]\BL^{n-5}&\mbox{ if }s=\frac{n-2}{2}, \\  
  [\BP^1][\BP^1]\BL^{n-4}  &\mbox{ if }1<s<\frac{n-2}{2}. 
\end{cases}
\]
\end{corollary}

\begin{corollary}
\label{cor:H2}
Define, for $i=1,2$,
\[
H_2^{[n],i}=\Set{[I]\in H_2^{[n]} |
\begin{array}{cc}
\mbox{the unique generator of }(\In I)_2 \\
\mbox{has }i\mbox{ distinct roots}
\end{array}
}.
\] 

Then, if $n=2k$ is even, we have the identities
\begin{align*}
 {[H_2^{[n],1}]}&= (1+(k-2)\BL)[\BP^1]\BL^{n-5}, \\  
 {[H_2^{[n],2}]}&= (k-2)[\BP^1]\BL^{n-3}.
\end{align*}
Instead, if $n=2k+1$ is odd, we have the identities
\begin{align*}
{[H_2^{[n],1}]} &= (1+(k-2)\BL)[\BP^1]\BL^{n-5},\\
{[H_2^{[n],2}]} &= (1+(k-2)[\BP^1])\BL^{n-3}.
\end{align*}
As a consequence, we get 
    \[ [H_2^{[n]}]=\BL^{n-5}\left(\left(\left\lfloor \frac{n}{2}\right\rfloor-2\right)\BL^3+(n-4)\BL^2+\left(\left\lfloor \frac{n}{2}\right\rfloor-1\right)\BL +1\right).
    \] 
\end{corollary}
 
 \subsection{The stratum \texorpdfstring{$H_3^{[n]}$}{}}\label{sec:stratum_H3} 
\begin{prop}
\label{prop:H3}
In $K_0(\Var_{\BC})$, we have the identities
\[
[H_3^{[n]}] = 
\begin{cases}
[\nested{n}_0] - \left(1+(k-2)[\BP^1]\BL + [\BP^1]\BL^2\right)[\BP^1]\BL^{n-5} & \mbox{ if }n=2k, \\ \\
[\nested{n}_0]- \left([\BP^2]+(k-2)[\BP^1]^2\BL +[\BP^1]^2\BL^2 \right)\BL^{n-5} & \mbox{ if }n=2k+1.
\end{cases}
\]
\end{prop}
\begin{proof}
The proof consists of a direct computation on the equality  
\[
[ H_3^{[n]}]=[\nested{n}_0]-[\mathsf C^{[n]} ]- [H_1^{[n]}]- [H_2^{[n]}].\qedhere
\]
\end{proof}
 \begin{lemma}\label{lemma:intercoeff}
    Fix a positive integer $n>3$ and write the motive $[\nested{n}_0]$ as a polynomial in $\BL$, \[
     [\nested{n}_0]=p^{[n]} (\BL) =\BL^{n-1}+\sum_{i=0}^{n-2}a_i^{[n]}\BL^i \in \BZ[\BL]\subset K_0(\Var_{\BC}),
     \] 
     for some $a_i^{[n]}\in\BZ$. Then, we have
\begin{equation}
\label{eqn:L-coefficients}
\begin{split}
a_{n-2}^{[n]}&=\left\lfloor\frac{n}{2}\right\rfloor,\\
a_{n-3}^{[n]}&=\left\lfloor\frac{n(n-6)}{12}\right\rfloor+\left\lfloor\frac{n-1}{2}\right\rfloor + 1 {  =\left\lfloor\frac{n^2+4}{12}\right\rfloor}.         
\end{split}
\end{equation}
\end{lemma}
 \begin{proof}
     The strategy of the proof consists in detecting which Hilbert--Samuel strata contribute to the coefficients in the statement. We do this by applying \Cref{thm:briarro}. The strata contributing to $a_{n-2}^{[n]}$ are all contained in $\mathsf C^{[n]} \amalg H_1^{[n]}\amalg H_2^{[n]}$. Thus, we can recover $a_{n-2}^{[n]}$ from \Cref{eq:curv}, \eqref{eq:H1} and  \Cref{cor:H2}. Now, by equation \eqref{eq:briarro}, given a Hilbert--Samuel function
      \[
     \bh = ( 1,2,\ldots,d, h_{d},h_{d+1},\ldots ), 
     \] 
with $h_d<d+1$ and $h_{i}\geqslant h_{i+1}$, for $i\geqslant d $, the motive $[H_{\bh}]$ has degree $n-3$
if and only if
 \[
n-3=\dim H_{\bh} = \lvert \bh \rvert-d-\sum_{i\geqslant d} \binom{h_{i-1}-h_{i}}{2}.
 \]
The contribution to $a_{n-3}^{[n]}$ coming from $\mathsf C^{[n]}\amalg H_1^{[n]}\amalg H_2^{[n]}$ is $0+1+n-4$.  Now, we only need to compute the contribution coming from $H_3^{[n]}$. To do so,  we can put $d=3$ and $n\geqslant 6$. Then, a Hilbert--Samuel stratum $H_{\bh}$ contributes to $a_{n-3}^{[n]}$ if and only if the identity
\[
\sum_{i=d}^{d+u-1}\binom{h_{i-1}-h_{i}}{2}=0
\]
is satisfied. This condition translates into the fact that the $\bh$ must be of the form
\[
\bh = ( 1,2,\underbrace{3,3\ldots,3}_r,\underbrace{2,2\ldots,2}_s,\underbrace{1,1\ldots,1}_t), 
\]
for some $r,s,t>0$ such that $3(r+1)+2s+t=n$.
Therefore $a^{[n]}_{n-3}$ agrees with the number of partitions of $n$ with exactly 3 parts. The second formula in \eqref{eqn:L-coefficients} then follows from~\cite[Prop.~1.8.6]{RPStanley-VolI} {  and standard manipulations}.
\end{proof}

\section{The motive of \texorpdfstring{$S_p^{[2,n]}$}{}}\label{sec:motive_2n}

In this subsection we shall compute the motives $[S_p^{[2,n]}]\in K_0(\Var_{\BC})$ and the associated generating function $\mathsf{Hilb}^{[2,\bullet]}(t)$, thus proving the first identity in \Cref{thm:intro-k=2and3}. 

Recall that the closed points of $S_p^{[2,n]} = \nested{2,n}_0$ correspond to pairs of nested subschemes $Z_1\subset Z_2\subset \BA^2$ such that $Z_2$ (and hence $Z_1$) is entirely supported at $0 \in \BA^2$.

\begin{theorem}
\label{thm:motivepunctual}  
Let $S$ be a smooth quasiprojective surface, $p \in S$ a closed point. For every $n \in \BZ_{\geqslant 2}$, there is an identity
\[
[S_p^{[2,n]}]= [S_p^{[n]}] [\BP^1] - \BL^{n-1} [\BP^1]
\]   
in $K_0(\Var_{\BC})$.
\end{theorem}
\begin{proof}
The basic observation is that every fat point of length $2$ is curvilinear.
Let us denote by $\mathsf C^{[2,n]}\subset \nested{2,n}_0$ the nested curvilinear locus and by $\mathsf C^{[n]}\subset \nested{n}_0$ the standard curvilinear locus, as in \Cref{sec:puttane-curvilinee}. Then, by \Cref{rmk:curvilinear-loci}\ref{curv-1}-\ref{curv-2}, we have  
\[
[\mathsf C^{[2,n]}]=[\mathsf C^{[n]}] = [\BP^1]\BL^{n-2}.
\] 
On the other hand, the morphism 
\[
\begin{tikzcd} 
\nested{2,n}_0\smallsetminus  \mathsf C^{[2,n]} \arrow{r}{\pr_2}&\nested{n}_0\smallsetminus \mathsf C^{[n]}
\end{tikzcd}
\]
is a trivial $\mathsf C^{[2]}$-bundle and $\mathsf C^{[2]}\cong \BP^1$. Therefore there is an identity
\[
[S_p^{[2,n]}] -  [\mathsf C^{[2,n]}] = [\BP^1] ([S_p^{[n]}]-[\mathsf C^{[n]}]).
\] 
Summarising, we obtain
\begin{align*}
[S_p^{[2,n]}]&=[S_p^{[2,n]}] - [ \mathsf C^{[2,n]}] + [ \mathsf C^{[2,n]}]  =\\
&=[\BP^1] ([S_p^{[n]}] - [ \mathsf C^{[n]}])  +[\BP^1]\BL^{n-2}  \\
&= [\BP^1] ([S_p^{[n]}] -[\BP^1]\BL^{n-2})  +[\BP^1]\BL^{n-2} \\
&= [\BP^1] ([S_p^{[n]}] - \BL^{n-1}),
\end{align*}  
as required.
\end{proof}

In \Cref{motive-punctual-3n} we shall prove the analogue of \Cref{thm:motivepunctual} for $S_p^{[3,n]}$. 

\begin{remark}
    Notice that \Cref{thm:motivepunctual} is consistent with 
    \[
    \dim\, \nested{2,n}_0=\dim \,\nested{n}_0.
    \]
    Indeed, by \cite{Bria1}, the motive  $[\nested{n}_0]$ is a monic polynomial  in $\BL$ of degree $n-1=\dim \,\nested{n}_0$ and the same holds for $\dim\, \nested{2,n}_0$ by \Cref{thm:motivepunctual}. However, the polynomial $[\nested{2,n}_0]$ is not monic. This is a motivic evidence for the reducibility result   presented in \cite[Corollary 7.5]{BULOIS}, see also \Cref{rem:BULOIS} and  \Cref{cor:irrcomp2n,cor:irrcomp3n}. 
\end{remark}

As an immediate  application of \Cref{thm:motivepunctual} we obtain the following.

\begin{corollary} 
\label{cor:irrcomp2n} 
For $n\geqslant 4$ the projective scheme $S_p^{[2,n]}$ has dimension $n-1$. Moreover, if we denote by $\nu^{[2,n]}$ the number of irreducible components of $S_p^{[2,n]}$ of maximal dimension, we have
\[
\nu^{[2,n]}=\left\lfloor\frac{n}{2}\right\rfloor.
\] 
\end{corollary}
\begin{proof}
Combine \Cref{lemma:intercoeff} and \Cref{thm:motivepunctual} with one another.
\end{proof}

\begin{remark}\label{rem:BULOIS}
It is worth mentioning that in \cite[Cor.~7.5]{BULOIS} a stronger result than \Cref{cor:irrcomp2n} is presented. Precisely, the authors show equidimensionality of $S_p^{[2,n]}$ and hence they deduce that the number $\lfloor n/2\rfloor$ is comprehensive of all the components. In \Cref{cor:irrcomp3n} we prove an analogue of \Cref{cor:irrcomp2n} for $S_p^{[3,n]}$ which, to the best of our knowledge, is not present in the literature.
\end{remark}

To prove the first identity in \Cref{thm:intro-k=2and3}, it remains to assemble the motives $[S_p^{[2,\bullet]}]$ computed in \Cref{thm:motivepunctual} into a generating function. Recall that, by work of G\"{o}ttsche \cite{Gottsche-motivic}, the motivic generating function
\[
\mathsf{Hilb}^\bullet(t) = \sum_{n\geqslant 0} \,[S_p^{[n]}]t^n = 1+t+(\BL+1)t^2+\cdots
\]
can be written in the explicit form
\[
\mathsf{Hilb}^\bullet(t) = \prod_{j\geqslant 1}\,\frac{1}{1-\BL^{j-1}t^j}.
\]
We next express the punctual nested series
\[
\mathsf{Hilb}^{[2,\bullet]}(t) = \sum_{n\geqslant 2}\,[S_p^{[2,n]}]t^{n}
\]
in terms of $\mathsf{Hilb}^\bullet(t)$, proving the first formula in \Cref{thm:intro-k=2and3}.
 
\begin{corollary}
\label{cor:k=2}
There is an identity
\[
\mathsf{Hilb}^{[2,\bullet]}(t)=[\BP^1]\mathsf{Hilb}^\bullet(t)+[\BP^1]\frac{t(\BL-1)-1}{1-\BL t}
\]
in $K_0(\Var_{\BC})\llbracket t \rrbracket$.
\end{corollary}

\begin{proof}
We have
\begin{align*}
\mathsf{Hilb}^{[2,\bullet]}(t) 
&= [\BP^1]\sum_{n\geqslant 2}\,\left([\nested{n}_0]-\BL^{n-1}\right)t^{n} \\
&=[\BP^1]\left(\sum_{n\geqslant 2}\,[\nested{n}_0]t^{n}-\sum_{n\geqslant 0}\,\BL^{n+1}t^{n+2}\right) \\ 
&=[\BP^1]\left(\mathsf{Hilb}^\bullet(t)-1-t-\frac{\BL t^2}{1-\BL t}\right) \\
&=[\BP^1]\left(\mathsf{Hilb}^\bullet(t)+\frac{t(\BL-1)-1}{1-\BL t}\right),
\end{align*}
as required.
\end{proof}

The series starts as
\begin{multline*}
    \mathsf{Hilb}^{[2,\bullet]}(t) = (1 + \BL)t^2 + (1 + \BL)^2 t^3 + (1 + \BL) (1 + \BL + 2 \BL^2) t^4 + (1 + \BL)^2 (1 + 
    2 \BL^2) t^5 \\
    + (1 + \BL) (1 + \BL + \BL^2 (2 + 3 \BL (1 + \BL))) t^6 + (1 + 
    \BL)^2  (1 + \BL^2 (2 + \BL + 3 \BL^2))  t^7 + \cdots
\end{multline*}

\section{The motive of \texorpdfstring{$S_p^{[3,n]}$}{}}\label{sec:nested_3n}
Similarly to Section \ref{sec:motive_2n}, we want to study the motive of the punctual nested Hilbert scheme $S_p^{[3,n]}=\nested{3,n}_0$. The (smooth) base case $n=4$ is rather special and already known in the literature, see \cite{Gottsche-motivic}. Let us then fix $n\geqslant 5$. The main result of this subsection is the second part of the statement of \Cref{thm:intro-k=2and3}, which we recall below.

\begin{theorem}
\label{motive-punctual-3n}  
Let $S$ be a smooth quasiprojective surface, $p \in S$ a closed point. For every $n \in \BZ_{\geqslant 5}$, there is an identity 
 { \begin{equation}
\label{eq:totale3n}
    [S_p^{[3,n]}]=[\BP^2][S_p^{[n]}] 
    - \BL^{n-3}[\BP^1]\left( \BL^2 + \left\lfloor \frac{n}{2}\right\rfloor\BL  -1  \right) 
\end{equation}}
in $K_0(\Var_{\BC})$.
\end{theorem} 
As a consequence of \Cref{motive-punctual-3n} we immediately get the following.

\begin{corollary}
\label{cor:irrcomp3n}
For $n\geqslant 4$ the projective scheme $S_p^{[3,n]}$ has dimension $n-1$. Moreover, if $n\ge 5$ and we denote by $\nu^{[3,n]}$ the number of irreducible components of $S_p^{[3,n]}$ of maximal dimension, we have
\[
\nu^{[3,n]}=    {\left\lfloor{\frac{n^2+4}{12}}\right\rfloor  }.
\]  
\end{corollary}
\begin{proof}
Combine \Cref{lemma:intercoeff} and \Cref{motive-punctual-3n} with one another.
\end{proof}

As before, our strategy to prove \Cref{motive-punctual-3n} will consist in finding a suitable stratification of $S_p^{[3,n]} = \nested{3,n}_0$. We start by writing
\[
\nested{3,n}_0= \mathsf C^{[3,n]} \amalg Y^{[3,n]}_{(1,2)}\amalg Y^{[3,n]}_{(1,1,1)},
\]
where 
\[
\mathsf C^{[3,n]}=\Set{[Z_1\subset Z_2]\in \nested{3,n}_0| \bh_{Z_2}(1)=1}
\]
is the nested curvilinear locus as defined in \Cref{sec:puttane-curvilinee}, and
\[
Y_\bullet^{[3,n]}=\Set{[Z_1\subset Z_2]\in \nested{3,n}_0\setminus \mathsf C^{[3,n]}| \bh_{Z_1}=\bullet}.
\]
The motives $[\mathsf C^{[3,n]}]$ and $[Y_{(1,2)}^{[3,n]}]$ are easily computed as there are bijiective morphisms
\[
\begin{tikzcd}[row sep = tiny]
    \mathsf C^{[3,n]} \arrow{r}{\pr_2} &  \mathsf C^{[n]}\\
    {[Z_1\subset Z_2]} \arrow[r,mapsto] &  {[Z_2]}
\end{tikzcd}\quad\mbox{ and }\quad\begin{tikzcd}[row sep = tiny]
    Y^{[3,n]}_{(1,2)} \arrow{r}{\pr_2} &  \nested{n}_0\setminus \mathsf C^{[n]}\\
    {[Z_1\subset Z_2]} \arrow[r,mapsto] &  {[Z_2]}.
\end{tikzcd}
\]
Hence, we have
\begin{equation}\label{eq:1pezzo3n}
    [\mathsf C^{[3,n]} \amalg Y^{[3,n]}_{(1,2)}]=[S_p^{[n]}].
\end{equation}
To compute the motive of $Y_{(1,1,1)}^{[3,n]}$, we decompose the space into three locally closed strata, namely
\[
Y_{(1,1,1)}^{[3,n]}=Y_{(1,1,1),1}^{[3,n]}\amalg Y_{(1,1,1),2}^{[3,n]}\amalg Y_{(1,1,1),3}^{[3,n]},
\]
where
\[
Y_{(1,1,1),i}^{[3,n]}=\Set{ [Z_1\subset Z_2 ]\in Y_{(1,1,1)}^{[3,n]} | \bh_{Z_2}(2)=i }.
\]
Again, the motives of two of these strata can be easily computed. Indeed, the second projection induces Zariski locally trivial fibrations
\[
\begin{tikzcd}[row sep = tiny]
    Y_{(1,1,1),1}^{[3,n]} \arrow[r,"\BA^1"] & H_1^{[n]}\\
    {[Z_1\subset Z_2]} \arrow[r,mapsto] &  {[Z_2]}
\end{tikzcd}\quad\mbox{ and }\quad\begin{tikzcd}[row sep = tiny]
   Y_{(1,1,1),3}^{[3,n]}\arrow[r,"\mathsf C^{[3]}"] &  H_3^{[n]} \\
    {[Z_1\subset Z_2]} \arrow[r,mapsto] &  {[Z_2]}
\end{tikzcd}
\]
with fibres $\BA^1$ and $\mathsf C^{[3]}$, respectively. This gives
\begin{equation}
\label{eq:2pezzo3n}
    [Y_{(1,1,1),1}^{[3,n]}] =[H_1^{[n]}]\BL\quad\mbox{ and }\quad [Y_{(1,1,1),3}^{[3,n]} ] =[H_3^{[n]}]   [\mathsf C^{[3]}].
\end{equation}
The motives of $H_i^{[n]}$ for $i=1,3$ were computed in \Cref{eq:H1} and \Cref{prop:H3} respectively. The last contribution we have to compute is the motive  of $Y_{(1,1,1),2}^{[3,n]}$. This is the content of the following lemma.
\begin{lemma}\label{lemma:Y1112}
For each $n \in \BZ_{\geq 5}$ there is an identity
   
    \[
    { 
    [Y_{(1,1,1),2}^{[3,n]}]=\BL^{n-4}\left( ({n-5} ) \BL^3 + \left(\left\lfloor \frac{n}{2}\right\rfloor+n-6  \right) \BL^2+ \left\lfloor \frac{n}{2}\right\rfloor \BL +1\right)}
    \] 
in $K_0(\Var_{\BC})$. 
\end{lemma}

\begin{proof}
We start by stratifying the locus $Y_{(1,1,1),2}^{[3,n]}$ as
    \[
    Y_{(1,1,1),2}^{[3,n]}=\coprod_{s=2}^{\left\lfloor\frac{n-1}{2}\right\rfloor}A_{1,s}\amalg\coprod_{s=2}^{\left\lfloor\frac{n-1}{2}\right\rfloor}A_{2,s},
    \]
    where
    \[
    A_{i,s}=\Set{ [Z_1\subset Z_2] \in Y_{(1,1,1),2}^{[3,n]} | [Z_2] \in H_{\bh_s}^i}.
    \]
The locus $H_{\bh_s}^i$ was defined in \Cref{def:H^i_s}. Now, for $i=1$ we consider the $\BA^1$-fibration
    \[
    \begin{tikzcd}[row sep = tiny]
        A_{1,s}\arrow[r,"\pi_s"]&  H_{\bh_s}^1\\
        {[I_1\supset I_2]}\arrow[r,mapsto]& {[I_2]}.
    \end{tikzcd}
    \]
The motive of $[H_{\bh_s}^1]$ has been computed in \Cref{prop:stratquad}. To see that the fibre of $\pi_s$ is an affine line we show first that the {  tangent space to each fibre of $\pi_s$ is one dimensional and then we exhibit an injective morphism from $\BA^1$ to each fibre. This ensures surjectivity as well and hence an equality in the Grothendieck ring of varieties, see \Cref{sec:K0(Var)}. Note that each fibre contains a unique homogeneous point: this trivialises the map $\pi_s$. 
In particular, if $[I_1\supset I_2]\in A_{1,s}$ has $I_1$ homogeneous, then we have 
\[
\mathsf T_{[I_1\supset I_2 ]} \pi_s^{-1}( [I_2]) \cong \mathsf T_{[I_1 ]}^{>0} (\BA^2 )^{[3]}   ,
\]
which is an immediate computation via \Cref{rem:nonneg-punctual}, see \Cref{rmk:novo}. Now, by \Cref{thm:briarro}, we have $\dim_{\BC} \mathsf T_{[I_1 ]}^{>0} (\BA^2 )^{[3]} =1$ which concludes the computation of the dimension of the tangent space. We now single out a copy of $\BA^1$  in each fibre of $\pi_s$.

}  
 Without loss of generality, we can suppose that $y^2$ is the unique degree 2 generator of $\In I_2$, which implies that the unique point $[I_1\supset I_2]$ with $I_1$ homogeneous in the fibre of $\pi_s^{-1}[I_2]$ is the point $[(x^3,y)\supset I_2]$. Notice that we have the equality of ideals
    \begin{equation}\label{eq:poscurv}
        (x^3,y+\alpha x^2)=(y+\alpha x^2,y^2,xy)+\Fm^3\subset R[\alpha],
    \end{equation}
      Since all the elements in $I_2$ have initial form of degree greater or equal than 3 or initial form  $y^2$, the ideal \eqref{eq:poscurv} contains $I_2$ for every choice of $\alpha \in \BC$.

    Let us now consider the strata $A_{2,s}$, for  $s=2,\ldots,\left\lfloor\frac{n-1}{2}\right\rfloor$ and let us denote by $\OA_{2,s}$ the locus parametrising pairs of homogeneous ideals
    \[
    \OA_{2,s}=\Set{[I_1\supset I_2] \in A_{2,s} | \In I_i=I_i,\mbox{ for }i=1,2}.
    \]
    Then, we have again a composition of Zariski locally trivial fibrations 
    \[
    \begin{tikzcd}[row sep = tiny]
        A_{2,s}\arrow[r,"\psi_s" ]&  \OA_{2,s}\arrow[r,"\varphi_s"]&  \OC^{[3]}\cong\BP^1\\
        {[I_1 \supset I_2]}\arrow[r,mapsto]& {[\In I_1\supset \In I_2]}\arrow[r,mapsto]& {[\In I_1]},
    \end{tikzcd}
    \]
    where $\OC^{[3]}\subset \mathsf C^{[3]}$ denotes the homogeneous locus. The fibres $G_s$ of $\psi_s$ and $F_s$ of $\varphi_s$ { have motives}
  {   \[
    [G_s]=\begin{cases} 
        \BL^{n-4}&\mbox{ if }  s=\frac{n-1}{2} ,\\
        \BL^{n-4}&\mbox{ if }  s=\frac{n-2}{2} ,\\
        \BL^{n-3}&\mbox{ if } 2\leqslant s < \left\lfloor\frac{n-1}{2}\right\rfloor,
    \end{cases}
    \qquad
    [F_s]=\begin{cases}
            \BL&\mbox{ if }   s = \frac{n-1}{2},\\
            [\BP^1]\BL&\mbox{ if }   s = \frac{n-2}{2},\\
        \BL +  \BL&\mbox{ if } 2\leqslant s < \left\lfloor\frac{n-1}{2}\right\rfloor.
    \end{cases}
    \] }
The first fibre is again deduced from \eqref{eq:poscurv} and  \Cref{prop:stratquad}.    

Let us explain the second fibre. Fix some point $[I_1\supset I_2]\in \OA_{2,s}$. Without loss of generality, we can suppose 
        \[
     I_1 = \langle (x^*)^2\rangle ^\perp=(y,x^3) , 
    \]
    see  \Cref{prop:gorapolar}. Then, the ideal $I_2$ satisfies
    \[
    I_2=\langle  (x^*)^{s},\ell^{n-s-1}\rangle^\perp\quad\mbox{ or }\quad I_2= \langle \ell^{s}, (x^*)^{n-s-1}\rangle^\perp,
    \]
    for some $\ell \in  \BP R_1^*\setminus [x^*]$. Now, when {  $2\leqslant s < \left\lfloor\frac{n-1}{2}\right\rfloor$}, the two choices for $I_2$ are not equivalent and the loci parametrising them are clearly disjoint. {   Note that in the case $s=\frac{n-2}{2}$,  one must take into account that the choice of the generator of maximal degree takes place in a $\BP^1$.} 
    Summarising, we have
    \[
    [A_{1,s}]=[H_{\bh_s}^1]\BL,
    \]
    and, by \Cref{prop:stratquad},
    \[
    [A_{2,s}]=\begin{cases}
       [\BP^1]\BL^{n-3} & \text{ if  }s=\frac{n-1}{2}\\
        [\BP^1]^2\BL^{n-3} & \text{ if  }s=\frac{n-2}{2}\\
       2[\BP^1]\BL^{n-2} & \text{ if  } 2 \le s<\floor{\frac{n-1}{2}}\\
    \end{cases}
    \]
    To get the statement it is enough to sum over all $2\le s \le \floor{\frac{n-1}{2}}$.
\end{proof}

By summing all the contributions computed so far, we get the proof of \Cref{motive-punctual-3n}.

\begin{proofof}{\Cref{motive-punctual-3n}} In order to prove the claim, we sum all the motivic contributions
\[
[\nested{3,n}_0]= [\mathsf C^{[3,n]}] \amalg [Y^{[3,n]}_{(1,2)}]\amalg [Y^{[3,n]}_{(1,1,1)}],
\]
from \eqref{eq:1pezzo3n},\eqref{eq:2pezzo3n}, \eqref{eq:H1}, \Cref{prop:H3,lemma:Y1112}.
\end{proofof}

\begin{corollary}
\label{cor:k=3}
There is an identity
\[\mathsf{Hilb}^{[3,\bullet]}(t)= [\BP^2]\mathsf{Hilb}^\bullet(t) 
-\frac{[\BP^2] - (\BL^3-1) t - (\BL^3-1)[\BP^1] t^2 - (\BL^2-\BL^5) t^3 - \BL^2 t^4  }{(1-\BL t)(1-\BL^2 t^2)}.  \]
\end{corollary}

\begin{proof}
The statement follows expanding the formula in \Cref{motive-punctual-3n} after adding in the terms corresponding to the nestings $[3,3]$ and $[3,4]$, namely the motives
\begin{align*}
[S_p^{[3,3]}] &= [S_p^{[3]}]=1+\BL+\BL^2,\\
[S_p^{[3,4]}] &= 1 + \BL (1 + \BL) (2 + \BL),
\end{align*}
the second one being computed in \cite{Gottsche-motivic}.
\end{proof}

The proof of \Cref{thm:intro-k=2and3} is now complete.

\part*{Euler characteristics}
We prove in this part of the paper Theorems \ref{thm:intro:Z_k}, \ref{thm:intro-Z_Dr} and \ref{thm:Z_Dr3-intro} from the Introduction. 

We introduce in  \Cref{subsec:nestedpart} some terminology and notation regarding (skew) Ferrers diagrams. We then move, step by step, towards the proof of \Cref{thm:intro:Z_k}, which says that, for every $\bk = (k_1,\ldots,k_s) \in \BZ_{\geq 0}^s$ and for every pointed smooth surface $(S,p)$, the series
\begin{equation}
\label{eqn:Z_k-general}
\FFZ_{\bk}(q)=\sum_{n\geqslant 0} \chi\left(S_p^{[n,n+k_1,n+k_1+k_2,\ldots,n+\sum_{i=1}^sk_i]}\right)q^n    
\end{equation} 
is a product of $\chi \mathsf{Hilb}^\bullet(q) = \prod_{m\geq 1}(1-q^m)^{-1}$ and a rational function with only roots of unity as poles. We will start in \Cref{subsec:bassecase} by solving the case $s=1$ of \Cref{thm:intro:Z_k}, that will be extended to the general case $s>1$ in \Cref{subsec:g.f.tantinest}. The higher rank version of this result (where the main character is the nested Quot scheme) will be given in \Cref{sec:C-D}, whereas the proofs of Theorems \ref{thm:intro-Z_Dr} and \ref{thm:Z_Dr3-intro} will be presented in \Cref{sec:C-D}.

\section{Nested partitions}\label{subsec:nestedpart}

In this section we review the definition of \emph{skew Ferrers diagram} and we set up the notation needed in later sections. We also recall in \Cref{sec:euler-char-combinatorial} the relation between $\bn$-nested Ferrers diagrams and the Euler characteristic of the punctual nested Hilbert scheme $S_p^{[\bn]}$ of a pointed smooth surface $(S,p)$.

\subsection{Nested Ferrers diagrams}
\label{sec:nested-Ferrers}
The monoid $(\BN^2,+)$ is endowed with its standard component-wise poset structure throughout. 

\begin{definition}
\label{def:partition}
Fix $n \in \BZ_{\geqslant 0}$. A \textit{Ferrers diagram} of size $n$ is a collection of $n$ points $\lambda =\set{ a_1,\ldots, a_n}\subset \BN^{2}$ such that if  $ a \in \BN^{2}$ satisfies $ a\le a_i$ for some $i=1,\ldots,n$, then $ a \in \lambda$. We denote by $\mathrm{P}^{[n]}$ the set of Ferrers diagrams of size $n$, and by $\chi^{[n]}$ the cardinality $\chi^{[n]}=\lvert \mathrm{P}^{[n]}\rvert$. A \emph{skew Ferrers diagram} is the set-theoretic difference of two Ferrers diagrams. We denote by $\mathrm{Q}^{[n]}$ the set of skew Ferrers diagrams of size $n$.
\end{definition} 

Ferrers diagrams are also known as \emph{Young diagrams}, and they are in natural bijection with classical  \emph{partitions}. 

\begin{notation}
\label{notation:linear-partitions}
When we identify a Ferrers diagram $\lambda \in \mathrm{P}^{[n]}$ with a partition, we may choose to represent it in any of the two following (equivalent) ways:
\begin{itemize}
    \item [\mylabel{linear-partition-1}{(1)}] $\lambda = (\lambda_1,\lambda_2,\ldots,\lambda_{l_\lambda})$ where $\lambda_1\geqslant \lambda_2\geqslant \cdots\geqslant \lambda_{l_\lambda}>0$ and the size of the partition is $|\lambda| = \sum_j\lambda_j = n$,
    \item [\mylabel{linear-partition-2}{(2)}] $\lambda = (1^{m_1(\lambda)}2^{m_2(\lambda)}\cdots n^{m_n(\lambda)})$ where the size of the partition is $|\lambda| = \sum_iim_i(\lambda) = n$. The exponents $m_i(\lambda)$ are also called multiplicities of the partition $\lambda$ --- some $m_i$ might be zero, in which case the term $i^{m_i}$ will be omitted. In particular, this notation implies that $m_j(\lambda) = 0$ for $j \geqslant n+1$. We write simply $m_i$ when it is clear which partition we are referring to.
\end{itemize}
\end{notation}

The notation in \ref{linear-partition-2} means that $\lambda$ consists of $m_i(\lambda)$ \emph{parts} of size $i$ for all $i=1,\ldots,n$, so that $l_\lambda = \sum_i m_i(\lambda)$, the total number of parts, counts the number of rows in the Ferrers diagram attached to $\lambda$. It is straightforward to pass from \ref{linear-partition-1} to \ref{linear-partition-2} and viceversa. 

\begin{example}
The Ferrers diagram
\[
\begin{matrix}
       \yng(4,3,3,1,1) \quad 
       \ytableausetup{boxsize=1.09em}
       \ytableausetup{boxframe=0.02em}\ytableausetup{aligntableaux=bottom}
   \end{matrix}
\]
 corresponds to $\lambda = (4,3,3,1,1) = (1^2 3^2 4^1)$, has size $12$, and has $l_\lambda = 2+2+1 = 5$ parts (rows). Notice that we employ the English notation for partitions, i.e. the orientation of the $x$-axis goes from left to right and the orientation of the $y$-axis goes from top to bottom.
\end{example}

Note that a \emph{skew} Ferrers diagram may be disconnected, which means that there is at least one pair of boxes which cannot be connected through a sequence of boxes sharing an edge with one another.
For later convenience we shall identify skew Ferrers diagrams that differ by a translation.

\begin{definition}
    Let $\lambda^{(1)},\lambda^{(2)}\in \mathrm Q^{[n]}$ be two skew Ferrers diagrams with the same number $k$ of connected components, and let
    \[
\lambda^{(j)}=\coprod_{i=1}^{k}\lambda_i^{(j)}
    \]
be their decomposition in their respective connected components. We declare that $\lambda^{(1)}$ and $\lambda^{(2)}$ are \textit{equivalent up to translation} if there is a permutation $\sigma\in\FS_k$, acting as a relabelling of the indices, and for each $i=1,\dots,k$ a vector $\bv_i\in\BZ^2$, such that $\lambda_i^{(1)}=\bv_i+\lambda^{(2)}_{\sigma(i)}$.
\end{definition}

We denote by $\mathscr Q^{[n]}$ the set of equivalence classes 
    \[
    \mathscr Q^{[n]}=\mathrm Q^{[n]}\big/\,\mathrm{translation}.
    \]
Recall that for $\lambda\in\mathrm P^{[n]}$, the transposed $\lambda^T\in\mathrm P^{[n]}$ is defined to be the reflection of $\lambda$ along the main diagonal, i.e.~the line $\set{(x,y)\in\BN^2|y=x}\subset \BN^2$.

\begin{notation}
\label{notation:everything}
We shall write $\mathrm{P} = \coprod_{n\geq 0}\mathrm{P}^{[n]}$ and $\mathscr Q = \coprod_{n \geq 0} \mathscr Q^{[n]}$ in what follows.
\end{notation}

\begin{definition}\label{def:transposed_diagram}
For a skew Ferrers diagram $[\lambda]\in\mathscr Q^{[n]}$ such that $\lambda=\mu\smallsetminus\nu$, for some $\mu,\nu\in\mathrm P$, define the transposed $[\lambda]^T\in\mathscr Q^{[n]}$ as $[\lambda]^T=[\mu^T\smallsetminus\nu^T]$.
\end{definition}

For the sake of conciseness, we will not make a notational distinction between skew Ferrers diagrams $\lambda\in\mathrm Q^{[n]}$ and equivalence classes $[\lambda]\in\mathscr Q^{[n]}$.

\begin{definition}
\label{def:nesteDpartition}
Fix a nondecreasing sequence $\bn=(n_1,\ldots,n_\ell) \in \BZ_{\geqslant 0}^\ell$. An $\bn$-\emph{nesting of Ferrers diagrams} is a filtration of Ferrers diagrams $\lambda_1\subset \cdots \subset \lambda_\ell\subset \BN^2$ such that $\lambda_i\in \mathrm{P}^{[n_i]}$, for all $i=1,\ldots,\ell$. We denote by $\mathrm{P}^{[\bn]}$ the set of $\bn$-nestings of Ferrers diagrams, and by $\chi^{[\bn]}$ the cardinality $\chi^{[\bn]}=\lvert \mathrm{P}^{[\bn]}\rvert$.
\end{definition}

\subsection{The Euler characteristic of \texorpdfstring{$S_p^{[\bn]}$}{}}
\label{sec:euler-char-combinatorial}

Let $(S,p)$ be any pointed smooth quasiprojective surface. It is a standard fact that the number $\chi^{[n]}$ of partitions of size $n$ agrees with the topological Euler characteristic of $S_p^{[n]}$. This is because the natural scaling action of the torus $\mathbf{T} = \BG_{m}^2$ on $\BA^2$ has a natural lift to $\nested{n}_0 \cong S_p^{[n]}$, and the $\mathbf{T}$-fixed locus consists of monomial ideals of colength $n$. These, in turn, correspond bijectively to partitions of $n$. In exactly the same way, one obtains the identity
\[
\chi^{[\bn]} = \chi\left(S_p^{[\bn]}\right).
\]
For future use, let us define the series
\begin{equation}
\label{eqn:Z-series}
\begin{split}
\ZZ(q)&=\sum_{n \geqslant 0} \chi^{[n]} q^n \\
\FFZ_D(q)&=\sum_{n \geqslant 0} \chi^{[n,n+D]} q^n
\end{split}
\end{equation}
in $\BZ \llbracket q\rrbracket$, the second series being defined for an arbitrary integer $D \geq 0$. The decoration `$\mathsf F$' in the notation, here and throughout, stands for `flag'. Clearly one has $\FFZ_0 = \ZZ$.

\begin{example}
\label{ex:Z_1}
Euler proved the relation
\[
\ZZ(q)=
\prod_{j\geqslant 1}\frac{1}{1-q^j}.
\]
The first nontrivial example of the series \eqref{eqn:Z-series} is $\FFZ_1$, which can be extracted from G\"{o}ttsche's motivic formula \cite[Sec.~5]{Gottsche-motivic} or from Cheah's work \cite[Thm.~3.3.3]{MR1616606}. It reads
\begin{equation}
\label{eqn:Z1}
\FFZ_1(q) = \frac{1}{1-q}\ZZ(q).
\end{equation}
\end{example}

\begin{remark}
\label{rmk:alternative}
There is, in fact, a direct combinatorial argument to prove \Cref{eqn:Z1}. The coefficient of $q^n$ in $\FFZ_1(q)$, namely the cardinality $\chi^{[n,n+1]}$ of $\mathrm{P}^{[n,n+1]}$, is the number of pairs $(\lambda,\lambda') \in \mathrm P^{[n]} \times \mathrm P^{[n+1]}$ such that $\lambda \subset \lambda'$. On the other hand, the coefficient of $q^n$ in $\ZZ(q)/(1-q)$ is 
\[
\sum_{k=0}^n \chi^{[k]},
\]
which is nothing but the number of Ferrers diagrams of size at most $n$. These two coefficients agree for the following reason. An element of $\mathrm{P}^{[n,n+1]}$ can be seen as a Ferrers diagram $\lambda' \in \mathrm P^{[n+1]}$ with one box marked, located in a position such that removing it one gets a Ferrers diagram $\lambda \in \mathrm{P}^{[n]}$. The required bijection is given by sending $\lambda'$ to the partition where the whole row with the marking is removed. The inverse is given by adding to a given $\nu \in \mathrm{P}^{[k]}$ a new row of $n + 1 - k$ boxes (in the unique position such that the result is a Ferrers diagram) with a marking in the last box of that row.
\end{remark}

\section{The Euler characteristic of the nested Hilbert scheme}
\label{sec:g.f.}

After providing an independent proof of \Cref{eqn:Z1} in \Cref{sec:D=1}, we introduce in \Cref{sec:madonna-scannata} a differential operator $T_\lambda$ attached to a skew Ferrers diagram $\lambda$. Together with the $D$-decorated partition function $F^{[D]}$ introduced in \Cref{def:F-partition-function}, this is the key ingredient for the proof of \Cref{thm:intro:Z_k}, via \Cref{cor:TD-applied-to-FD}, which expresses the series $\FFZ_D$ in terms of $F^{[D]}$. In \Cref{rmk:transpose} we prove a useful transposition invariance property for the series $T_\lambda . F^{[D]}\big|_{\vec{y} = 1}$, whose enumerative meaning is pinpointed in \Cref{thm:skewferret}.

\subsection{The series \texorpdfstring{ $\FFZ_1(q)$}{} revisited}
\label{sec:D=1}
We present here an independent proof of the formula \eqref{eqn:Z1} for $\FFZ_1(q)$, cf.~\Cref{Z1-again}. This serves as an illustration of the general technique that we shall develop in order to prove \Cref{thm:intro:Z_k}. Define
\begin{equation}
\label{F-series}
	F(x;q) = x \cdot \prod_{i\geqslant 1} \left(1 + x \sum_{m_i \geqslant 1} (q^{i})^{m_i} \right).
\end{equation}
It recovers the generating series $\ZZ(q)$ at $x=1$, i.e.~one has
\[
F(1;q) = \ZZ(q).
\]

\begin{lemma} 
\label{lemma:Z1-and-derivative}
The series $\FFZ_1(q)$ can be expressed in terms of $F(x;q)$ via the identity
\begin{equation}
    \FFZ_1(q) = \frac{\mathrm{d}}{\mathrm{d}x} F(x;q) \Bigg{|}_{x=1}.
\end{equation}

\begin{proof}
Let us write a partition $\lambda$ in terms of its multiplicities (cf.~\Cref{notation:linear-partitions}\ref{linear-partition-2}) as 
$$
\lambda = \left(1^{m_1(\lambda)} 2^{m_2(\lambda)} 3^{m_3(\lambda)} \cdots \right).
$$
Notice that these multiplicities arise from the indices of expansion for 
$$
\ZZ(q) = \prod_{j\geqslant 1} \frac{1}{1-q^j} = \prod_{j\geqslant 1} \sum_{m_j \geqslant 0} \left(q^{j}\right)^{m_j}
$$
when expanding $\ZZ(q)$ as a $q$-series. We now describe the addition of one box to the Ferrers diagram of a partition in term of its multiplicities: there is a bijection between the number of single boxes that can be added to a Ferrers diagram turning it into another Ferrers diagram and the number of its nonzero multiplicities, plus one. This is true because there is an admissible position for a box addition precisely where the parts of the diagram change size, except for the very first admissible box addition which is granted to all Ferrers diagrams.

\begin{figure}[H]
\begin{tikzpicture}
\draw[<->] (1.8,0)--(2.8,0);
\draw[<->] (6.8,0)--(7.8,0);
    \node at (0,0) {
\begin{tikzpicture}
\node at (2.5/3,4/3) {\Large$\lambda$};
    \draw (0,0)--(2/3,0)--(2/3,1/3)--(3/3,1/3)--(3/3,2/3)--(6/3,2/3)--(6/3,4/3)--(9/3,4/3)--(9/3,7/3)--(0,7/3)--cycle;
\end{tikzpicture}};
    \node at (4.8,0) {
\begin{tikzpicture}
    \draw (0,0)--(2/3,0)--(2/3,1/3)--(0,1/3)--cycle;
    \draw(0,1/3+0.1)--(3/3,1/3+0.1)--(3/3,2/3+0.1)--(0,2/3+0.1)--cycle;;
    \draw (0,2/3+0.2)--(6/3,2/3+0.2)--(6/3,4/3+0.2)--(0,4/3+0.2)--cycle;
    \draw (0,4/3+0.3)--(9/3,4/3+0.3)--(9/3,7/3+0.3)--(0,7/3+0.3)--cycle;
    \node[left] at (0.1,- 0.05) {\tiny $5\to$};
    \node[left] at (0.1,1/3+ 0.05) {\tiny $4\to$};
    \node[left] at (0.1,2/3+ 0.15) {\tiny $3\to$};
    \node[left] at (0.1,4/3+ 0.25) {\tiny $2\to$};
    \node[left] at (0.1,7/3+ 0.35) {\tiny $1\to$};
\end{tikzpicture}};
    \node at (10,0) {
\begin{tikzpicture}
    \draw (0,0)--(2/3,0)--(2/3,1/3)--(3/3,1/3)--(3/3,2/3)--(6/3,2/3)--(6/3,4/3)--(9/3,4/3)--(9/3,7/3)--(0,7/3)--cycle;
    \draw[pattern=north east lines] (0,-1/3)--(1/3,-1/3)--(1/3,0)--(0,0)--cycle;
    \draw[pattern=north east lines] (0+2/3,-1/3+1/3)--(1/3+2/3,-1/3+1/3)--(1/3+2/3,0+1/3)--(0+2/3,0+1/3)--cycle; 
    \draw[pattern=north east lines] (0+3/3,-1/3+2/3)--(1/3+3/3,-1/3+2/3)--(1/3+3/3,0+2/3)--(0+3/3,0+2/3)--cycle; 
    \draw[pattern=north east lines] (0+6/3,-1/3+4/3)--(1/3+6/3,-1/3+4/3)--(1/3+6/3,0+4/3)--(0+6/3,0+4/3)--cycle; 
    \draw[pattern=north east lines] (0+9/3,-1/3+7/3)--(1/3+9/3,-1/3+7/3)--(1/3+9/3,0+7/3)--(0+9/3,0+7/3)--cycle; 
    \node[right] at  ( 1/3 ,-1/3  ) {\tiny 5};
    \node[right] at  ( 1/3+2/3,-1/3+1/3 ) {\tiny 4};
    \node[right] at  ( 1/3+3/3,-1/3+2/3 ) {\tiny 3};
    \node[right] at  ( 1/3+6/3,-1/3+4/3 ) {\tiny 2};
    \node[right] at  ( 1/3+9/3,-1/3+7/3 ) {\tiny 1};
\end{tikzpicture}};
\end{tikzpicture}
  \caption{On the left the Ferrers diagram of a partition $\lambda$. In the center its decomposition in blocks corresponding to each of nonzero multiplicities $m_i(\lambda)$, and marked all possible delimiters between blocks, including the two extremes top and bottom. On the right all possible additions of a single box in $\lambda$: they correspond one-to-one to blocks delimiters, i.e. to nonzero multiplicities $m_i(\lambda)$, plus one.}
  \label{fig:decomposition_D=1}
\end{figure}

Hence we simply need to keep track of how many nonzero multiplicities there are in every given partition, and we do that at the level of the generating series $\mathsf Z(q)$ by decorating it with an additional variable $x$ that multiplies all nonzero multiplicities, times an extra power of $x$. This $x$-decoration of $\mathsf Z(q)$ results in the definition of 
$$
F(x;q) = x \cdot \prod_{i\geqslant 1} \left(1 + x \sum_{m_i \geqslant 1} (q^{i})^{m_i} \right).
$$
Expanding $F$ in the $x$-variable as 
$$
F(x;q) = \sum_{k\geqslant 0} f_{k}(q) x^{k+1},
$$
the coefficients $f_k(q)$ are the generating series of all those partitions with exactly $k$ nonzero multiplicities, and hence $k+1$ admissible box additions, which is the partition weight we are after. In order to weigh each of these partitions by $k+1$ it suffices to take the first derivative in $x$, yielding
$$
\frac{\mathrm{d}}{\mathrm{d}x} F(x;q) = \sum_{k\geqslant 0} (k+1)f_{k}(q) x^{k}.
$$
We can then safely set $x=1$ to obtain the wanted $q$-generating series. 
\end{proof}
\end{lemma}

We now employ the previous lemma to compute $\FFZ_1(q)$ explicitly. 

\begin{prop} 
\label{Z1-again}
There is an identity
\begin{equation}
\FFZ_1(q) = \frac{1}{1-q}\ZZ(q).
\end{equation}
\begin{proof}
It follows from \Cref{lemma:Z1-and-derivative} and the computation of the derivative:
\begin{align}
	\frac{\mathrm{d}}{\mathrm{d}x} F(x;q) \Bigg{|}_{x=1} &= \frac{\mathrm{d}}{\mathrm{d}x} \left[x \cdot \prod_{i\ge 1}  \left(1 + x \sum_{m_i \geqslant 1} (q^{i})^{m_i} \right) \right] \Bigg{|}_{x=1}
	\nonumber\\
	&= \prod_{i\ge 1} \left(1 + x \sum_{m_i \geqslant 1} (q^{i})^{m_i} \right) + x\cdot  \sum_{j\ge 1} \left(\sum_{m_j \geqslant 1} (q^{j})^{m_j} \right) \prod_{ \substack{i\ge1 \\i \neq j}}\left(1 + x \sum_{m_i \geqslant 1} (q^{i})^{m_i} \right) \Bigg{|}_{x=1}
	\nonumber\\
	&= \ZZ(q) + \sum_{j\ge1} \frac{q^j}{1-q^j} \prod_{ \substack{i\ge1 \\i \neq j}} \left(1 + \sum_{m_i \geqslant 1} (q^{i})^{m_i} \right)
	\nonumber\\
	&= \ZZ(q) + \sum_{j\ge1} q^j \cdot \ZZ(q)
	\nonumber\\
	&= \frac{1}{1-q}\ZZ(q).
    \nonumber
\qedhere
\end{align}
\end{proof}
\end{prop}

\subsection{The skew Ferrers operators}
\label{sec:madonna-scannata}

In this section we start setting up the main ingredients for the proof of \Cref{thm:intro:Z_k}.

\begin{notation-construction}
\label{not:SK-paths}
Fix a connected skew Ferrers diagram $\lambda$ consisting of $D$ boxes. Notice that each such diagram has a unique north-easternmost point $P_{\NE}(\lambda)$ and a unique south-westernmost point $P_{\SW}(\lambda)$. There are exactly two paths on the boundary of the diagram connecting these two points: we will refer to them as the \emph{north-west path} \textemdash\; the path that starting from $P_{\NE}(\lambda)$ has first a step west \textemdash\; and the \emph{south-east path} \textemdash\; the one that starting from $P_{\NE}(\lambda)$ has first a step south.

\begin{figure}[H]
\begin{tikzpicture}
    \draw (0,0)--(3/4,0)--(3/4,2/4)--(5/4,2/4)--(5/4,6/4)--(7/4,6/4)--(7/4,8/4)--(9/4,8/4)--(9/4,11/4)--(14/4,11/4)--(14/4,12/4)--(4/4,12/4)--(4/4,4/4)--(0,4/4)--cycle;

    \draw[very thick] (4/4,12/4)--(4/4,4/4);
    \draw[very thick] (14/4,12/4)--(4/4,12/4);
    \draw[very thick] (4/4,4/4)--(0,4/4);
    \draw[very thick]  (0,4/4)--(0,0);

    \node[right] at (14/4-0.1 , 12/4 +0.2 ) {\small$P_{\NE}(\lambda)$};
    \node[left] at ( 0 , -0.1  ) {\small$P_{\SW}(\lambda)$};

\end{tikzpicture}
  \caption{The north-west path (thick) and the south-east path (thin) of a skew Ferrers diagram $\lambda$.}
  \label{fig:boat2}
\end{figure} 

The definition of a skew Ferrers diagram ensures that the north-west path only consists of steps to the west and steps to the south, with the first step towards west and the last step towards south. Hence the north-west path can be described by
\[
\ell_1, \dots, \ell_M\,\in\,\BZ_{\geqslant 1}
\]
horizontal lengths of $\lambda$ (the consecutive numbers of steps towards west) and by 
\[
v_1, \dots, v_M\,\in\,\BZ_{\geqslant 1}
\]
vertical lengths of $\lambda$ (the consecutive numbers of steps towards south), where the path starts taking $\ell_1$ steps west and ends taking $v_M$ steps south.

\begin{figure}[H] 
\begin{tikzpicture}
    \draw (0,0)--(3/4,0)--(3/4,2/4)--(5/4,2/4)--(5/4,6/4)--(7/4,6/4)--(7/4,8/4)--(9/4,8/4)--(9/4,11/4)--(14/4,11/4)--(14/4,12/4)--(4/4,12/4)--(4/4,4/4)--(0,4/4)--cycle;

    \node[above] at ( 9/4 , 12/4 +0.2  ) {\small$\ell_1$};
    \node[above] at ( 2/4-0.1 , 4/4 +0.2 ) {\small$\ell_2$};
    \node[left] at ( 4/4-0.2 , 8/4+0.1  ) {\small$v_1$};
    \node[left] at ( 0-0.2 , 2/4+0.1  ) {\small$v_2$};

    \draw[thick,->] (4/4-0.2,12/4-0.1)--(4/4-0.2,4/4+0.3);
    \draw[thick,->] (14/4-0.1,12/4+0.2)--(4/4+0.1,12/4+0.2);
    \draw[thick,->] (4/4-0.2,4/4+0.2)--(0+0.1,4/4+0.2);
    \draw[thick,->]  (0-0.2,4/4-0.1)--(0-0.2,0+0.1);
\end{tikzpicture}
  \caption{Vertical and horizontal lengths on the north-west path of a skew Ferrers diagram $\lambda$.}
  \label{fig:boat1}
\end{figure}   

\begin{remark}
\label{rmk:when-M-is-1}
In the above construction, $\lambda$ is a Ferrers diagram (not just a skew Ferrers diagram) if and only if $M = 1$.
\end{remark}

For future use, we set
\begin{equation}
\label{eqn:NW-len}
L = \sum_{i=1}^M \ell_i, \qquad \qquad V = \sum_{i=1}^M v_i, \qquad \qquad L_{\NW} = V + L.
\end{equation}
We refer to $L_{\NW}$  as the length of the north-west path of $\lambda$. Notice that 
\begin{equation}
 L_{\NW} - 1 \leqslant D,
\end{equation}
and the bound is sharp (e.g. the one-part partition realises the equality). It should be noted that the quantities $V,L,M,L_{\NW}$ depend on $\lambda$, but we omit to write the dependence explicitly to simplify the notation throughout. Moreover, for disconnected skew Ferrers diagrams we define their quantities $V, L, M, L_{\mathrm{NW}}$ as the sum, or the union, of those of their connected components.

Finally, for a pair $(j,h) \in \BN \times \BN$, we introduce a formal variable
\[
y_{j,h},
\]
along with a differential operator
\[
\frac{\mathrm{d}}{\mathrm{d} y_{j,h}}.
\] 
 \end{notation-construction}

\begin{definition}
\label{def:SK-operator}
Relying on \Cref{not:SK-paths}, we define the \emph{skew Ferrers operator} relative to a connected skew Ferrers diagram $\lambda$ as the differential operator
\begin{align} 
\label{eq:defTlambda}
    T_{\lambda} &= 
    \sum_{j \geq 0} \prod_{k=1}^M \left[ \left( \frac{\mathrm{d}}{\mathrm{d}y_{j+\sum_{i=k+1}^M \ell_i,v_k}} - \frac{\mathrm{d}}{\mathrm{d}y_{j + \sum_{i=k+1}^M \ell_i,v_k+1}} \delta_{k < M}\right) \prod_{p=0}^{\ell_k-2} \frac{\mathrm{d}}{\mathrm{d}y_{j+\sum_{i=k}^{M} \ell_i - p - 1, 0}}\right].
\end{align}
\end{definition}

Here we use the convention that empty sums vanish. Clearly all such differential operators commute with each other. 

Skew Ferrers operators corresponding to disconnected skew Ferrers diagrams $\lambda$ with connected components $\lambda^{(1)}, \ldots, \lambda^{(N)}$ are defined in terms of the operators of their connected components $T_{\lambda^{(1)}}, \ldots, T_{\lambda^{(N)}}$ in the following way. Let $S_{[j, j+ L-1]}$ be the $j$-th summand in \Cref{eq:defTlambda}, so that for a connected skew Ferrers diagram one has $T_{\lambda} = \sum_{j \ge 0} S_{[j, j+L-1]}$. Notice that $S_{[j, j+ L-1]}$ is a homogeneous polynomial of degree $L$ in the derivatives 
\[
\frac{\mathrm d}{\mathrm dy_{t, \bullet}}, \qquad\qquad t \in \{ j, \dots, j + L-1\},
\]
where each summand of $S_{[j, j+ L-1]}$ contains exactly one $\frac{\mathrm d}{\mathrm dy_{t, \bullet}}$ for each $t \in \{ j, \dots, j + L-1\}$ and for some value of the second indices of the $y$-variables. 

Given $\lambda = \lambda^{(1)}\amalg \cdots \amalg \lambda^{(N)}$ as above, define
\[
\CI_\lambda = \Set{(I^{(1)},\ldots,I^{(N)}) |
\begin{array}{ccc}
I^{(h)} \subset \BN \mbox{ is an interval of length }L(\lambda^{(h)}) \\
\mbox{for all }h=1,\ldots,N,\mbox{ and }I^{(i)}\cap I^{(j)} = \emptyset \mbox{ for }i \neq j
\end{array}
}.
\]
For $i = 1, \ldots, N$, let $S^{(i)}_I$ be the summand of $T_{\lambda^{(i)}}$ corresponding to the interval $I$ of first indices according to the discussion above, all of whose elements are consecutive nonnegative integers and $\lvert I\rvert = L(\lambda^{(i)})$.  The operator $T_{\lambda}$ is defined by 
\begin{equation}
T_{\lambda} = \frac{1}{\lvert\Sym(\lambda)\rvert}\sum_{(I^{(1)},\ldots,I^{(N)}) \in \CI_\lambda} S_{I_1}^{(1)} \cdots S_{I_N}^{(N)},
\end{equation}
where $\lvert\Sym(\lambda)\rvert = \prod_{\rho \in \mathscr{Q}} \mathfrak s_{\rho}(\lambda)!$, in which $\mathscr Q$ was introduced in \Cref{notation:everything} and $\mathfrak s_{\rho}(\lambda)$ is the number of connected components of $\lambda$ equal to $\rho$. Again, the operators $S_{I_i}^{(i)}$ relative to disjoint indices subsets commute with each other.

\begin{example}
The operators associated to skew Ferrers diagrams with up to 3 boxes are the following.
\begin{align*}
    T_{\scalebox{0.25}{\yng(1)}} &= \sum_{j\ge0} \frac{\mathrm{d}}{\mathrm{d} y_{j,1}} 
    \qquad \qquad \qquad     
    &T_{\scalebox{0.25}{\yng(1)}, \, \scalebox{0.25}{\yng(1)}, \, \scalebox{0.25}{\yng(1)}} &= \sum_{0 \leqslant j_1 < j_2 < j_3 } \frac{\mathrm{d}}{\mathrm{d} y_{j_1,1}} \frac{\mathrm{d}}{\mathrm{d} y_{j_2,1}} \frac{\mathrm{d}}{\mathrm{d} y_{j_3,1}} 
    \\
    T_{\scalebox{0.25}{\yng(1)},\, \scalebox{0.25}{\yng(1)}} &= \sum_{0 \leqslant j_1 < j_2} \frac{\mathrm{d}}{\mathrm{d} y_{j_1,1}} \frac{\mathrm{d}}{\mathrm{d} y_{j_2,1}}
    \qquad \qquad \qquad
    &T_{\scalebox{0.25}{\yng(1)},\, \scalebox{0.25}{\yng(1,1)}} &= \sum_{\substack{0 \leqslant j_1, j_2 \\ j_1\neq j_2}} 
        \frac{\mathrm{d}}{\mathrm{d}y_{j_1,1}}\frac{\mathrm{d}}{\mathrm{d}y_{j_2,2}}
    \\
    T_{\scalebox{0.25}{\yng(2)}} &= \sum_{j\ge0} \frac{\mathrm{d}}{\mathrm{d}y_{j,1}}\frac{\mathrm{d}}{\mathrm{d}y_{j+1,0}} 
    \qquad \qquad \qquad
    & T_{\scalebox{0.25}{\yng(1)},\, \scalebox{0.25}{\yng(2)}} &= \sum_{\substack{0 \leqslant j_1, j_2 \\ j_1 \notin \{j_2-1,j_2\}}} \frac{\mathrm{d}}{\mathrm{d} y_{j_1,1}} \frac{\mathrm{d}}{\mathrm{d} y_{j_2-1,1}} \frac{\mathrm{d}}{\mathrm{d} y_{j_2,0}}
    \\
    T_{\scalebox{0.25}{\yng(1,1)}} &= \sum_{j\ge 0} \frac{\mathrm{d}}{\mathrm{d} y_{j,2}} 
    \qquad \qquad \qquad
    &T_{\scalebox{0.25}{\yng(2,1)}} &= \sum_{j \ge 0} \frac{\mathrm{d}}{\mathrm{d}y_{j,2}} \frac{\mathrm{d}}{\mathrm{d}y_{j+1,0}}
    \\
    T_{\scalebox{0.25}{\yng(1,1,1)}} &= \sum_{j\ge0} \frac{\mathrm{d}}{\mathrm{d}y_{j,3}} 
    \qquad \qquad \qquad
    & T_{\scalebox{0.25}{\yng(3)}} &= \sum_{j \ge 0} \frac{\mathrm{d}}{\mathrm{d}y_{j,1}} \frac{\mathrm{d}}{\mathrm{d}y_{j+1,0}} \frac{\mathrm{d}}{\mathrm{d}y_{j+2,0}}
    \\
    &
    \qquad \qquad \qquad
    & T_{\rotatebox[origin=c]{90}{\scalebox{0.25}{\yng(1,2)}}} &= \sum_{j \ge 0} \frac{\mathrm{d}}{\mathrm{d}y_{j,1}}\left(\frac{\mathrm{d}}{\mathrm{d}y_{j+1,1}} - \frac{\mathrm{d}}{\mathrm{d}y_{j+1,2}}\right).
\end{align*}
\end{example}

\begin{remark}\label{rem:Tdependence}
It is clear from \Cref{eq:defTlambda} that $T_\lambda$ only depends on the sequence $(\ell_1,v_1,\dots,\ell_M,v_M)$. In particular, whenever two skew Ferrers diagrams $\lambda\in\mathscr Q^{[n]}$, $\lambda'\in\mathscr Q^{[n']}$ share the same north-west path, one has $T_\lambda=T_{\lambda'}$. See also \Cref{rmk:transpose}.
\end{remark}

\subsection{\texorpdfstring{The series $\FFZ_D$ in terms of the $D$-decorated partition function}{}}

The goal of this section is to prove an analogue of \Cref{lemma:Z1-and-derivative} for the series $\FFZ_D(q)$ introduced in \Cref{eqn:Z-series}, for all $D\ge 1$. This will be achieved in \Cref{cor:TD-applied-to-FD} as an application of the main theorem of this section, namely \Cref{thm:skewferret}. A few examples of generating functions are collected in \Cref{sec:numerics}.

We first have to define an analogue of the series $F(x;q)$ defined in \Cref{F-series}. We do it as follows.

\begin{definition}
\label{def:F-partition-function}
For $D\in \BZ_{>0}$, define the $D$-\emph{decorated partition function} as
\[
F^{[D]}(\vec{y}; q) = 
    \left(\prod_{k \geq 1}^D y_{0,k} \right)
    \prod_{j\ge1} 
    \left( 
    y_{j,0} + y_{j,1}q^j + y_{j,1}y_{j,2} q^{2j} + \dots +
    \left(\prod_{k=1}^D y_{j,k}\right) \sum_{m_j \ge D} q^{jm_j}
    \right).
\]
\end{definition}

\begin{remark}
\label{rmk:completion_FD}
In principle, one could view $F^{[D]}(\vec y;q)$ as an element of $\widehat A\llbracket q\rrbracket$, where $\widehat A$ is the completion 
    \[
    \widehat A=\lim_{\substack{\longleftarrow\\n}}\BZ\left[\vec y\right]\big/I^n
    \]
    with respect to the augmentation ideal
    \[
    I=\sum_{j\ge 1}\left(y_{j,0}-1,y_{j,1},\dots,y_{j,D}\right)\subset\BZ\left[\vec y\right].
    \]
\end{remark}

Sometimes in the following we will use the shorthand $F = F^{[D]}$, with the convention that $F$ is decorated with as many variables as needed depending on the operator which is applied to it. We also use $F \big{|}_1$ as an abbreviation of $F \big{|}_{\vec{y} = 1}$, which means to set the values of all variables to 1.

For a possibly disconnected skew Ferrers diagram $\lambda \in \mathscr Q^{[D]}$, define
\[
A_{n-D}(\lambda)
\]
to be the number of ways to add $\lambda$ to all partitions of size $n-D$. Form the generating function
\begin{equation}
\label{Z-lambda}
\FFZ_\lambda(q) = \sum_{n\geqslant D} A_{n-D}(\lambda) q^n\in\BZ\llbracket q\rrbracket.    
\end{equation}

\begin{theorem} 
\label{thm:skewferret} 
Fix a possibly disconnected skew Ferrers diagram $\lambda \in \mathscr{Q}^{[D]}$. Then, there is an identity
\begin{equation}
\FFZ_\lambda(q) = T_{\lambda}.F^{[D]}\Bigg|_{\vec{y} = 1}.
\end{equation} 
\end{theorem}

\begin{proof}
    The proof relies on the definition of the skew Ferrers operators in \Cref{eq:defTlambda} in relation to how the $y$-variables decorate the partition function $F^{[D]}$. 
    
    We start by looking at the case of $\lambda$ being a connected skew Ferrers diagram. The idea is that the differential operators describe step-by-step the north-west boundary of the skew Ferrers diagram. Moreover, the $y$-variables decorate $F$ in such a way that the application of the operator $T_{\lambda}$ selects from the partition function $F$ the partitions $\rho$ in which $\lambda$ can be added, and weighs each such $\rho$ by the number of different possible ways of doing that \textemdash \; i.e. the number of south-east subpaths of $\rho$ that coincide with the north-west path of $\lambda$. 
    
    After the application of the operator $T_{\lambda}$ the $y$-variables are set to one, returning the $q$-generating series enumerating the number of ways of adding the diagram $\lambda$ in any partition. We are left with checking that the particular sum of products of derivatives indeed does the promised job.

    Let us study the interaction between the operator $T_{\lambda}$ and the $y$-variables decoration in $F^{[D]}$. Let us first notice that there is a difference between the role of the $y_{j,0}$  and of the $y_{j,h}$ variables.
    \begin{enumerate}
    \item Differentiating by $y_{j,0}$ means selecting from the generating series $F$ exactly those partitions $\rho$ that do not have parts equal to $j$, that is, imposing the condition $m_j(\rho) = 0.$
    \item Differentiating by $y_{j,h}$ for $h \geqslant1$ means instead selecting from the generating series $F$ exactly those partitions $\rho$ that have \textit{at least} $h$ parts equal to $j$ (as opposed to \textit{exactly} $h$ parts), that is, imposing the condition $m_j(\rho) \geqslant h.$
    \end{enumerate}
Summarising, one can think of the derivatives in the $y$-variables as having the following effect on the generating series $F$ of all partitions $\rho$:
\[
\frac{\mathrm{d}}{\mathrm{d}y_{j,0}} \;\;\; \longleftrightarrow \;\;\;  m_j(\rho) = 0,
\qquad \qquad 
\frac{\mathrm{d}}{\mathrm{d}y_{j,h}} \;\;\; \longleftrightarrow \;\;\;  m_j(\rho) \geqslant h \;\; \text{ for } \;\; h\geqslant1,
\]
where the double arrows mean that the application of the differential operator to the partition function annihilates all the partitions that do not satisfy the property on the right. The composition of operators corresponds to taking the intersection of the resulting conditions imposed on the partitions. This shows that operators and generating series are compatible in the sense that there is a well-defined vocabulary that associates derivatives and conditions on partitions. In order to prove that $T_{\lambda}$ actually describes the north-west boundary of $\lambda$, we can now reason purely in terms of which conditions we need to impose on the multiplicities $m_j(\rho)$. 
    
Let us hence describe the conditions to be imposed on the $m_j(\rho)$ so that a subpath of the south-east path of $\rho$ is equal to the north-west path of $\lambda$. We start from the first horizontal length $\ell_1$. Note that every step to the west after the first one requires the vanishing of some $m_j(\rho)$, therefore, if we need $\ell_1$ west steps starting from parts of $\rho$ equal to $j$, we need to impose the conditions
    $$
    m_{j-1}(\rho) = m_{j-2}(\rho) = \dots = m_{j - \ell_1 + 1}(\rho) = 0,
    $$
    which corresponds to the composition of the operators
    $$
    \frac{\mathrm{d}}{\mathrm{d}y_{j-\ell_1+1,0}}\cdots \frac{\mathrm{d}}{\mathrm{d}y_{j-1,0}}.
    $$
    The reason why there are $\ell_1 - 1$ conditions and not $\ell_1$ is because the first west step is granted by default when   part length changes in a Ferrers diagram. See \Cref{fig:horiz_lengths}.
    
    \begin{figure}[H]
\begin{tikzpicture}
    \draw[ultra thick] (1/2,-8/2) --(0,-8/2)-- (0,0)--(5/2,0)--(5/2,-4/2)--(1/2,-4/2)  --cycle;
    \draw[ultra thick] (4/2,-4/2)--(4/2,-5/2)--(2/2,-5/2)--(2/2,-6/2)--(1/2,-6/2)--(1/2,-4/2) --cycle;
    \draw[ultra thin] (0,-1/2)--(5/2,-1/2);
    \draw[ultra thin] (0,-2/2)--(5/2,-2/2);
    \draw[ultra thin] (0,-3/2)--(5/2,-3/2);
    \draw[ultra thin] (0,-4/2)--(5/2,-4/2);
    \draw[ultra thin] (0,-5/2)--(2/2,-5/2);
    \draw[ultra thin] (0,-6/2)--(1/2,-6/2);
    \draw[ultra thin] (0,-7/2)--(1/2,-7/2);
    
    \draw[ultra thin] (1/2,0)--(1/2,-4/2);
    \draw[ultra thin] (2/2,0)--(2/2,-5/2);
    \draw[ultra thin] (3/2,0)--(3/2,-5/2);
    \draw[ultra thin] (4/2,0)--(4/2,-4/2);
\end{tikzpicture} 
  \caption{Let $\rho = (5^{4} 4^0 3^0 2^0 1^{4})$ and let $\lambda$ (for simplicity also a partition $\lambda = (3^1 2^0 1^1)$) with $\ell_1 = 3$ and $v_1 = 2$. The first $\ell_1 = 3$ steps west are described for some $j$ by the operator $\frac{\mathrm{d}}{\mathrm{d}y_{j-2,0}}\frac{\mathrm{d}}{\mathrm{d}y_{j-1,0}}$, which selects the partition $\rho$ for the two values of $j=5$ \textemdash\; since $m_{4}(\rho) = m_{3}(\rho) = 0$ \textemdash\; and of $j=4$ \textemdash\; since $m_{3}(\rho) = m_{2}(\rho) = 0$. The summand for $j=5$ will not survive the application of differential operator for the vertical part, but the $j=4$ summand will. Notice that the $j=4$ summand requires $m_{3}(\rho) = m_{2}(\rho) = 0$ but the fact that $m_4(\rho)$ additionally vanishes does not spoil the application of the operator.
  }
  \label{fig:horiz_lengths}
\end{figure} 

    We now have to describe $v_1$ (consecutive) steps south, which corresponds  to having exactly $v_1$ parts in $\rho$ of the correct length, this being determined by the previous operator in such a way that the first indices of the $y$-variables are consecutive. Recall however that $y$-variables with positive second index correspond to taking that amount of parts \textit{or more}, hence we need to impose (see \Cref{fig:vert_lengths}) the two conditions
    \[
    m_{j-\ell_1}(\rho) \geqslant v_1 \;\; \text{ and } \;\; m_{j-\ell_1}(\rho) \ngeq v_1 + 1 \qquad\qquad  \left(\text{ i.e. } m_{j-\ell_1}(\rho) = v_1 \right),
    \]
    which can be achieved by taking the difference of the corresponding two operators 
    \[
    \frac{\mathrm{d}}{\mathrm{d}y_{j-\ell_1,v_1}} \; - \; \frac{\mathrm{d}}{\mathrm{d}y_{j-\ell_1,v_1+1}}.
    \]
    
\begin{figure}[H]\begin{tikzpicture}
        \draw (0,0)--(1/3,0)--(1/3,5/3)--(7/3,5/3)--(7/3,9/3)--(12/3,9/3)--(12/3,11/3)--(0,11/3)--cycle;
        \draw (1/3,2/3)--(8/3,2/3)--(8/3,5/3)--(10/3,5/3)--(10/3,7/3) --(11/3,7/3)--(11/3,9/3) ;
        \draw[|-|] (0.2,7/3)--(7/3-0.2,7/3);
        \node[above] at (3.5/3,7/3) {\tiny $k=7$};
        \draw[|-|] (7/3+0.2,5/3+0.2)--(7/3+0.2,9/3-0.2);
        \node[right] at (7/3+0.2,7/3) {\tiny $v_1=4$};
    \end{tikzpicture} 
  \caption{Here we see that $\lambda$ is such that $v_1=4$, and that $v_1$ is not the last vertical length of the skew Ferrers diagram (i.e. $M >1$). Therefore one of the conditions for $\rho$ is that the south-east path of $\rho$ needs to have exactly $4$ south steps at some point. If $v_1$ had been the last vertical length of $\lambda$, then $\rho$ could have had \textit{at least} $4$ steps south along the boundary matching $v_1$ and still be compatible. In order to have $4$ south steps, $\rho$ must have a multiplicity equal to $4$, in this case $m_7(\rho) = 4$. This corresponds to the operator $\frac{\mathrm{d}}{\mathrm{d}y_{7,4}} - \frac{\mathrm{d}}{\mathrm{d}y_{7,5}}$.
  }
  \label{fig:vert_lengths}
\end{figure}

One continues the same argument for the subsequent $\ell_2$ steps west and $v_2$ steps south, and so on, up to the last $\ell_M$ steps west and $v_M$ steps south. Notice that for the latter we do not need exactly $v_M$ steps but they could be possibly more, hence there is no need to take the difference of the two operators. That is why we set the last vertical operator to be zero through the Kronecker delta (for which we use the notational convention introduced in \Cref{sec:conventions}).

Finally, summing over $j$ accounts for all possible additions. The prefactor $\prod_{k=1}^D y_{0,k}$ accounts for the additions of a skew diagram $\lambda$ so that the segment of length $v_{M}$ is adjacent to the wall instead of to the partition $\rho$, and it specialises to the multiplication of the single power of $x$ in \Cref{lemma:Z1-and-derivative}.

The theorem clearly extends to disconnected skew Ferrers diagrams since their operators
are defined in terms of the summands of the operators of their connected components. This concludes the proof of the theorem.
\end{proof}

\begin{definition}
Fix a positive integer $D$. Define the $D$-\emph{th skew Ferrers operator} $T_D$ as
\begin{equation}
    T_D = \sum_{\lambda\in \mathscr{Q}^{[D]}} T_{\lambda}.
\end{equation}
\end{definition}

\begin{corollary} 
\label{cor:TD-applied-to-FD}
Fix a positive integer $D$. The operator $T_D$ acting on the $D$-decorated partition function produces the generating series of the Euler characteristics of the $[n,n+D]$-nested punctual Hilbert schemes on a smooth surface $S$, i.e.
    \begin{equation}
        \FFZ_D(q)=T_D.F^{[D]} \Bigg{|}_{\vec{y} = 1} .
    \end{equation}
\end{corollary}
    \begin{proof}
    The claim follows from \Cref{thm:skewferret} and the fact that for any pair of nested partitions $\lambda_n\subset\lambda_{n+D}$, for $\lambda_n\in\mathrm P^{[n]}$ and $\lambda_{n+D}\in\mathrm P^{[n+D]}$, their set-theoretic difference $\lambda_{n+D}\smallsetminus\lambda_{n}$ has class $[\lambda_{n+D}\smallsetminus\lambda_{n}]\in\mathscr Q^{[D]}$. 
    \end{proof}

\subsection{Transposition invariance of \texorpdfstring{$\FFZ_\lambda$}{}}
\label{rmk:transpose}

Given a possibly disconnected skew Ferrers diagram $\lambda$ of size $D$, and its transpose $\lambda^T$, one has in general that
\begin{equation}
    T_{\lambda}(\vec{y}) \neq T_{\lambda^{T}}(\vec{y})
\end{equation}
as differential operators, unless $\lambda$ and $\lambda^{T}$ have equal north-west path, see \Cref{rem:Tdependence}. However, in the next result we observe a symmetry property of the operators $T_{\lambda}$, \emph{after} their application to $F$ and the specialisation of the $y$-variables; we refer to this property as \textit{transposition invariance}.

\begin{lemma}[Transposition invariance] 
\label{lem:tinvariance} 
Fix a possibly disconnected skew Ferrers diagram $\lambda\in\mathscr Q^{[D]}$. There is an identity of $q$-series
\begin{equation}
    \FFZ_{\lambda}(q) = \FFZ_{\lambda^{T}}(q),
\end{equation}
where $\lambda^T$ denotes the transposed diagram as in \Cref{def:transposed_diagram}, and $\FFZ_\lambda$ was defined in \Cref{Z-lambda}.
\end{lemma}

\begin{proof}
It follows from the fact that transposition is a global involution on the set of all partitions.
\end{proof}

\section{Proof of \texorpdfstring{\Cref{thm:intro:Z_k}}{}}
\label{sec:proof-B}

In this section we prove \Cref{thm:intro:Z_k} as announced in the introduction. We start in \Cref{subsec:bassecase} by analysing the structure of the ratio $\FFZ_D / \ZZ$ for $D\ge 1$, which settles the case $s=1$. We explicitly compute the cases $D=2,3$ in \Cref{sec:in-action}. In \Cref{subsec:g.f.tantinest} we deduce a general argument to describe the ratio $\FFZ_{\bk} / \ZZ$ for an arbitrary $\bk\in\BZ^s_{\ge 0}$. 

\subsection{The base case \texorpdfstring{$s=1$}{}} \label{subsec:bassecase}
Recall the definition of the series $\FFZ_\lambda$ introduced in \eqref{Z-lambda}. The first and main step of the proof of \Cref{thm:intro:Z_k} is the following proposition on the ratio $\FFZ_\lambda/\ZZ$.

\begin{prop} 
\label{prop:Zlambda} 
Fix a possibly disconnected skew Ferrers diagram $\lambda\in\mathscr Q^{[D]}$. Then, there is an identity
\begin{equation}
\label{eqn:ratio-formula}
    \frac{\FFZ_{\lambda}(q)}{\ZZ(q)} 
    =
    \frac{\widetilde{\mathsf P}_{\lambda}(q)}{\prod_{i=\max(L,V)}^{L_{\NW}(\lambda) - 1} (1 - q^{i})},
\end{equation} 
where $\widetilde{\mathsf P}_{\lambda}(q)$ is a polynomial divisible by $q^B$, of degree bounded by 
$\binom{L_{\NW}-1}{2} + B $, for $B = \sum_{i=1}^M \ell_i \left( \sum_{k=1}^{i-1} v_k \right)$.
\begin{proof}

The proof follows from a careful analysis of the structure of the skew Ferrers differential operators and how they act on $F^{[D]}$. Let us first restrict ourselves to connected skew Ferrers diagrams and let us explicitly compute the action of an arbitrary such skew Ferrers $\lambda$ diagram on $F^{[D]}$: 
\begin{equation}
\label{eq:applicationTlambda}
\begin{split} 
    \tilde{T}_{\lambda} .F^{[D]} \big{|}_1 &= 
    \sum_{j\ge 1} \prod_{k=1}^M \left[ \left( \frac{\mathrm{d}}{\mathrm{d}y_{j+\sum_{i=k+1}^M \ell_i,v_k}} - \frac{\mathrm{d}}{\mathrm{d}y_{j + \sum_{i=k+1}^M \ell_i,v_k+1}} \delta_{k < M}\right) \prod_{p=0}^{\ell_k - 2} \frac{\mathrm{d}}{\mathrm{d}y_{j+\sum_{i=k}^{M} \ell_i - p - 1, 0}}\right].F^{[D]}\big{|}_1
     \\
    &= 
    \sum_{j\ge 1 } \frac{1}{(1 - q^j)} \prod_{k=1}^M \left[ (1 - q^{(j+\sum_{i=k+1}^M \ell_i)})q^{(j+\sum_{i=k+1}^M \ell_i)v_k} \prod_{p=0}^{\ell_k - 2} (1 - q^{j+\sum_{i=k}^{M} \ell_i - p - 1}) \right] \cdot \ZZ(q)
    \\ 
    &=
    \left[ \sum_{j\ge 1} q^{S(j, \lambda)} \prod_{p=1}^{L - 1}  (1 - q^{j+p})\right] \cdot \ZZ(q)
\end{split}
\end{equation}
where $S(j,\lambda)$ is the following degree one polynomial in $j$:
\begin{equation}
S(j,\lambda) = j V + B; \qquad \qquad B = \sum_{i=1}^M \ell_i \left( \sum_{k=1}^{i-1} v_k \right),
\end{equation}
and $\tilde{T}_{\lambda}$ differs from the operator $T_{\lambda}$ by the $j=0$ summand. Removing the $j=0$ cannot decrease the bound, in fact it increases it, making it less sharp, as removing the initial term of a geometric series increase the degree of the polynomial at the numerator, e.g. $\frac{1}{1-q} - 1 = \frac{q}{1-q}$. Even if the bound loses sharpness, the removal of the $j=0$ term makes the computation much more manageable.

Observe that each $j$-summand is a polynomial in $q$ whose exponents $S(j,\lambda)$ are 
polynomials in $j$ of degree exactly one (but not necessarily homogeneous). This condition is crucial for the statement to be true, and it allows us to compute further

\begin{align}
\label{eqn:bagascia}
\frac{\FFZ_{\lambda}(q)}{\ZZ(q)}&= \nonumber
\sum_{j\ge0} q^{S(j, \lambda)} \prod_{p=1}^{L - 1}  (1 - q^{j+p})
    \\
    &=
    \sum_{j\ge 0} \left( (-1)^{L-1} q^{S(j, \lambda) + (L-1)(j + \frac{L}{2})} + \dots + q^{S(j, \lambda)}\right)
    \nonumber \\ 
    &= \sum_{j\ge 0}  \sum_{t=0}^{L-1} (-1)^t q^{j(V+t) + B} \!\!\!\!\!\! \sum_{1 \leqslant i_1 < \dots < i_t \leqslant L-1 } \!\!\!\!\!\! q^{i_1 + \dots + i_t}
    \nonumber \\
    &= \sum_{j\ge0} \sum_{t=0}^{L-1} (-1)^t \sum_{k=0}^{U} p_{\binom{t+1}{2} + k,t} q^{j(V+t) + B + \binom{t+1}{2} + k}
    \nonumber \\ 
    &= q^B \sum_{t=0}^{L-1} (-1)^t \frac{ \sum_{k=0}^{U} p_{\binom{t+1}{2} + k,t} q^{\binom{t+1}{2} + k}}{1 - q^{V+t}} 
    \nonumber \\
    &= \frac{\mathsf P'_{\lambda}(q)}{(1 - q^{V})(1 - q^{V+1}) \cdots (1 - q^{V + L - 1})} 
\end{align}
where $p_{n,t}$ is the number of partitions of size $n$ of length $t$, and
\[
U = \binom{L}{2} - \binom{L-t}{2} - \binom{t+1}{2}.
\]
Note that the polynomial $\mathsf P'$ is essentially the polynomial $\widetilde{\mathsf P}$ in the statement, although it may in principle differ by a suitable product of polynomials of the form $(1-q^d)$ depending on whether $L \leq V$.

Recalling that $V + L = L_{\NW}$, we see that $\FFZ_{\lambda}(q)/\ZZ(q)$ is a rational function in $q$ with denominator dividing $\prod_{i=V}^{L_{\NW}-1}(1 - q^i)$. By transposition invariance in \Cref{lem:tinvariance} we can compute without loss of generality either with $\lambda$ or its transpose $\lambda^{T}$. Out of the two choices we pick the one satisfying the condition 

$$L \le V.$$ 

This refines the denominator showing that it divides $\prod_{i=\max(V,L)}^{L_{\NW}-1}(1 - q^i)$. The argument so far proves the rationality statement and the claim about the denominator.

Let us now study the numerator, that is, the polynomial $\mathsf P'_{\lambda}(q)$. From the second to last step of the computation above it is clear that $\mathsf P'_{\lambda}(q)$ is divisible by $q^B$. Recall that the transposition of the skew diagram at the level of the vertical and horizontal lengths takes the form 
$$
\sigma: \ell_i \longleftrightarrow v_{M+1-i}, \qquad \qquad i=1, \dots, M.
$$
This means that the bilinear form $B$ is invariant under $\sigma$, so it remains uneffected by the assumption $L \le V$. Notice moreover that transposition invariance from \Cref{lem:tinvariance} at the numerator takes the form
\[
\mathsf P_{\lambda}(q) = \prod_{j=1}^{V-1}(1-q^j)\mathsf P_{\lambda}'(q) = \prod_{j=1}^{L-1}(1-q^j)\mathsf P_{\lambda^T}'(q).
\]
We now bound the degree of $\mathsf P_{\lambda}(q)$. To this purpose we need to select the maximum among all powers of $q$ after multiplying and dividing by the denominator factors missing in the last expression of \eqref{eqn:bagascia} so to obtain precisely the denominator appearing in \eqref{eqn:ratio-formula}. In other words, we need to compute the quantity
\begin{equation}
\max_{\substack{1 \leqslant t \leqslant L-1 \\ 0 \leqslant k \leqslant U}}
\left\{ \binom{L_{\NW}}{2} - V - t + \binom{t+1}{2} + k + B\right\}.
\end{equation}

The maximum is realised by taking both the maximum index $t$ and the maximum index $k$, from which we obtain the following quantity.
\begin{align*}
\deg_q \mathsf P_{\lambda}(q) &\leqslant \binom{L_{\NW}}{2} - V - t + \binom{t+1}{2} + k + B \Big{|}_{k=U(t), t=L-1}
\\
&=\binom{L_{\NW}}{2} - V - t  + \binom{L}{2} - \binom{L-t}{2}  + B  \Big{|}_{t=L-1}
 \\
 &=\binom{L_{\NW}}{2} - (L_{\NW} - 1) + \binom{L}{2} + B
 \\
 &= \binom{L_{\NW}-1}{2} + \binom{L}{2} + B.
\end{align*} 

Therefore $\deg_q \mathsf P'_{\lambda}(q) \leqslant \binom{L_{\NW}-1}{2} + \binom{L}{2} - \binom{V}{2} + B,$ which by the hypothesis $L \le V$ can be bounded as
\begin{align*}
\deg_q \mathsf P'_{\lambda}(q) &\leqslant \binom{L_{\NW}-1}{2} + B.
\end{align*} 
It is easy to see that the number of disconnected components does not raise or lower the degree from \Cref{eq:applicationTlambda}. This concludes the proof of the degree bound for $\mathsf P_{\lambda}'(q)$ and of the proposition.
\end{proof}
\end{prop}

\begin{remark}
By \Cref{rmk:when-M-is-1}, the number $B$ is zero if and only if $\lambda$ is a partition.
\end{remark}

The next theorem settles the case $s=1$ of \Cref{thm:intro:Z_k}.

\begin{theorem}
\label{thm:ZD}
For every $D \in \BZ_{\geq 0}$ there is an identity
\begin{equation}
\label{eqn:ratio}
    \frac{\FFZ_D(q)}{\ZZ(q)} = \frac{\PP_{D}(q)}{\prod_{j=1}^D (1 - q^j)},
\end{equation}
where $\PP_{D}(q)$ is a polynomial of degree at most $\binom{D}{2} + \binom{D-1}{2} + \frac{D^2}{4} = \frac{5}{4}D^2 - 2D + 1$. 
\begin{proof}
Since each operator $T_D$ is a finite sum over the skew Ferrers operators of $D$ boxes, and the properties to be proven are compatible with the sum, the proof follows from \Cref{prop:Zlambda} by summing the contributions for all skew Ferrers diagrams $\lambda$ and by noticing that $L_{\NW} - 1$ is a sharp bound for $D$ among all $\lambda$. For instance, the skew Ferrers diagram $\lambda = (D)$ realises the equality. Similarly $\max(L,V)$ can be bounded by $D-1$ and $B$ by $D^2/4$. This concludes the proof of the theorem. 
\end{proof}
\end{theorem}

\subsection{The construction in action}
\label{sec:in-action}
In this section we perform the explicit calculations underlying \Cref{thm:ZD}, needed to compute $\FFZ_D$ for $D=2,3$, thus deriving explicitly the polynomials $\mathsf P_2$ and $\mathsf P_3$, cf.~\Cref{sec:series1-2} below. See \Cref{sec:some-polynomials} for more examples. In \Cref{sec:poly-values} we compute $\mathsf P_\lambda$ for $\lambda$ a full partition or a fully disconnected skew Ferrers diagram (cf.~\Cref{lemma:specialcases}), as well as the values at 0 and 1 of the polynomials $\mathsf P_\lambda$ and $\mathsf P_D$, cf.~Lemmas \ref{lem:Plambda0110} and \ref{lemma:specialisation}.

\subsubsection{\texorpdfstring{The series $\FFZ_2$ and $\FFZ_3$}{}}
\label{sec:series1-2}

\begin{example}[The case $D=2$] 
\label{ex:D=2}
    Let us see an application of our methods by computing explicitly $\FFZ_2(q)$. For brevity we will display only the variables needed for the application of each skew Ferrers operator and specialise to one the others directly.
    
    For $D=2$ one computes the following three terms:
\begin{align*}
T_{\scalebox{0.25}{\yng(2)}}.F\Bigg{|}_1 &= \sum_{\ell \ge 0} \frac{\mathrm{d}}{\mathrm{d}y_{\ell,1}} \frac{\mathrm{d}}{\mathrm{d}y_{\ell+1,0}} \left(\prod_{k=1}^2 y_{0,k} \right) \prod_{j\ge1}\left(y_{j,0} + y_{j,1} \frac{q^{j}}{1 - q^j}\right)\Big{|}_{y=1} 
        \\
        &=
        \sum_{\ell \ge 1} \frac{q^{\ell}}{1 - q^\ell} \prod_{\substack{j \ge 1\\ j\neq \ell, \ell + 1}} \frac{1}{1 - q^j} + (1-q)\ZZ(q)
        \\
        &= \left(\sum_{\ell \ge 1} (1 - q^{\ell+1})q^{\ell} + 1 - q \right) \ZZ(q) \\
        &= \left(\frac{q}{1-q} - \frac{q^3}{1-q^2} + 1 - q \right)\ZZ(q)
        \\
        &= \frac{1}{(1-q^2)} \ZZ(q),
        \\
        T_{\scalebox{0.25}{\yng(1,1)}}.F\Bigg{|}_1 &= \sum_{\ell \ge 0} \frac{\mathrm{d}}{\mathrm{d}y_{\ell,2}} \left(\prod_{k=1}^2 y_{0,k} \right) \prod_{j\ge1}\left(1 + q^j + y_{j,2} \frac{q^{2j}}{1 - q^j}\right)\Big{|}_{y=1} 
        \\
        &=
        \sum_{\ell \ge 1} \frac{q^{2\ell}}{1 - q^\ell} \prod_{\substack{j\ge 1\\j \neq \ell}} \frac{1}{1 - q^j} + \ZZ(q)
        \\
        &=
        \left(\sum_{\ell \geq 1} q^{2\ell} + 1\right) \ZZ(q)
        \\
        &= \frac{1}{(1-q^2)} \ZZ(q),
        \\
        T_{\scalebox{0.25}{\yng(1)}, \, \scalebox{0.25}{\yng(1)}}.F\Bigg{|}_1 &= \sum_{0 \leqslant \ell < k}\frac{\mathrm{d}}{\mathrm{d}y_{\ell,1}}\frac{\mathrm{d}}{\mathrm{d}y_{k,1}} 
        \left(\prod_{k=1}^2 y_{0,k} \right)
        \prod_{j\ge 1}\left(1 + y_{j,1} \frac{q^j}{1 - q^j}\right)\Big{|}_{y=1} 
        \\
        &=
        \sum_{1 \leqslant \ell < k}\frac{q^k}{1-q^k}\frac{q^\ell}{1-q^{\ell}} 
        \prod_{\substack{j\ge 1\\j \neq \ell,k}}\left(\frac{1}{1 - q^j}\right) + \frac{q}{1-q}\ZZ(q)
        \\
        & = \left(\sum_{\ell\ge 1} q^{\ell} \sum_{k \leq \ell + 1} q^k + \frac{q}{1-q} \right) \ZZ(q) 
        \\
        & = \left(\sum_{\ell\ge 1} q^{2\ell} + 1 \right) \frac{q}{1-q} \ZZ(q) 
        \\
        &= \frac{q}{(1-q)(1-q^2)} \ZZ(q).
    \end{align*}
    By summing up all contributions together one obtains
    \begin{equation}
        \FFZ_2(q) = \left( \frac{q}{(1-q)(1-q^2)} + 2\frac{1}{(1-q^2)}\right) \ZZ(q) = \frac{2 -q}{(1-q)(1-q^2)} \ZZ(q).
    \end{equation}

    Notice that $T_{\scalebox{0.25}{\yng(1,1)}}.F\big{|}_1 = T_{\scalebox{0.25}{\yng(2)}}.F\big{|}_1$ can be concluded a priori since the skew Ferrers diagrams are related by transposition, and so are the partitions to which they are applied, see \Cref{lem:tinvariance}. Therefore, it is enough to compute any of the two instead of both. We however performed the computation to showcase consistency of the operators definition in at least one example. 
\end{example}

\begin{example}[The case $D=3$] 
\label{ex:D=3}
Extending \Cref{ex:D=2}, we now compute the $D=3$ case. Using the decomposition into possibly disconnected skew Ferrers diagrams, one computes the following seven terms. This time we show only the minimal number of required computations, namely one for each class of skew Ferrers diagrams up to tranposition, see \Cref{lem:tinvariance}.  The disconnected terms are:
        \begin{align*}
        T_{\scalebox{0.25}{\yng(1)}, \, \scalebox{0.25}{\yng(1)}, \, \scalebox{0.25}{\yng(1)}}.F\Bigg{|}_1 &= \sum_{0 \leqslant \ell < k < p}\frac{\mathrm{d}}{\mathrm{d}y_{\ell,1}}\frac{\mathrm{d}}{\mathrm{d}y_{k,1}}\frac{\mathrm{d}}{\mathrm{d}y_{p,1}} 
        \left( \prod_{k=1}^3 y_{0,k} \right) \prod_{j\ge 1}\left(1 + y_{j,1} \frac{q^j}{1 - q^j}\right)\Big{|}_{y=1} 
        \\
        &=
        \sum_{1 \leqslant \ell < k < p}\frac{q^k}{1-q^k}\frac{q^\ell}{1-q^{\ell}} \frac{q^p}{1-q^p} \prod_{\substack{j\ge 1\\j \neq \ell,k}}\left(\frac{1}{1 - q^j}\right)
        +
        \sum_{1 \leqslant k < p}\frac{q^k}{1-q^k}\frac{q^p}{1-q^p} \prod_{\substack{j\ge 1\\j \neq k}}\left(\frac{1}{1 - q^j}\right)
        \\
        &= \sum_{\ell\ge 1} q^\ell \sum_{k\ge\ell + 1} \frac{q^{2k+1}}{1-q} \ZZ(q) + \sum_{k \geq 1} q^k \sum_{p \geq k+1} q^p \ZZ(q)
        \\
        &= \frac{q^6}{(1-q)(1-q^2)(1-q^3)} \ZZ(q) + \frac{q^3}{(1-q)(1-q^2)} \ZZ(q),
        \\
        &= \frac{q^3}{(1-q)(1-q^2)(1-q^3)} \ZZ(q),
        \\
        T_{\scalebox{0.25}{\yng(1)}, \, \scalebox{0.25}{\yng(2)}}.F\Bigg{|}_1 &= T_{\scalebox{0.25}{\yng(1)}, \, \scalebox{0.25}{\yng(1,1)}}.F\Bigg{|}_1 \\
        &= 
        \sum_{\substack{\ell,k \ge 0 \\ \ell \neq k}} 
        \frac{\mathrm{d}}{\mathrm{d}y_{k,1}}\frac{\mathrm{d}}{\mathrm{d}y_{\ell,2}} \left( \prod_{p=1}^3 y_{0,p} \right)
        \prod_{j\ge1}\left(1 + y_{j,1} q^j + y_{j,2} \frac{q^{2j}}{1 - q^j}\right)\Big{|}_{y=1}
        \\
        &=
        \sum_{\substack{\ell,k \ge 1 \\ \ell \neq k}} \frac{q^{k}}{1 - q^k} \frac{q^{2\ell}}{1 - q^\ell} \prod_{\substack{j\ge 1 \\j \neq \ell,k}} \frac{1}{1 - q^j}  + \left( \sum_{k \geq 1} q^k + \sum_{\ell \geq 1} q^{2j} \right) \ZZ(q)
        \\
        &=
        \left(\sum_{\ell \ge 1} q^{2\ell} \sum_{\substack{k\ge 1 \\ k \neq \ell}} q^k\right) \ZZ(q) + \frac{q}{1-q} \ZZ(q) + \frac{q^2}{1-q^2} \ZZ(q)
        \nonumber 
        \\
        &= \left(\frac{q^3}{(1-q)(1-q^2)} - \frac{q^3}{(1-q^3)} + \frac{q}{1-q} + \frac{q^2}{1-q^2} \right)\ZZ(q). 
        \\
        &= \left(\frac{1}{(1-q)(1-q^2)} - \frac{1}{(1-q^3)}\right)\ZZ(q). 
        \end{align*}
The connected terms are:
        \begin{align*}
        T_{\scalebox{0.25}{\yng(2,1)}}.F\Bigg{|}_1 =& \sum_{\ell \ge 0} \frac{\mathrm{d}}{\mathrm{d}y_{\ell,2}} \frac{\mathrm{d}}{\mathrm{d}y_{\ell+1,0}} 
        \left(\prod_{k=1}^3 y_{0,k} \right)
        \prod_{j\ge 1}\left(y_{j,0} + y_{j,1}q^j + y_{j,1}y_{j,2} \frac{q^{2j}}{1 - q^j}\right)\Big{|}_{y=1} 
        \\
        =& \sum_{\ell \ge 1} \frac{q^{2\ell}}{(1 - q^{\ell})} \prod_{\substack{j\ge1\\j \neq \ell, \ell+1}}\frac{1}{1 - q^j} + (1-q)\ZZ(q)
        \\
        &=
        \left(\sum_{\ell \ge 1} q^{2\ell}(1 - q^{\ell+1}) + 1 - q \right)\ZZ(q)
        \\
        &=
        \left( \frac{q^2}{1 - q^2} - \frac{q^4}{1-q^3} + 1 - q \right)\ZZ(q)
        \\
        &= \frac{1}{(1 + q)(1 - q^3)} \ZZ(q),
        \\
        T_{\scalebox{0.25}{\yng(3)}}.F\Bigg{|}_1 
        &=
        T_{\scalebox{0.25}{\yng(1,1,1)}}.F\Bigg{|}_1
        = \sum_{\ell \ge 0} \frac{\mathrm{d}}{\mathrm{d}y_{\ell,3}} \left(\prod_{k=1}^3 y_{0,k} \right)\prod_{j\ge 1}\left(y_{j,0} + \dots + y_{j,3} \frac{q^{3j}}{1 - q^j}\right)\Big{|}_{y=1} 
         \\
        &=\sum_{\ell \ge 1} \frac{q^{3\ell}}{1 - q^\ell} \prod_{\substack{j\ge 1\\j \neq \ell}} \frac{1}{1 - q^j} + \ZZ(q)
        \\
        &= \left(1 + \sum_{\ell \ge 1} q^{3\ell} \right) \ZZ(q) 
        \\
        &= \frac{1}{(1-q^3)} \ZZ(q) ,
        \\
        T_{\rotatebox[origin=c]{90}{\scalebox{0.25}{\yng(1,2)}}}.F\Bigg{|}_1 &= \sum_{\ell \ge 0} \frac{\mathrm{d}}{\mathrm{d}y_{\ell,1}} \left(\frac{\mathrm{d}}{\mathrm{d}y_{\ell+1,1}} - \frac{\mathrm{d}}{\mathrm{d}y_{\ell+1,2}}\right) \left(\prod_{k=1}^3 y_{0,k} \right) 
        \prod_{j\ge 1}\left(1 + y_{j,1} q^j + y_{j,1}y_{j,2}\frac{q^{2j}}{1 - q^j}\right)\Big{|}_{y=1} 
        \\
        &= \sum_{\ell \ge  1} \frac{q^{\ell}}{1 - q^{\ell}} \left(\frac{q^{\ell+1}}{1-q^{\ell+1}}\prod_{\substack{ j\ge 1\\j \neq \ell, \ell+1}}\frac{1}{1-q^j} - \frac{q^{2\ell+2}}{1-q^{\ell+2}} \prod_{\substack{ j\ge 1\\j \neq \ell, \ell+2}}\frac{1}{1-q^j} \right) + (q-q^2)\ZZ(q)
        \\
        &= \left(\frac{q^3}{1-q^2} - \frac{q^5}{1-q^3} + q-q^2 \right)\ZZ(q)
        \\
        &= \frac{q}{(1+q)(1-q^3)}\ZZ(q).
    \end{align*}
    By summing up all contributions together one obtains
    \begin{align*}
        \FFZ_3(q) 
        &= \Big[ \frac{q^3}{(1-q)(1-q^2)(1-q^3)} + 2\frac{1}{1-q^3} + 2 \left(\frac{1}{(1-q)(1-q^2)} - \frac{1}{(1-q^3)}\right) 
        \\
        & \qquad \qquad \qquad + \frac{q}{(1+q)(1-q^3)} + \frac{1}{(1+q)(1-q^3)} \Big] \ZZ(q)
        \\ 
        &= \frac{3 - q - q^2}{(1-q)(1-q^2)(1-q^3)} \ZZ(q).
        \nonumber
    \end{align*}
\end{example}

\begin{remark}
The connected cases can alternatively be computed directly employing the chain of equalities leading to equation \eqref{eqn:bagascia} in the proof of Proposition \ref{prop:Zlambda}.
\end{remark}

\subsubsection{\texorpdfstring{Some values of the polynomials $\mathsf P_\lambda$ and $\mathsf P_D$}{}}
\label{sec:poly-values}

For a skew Ferrers diagram $\lambda$, we write $\mathsf P_{\lambda}(q) = \widetilde{\mathsf P}_{\lambda}(q) \cdot \prod_{i = 1}^{\max(L,V)-1}(1 - q^i)$ where $\widetilde{\mathsf P}_{\lambda}(q)$ is the polynomial occurring in \Cref{prop:Zlambda}. Let us compute explicitly some polynomials $\mathsf P_{\lambda}(q)$ for special cases of skew Ferrers diagrams $\lambda$ and any size $D$.

\begin{lemma}\label{lemma:specialcases} For arbitrary positive integers $D$, consider $\lambda$ equal to the one-part partition $\lambda = (D)$ and its transpose $\lambda = (1,1, \dots, 1)$, as well as the disconnected diagram $\lambda = (1),(1), \dots, (1)$. Then we have
\begin{equation}
\frac{\FFZ_{\scalebox{0.25}{\yng(1,1,1,1)}}}{\ZZ(q)} = \frac{\FFZ_{\scalebox{0.25}{\yng(4)}}}{\ZZ(q)} = \frac{1}{1-q^D},
\qquad \text{ or equivalently } \qquad 
\mathsf P_{\scalebox{0.25}{\yng(1,1,1,1)}}(q) = \prod_{i=1}^{D-1} (1 - q^i),
\end{equation}
and 
\begin{equation}
\frac{\FFZ_{\scalebox{0.25}{\yng(1)}, \, \scalebox{0.25}{\yng(1)}, \, \dots, \scalebox{0.25}{\yng(1)}}}{\ZZ(q)} = \frac{q^{\binom{D}{2}}}{\prod_{i=1}^{D}(1-q^i)}, 
\qquad \text{ or equivalently } \qquad 
\mathsf P_{\scalebox{0.25}{\yng(1)}, \, \scalebox{0.25}{\yng(1)}, \, \dots, \scalebox{0.25}{\yng(1)}}(q) = q^{\binom{D}{2}}.
\end{equation}
\end{lemma}

\begin{proof}
The first equation immediately follows from the proof of \Cref{thm:ZD}, in particular specialising \Cref{eq:applicationTlambda} for $M=1$, $\ell_1 = 1$ and $v_1 = D$. The second equation is obtained by generalising the first term of \Cref{ex:D=3} to $D$ boxes, which leads to 
\[
\frac{\FFZ_{\scalebox{0.25}{\yng(1)}, \, \scalebox{0.25}{\yng(1)}, \, \dots, \scalebox{0.25}{\yng(1)}}}{\ZZ(q)}
=
\left(
\sum_{0 \leqslant k_1 < k_2 < \dots < k_D} q^{k_1 + k_2 + k_3 + \dots k_D}
\right)
= \frac{q^{\binom{D}{2}}}{\prod_{i=1}^D(1 - q^i)}.
\qedhere
\]
\end{proof}

\begin{remark}
\label{rem:0110}
Notice that both of the polynomials $\mathsf P_{\lambda}(q)$ in \Cref{lemma:specialcases} are of degree $\binom{D}{2}$. However, we have
    \begin{align*}
        \mathsf P_{\scalebox{0.25}{\yng(4)}}(0) &= 1, 
        \qquad \qquad 
        \quad \mathsf P_{\scalebox{0.25}{\yng(4)}}(1) = 0, 
        \\
        \mathsf P_{\scalebox{0.25}{\yng(1)}, \, \scalebox{0.25}{\yng(1)}, \, \dots, \scalebox{0.25}{\yng(1)}}(0) &= 0, 
        \qquad \qquad 
        \mathsf P_{\scalebox{0.25}{\yng(1)}, \, \scalebox{0.25}{\yng(1)}, \, \dots, \scalebox{0.25}{\yng(1)}}(1) = 1.
    \end{align*}
\end{remark}

In fact, behind this phenomenon, there is some structure, which we describe in the following two lemmas.
 
\begin{lemma} \label{lem:Plambda0110} 
For an arbitrary skew Ferrer diagram $\lambda$ of size $D$ we have
\begin{align*}
    \mathsf P_{\lambda}(0) &= \delta_{\lambda \in \mathrm{P}},
        \\
    \mathsf P_{\lambda}(1) &= \delta_{\lambda = \scalebox{0.25}{\yng(1)}, \, \scalebox{0.25}{\yng(1)}, \, \dots, \scalebox{0.25}{\yng(1)}},
\end{align*}
where we employ the notation for Kronecker deltas as in \Cref{sec:conventions}.

\begin{proof}
From \Cref{prop:Zlambda} the value $\mathsf P_{\lambda}(0)$ enumerates the ways of adding the skew Ferrers diagram $\lambda$ to the empty partition $\rho = (0)$. Clearly, this can be done if and only if $\lambda$ is already a partition, and then in exactly one way.

For the second equation, by \Cref{rem:0110} it is enough to show that if $\lambda$ is not the disjoint union of $D$ boxes, then $\mathsf P_{\lambda}(1)$ must be zero. This is the same as proving that every connected component of size greater than 1 provides a multiplicative factor of $(1 - q^i)$ for some $i$ to the polynomial $\mathsf P_{\lambda}(q)$. In fact it follows from the structure of \Cref{eq:applicationTlambda} that any connected component with horizontal total length $L$ or vertical total length $V$ greater than one produces such a factor. Therefore the only connected components possibly not vanishing at $q=1$ must have $L=M=1$, that is, they are diagrams of a single box. This concludes the proof of the lemma.
\end{proof}
\end{lemma}

\begin{lemma}
\label{lemma:specialisation}
For an arbitrary positive number $D$ we have
\begin{align}
    \mathsf P_D(0) &= p_D, \notag \\
    \mathsf P_D(1) &= 1 \label{eqn:pD=1}.
\end{align}
Here $p_D$ is the number of partitions of size $D$.
\begin{proof}
    It follows from \Cref{lem:Plambda0110} by summing over all skew Ferrers diagrams of size $D$.
\end{proof}
\end{lemma}

\begin{corollary}
\label{cor:pole-in-1}
The ratio $\FFZ_D(q) / \ZZ(q)$ has a pole at $q=1$ of order precisely $D$.
\begin{proof}
Combine \Cref{eqn:pD=1} and \Cref{eqn:ratio} with one another.
\end{proof}
\end{corollary}

\begin{lemma}
\label{lem:firstcoefficient}
The first coefficients of the polynomial $\mathsf{P}_D$ read
\begin{align*} 
 \operatorname{Coeff}_{q^1} \mathsf{P}_D(q) &= p_{D+1} - 2p_D, \quad D \geq 1,
 \\
 \operatorname{Coeff}_{q^2} \mathsf{P}_D(q) &= 2p_{D+2} - 2p_{D+1} - p_D - 2, \quad D \geq 2.
  \\
 \operatorname{Coeff}_{q^3} \mathsf{P}_D(q) &= 3p_{D+3} -4p_{D+2} - p_{D+1} + 2p_D - 2 - D - \delta_{D \equiv 1 \pmod 2}, \quad D \geq 3.
\end{align*}
Here $p_D$ is the number of partitions of size $D$.
\begin{proof}
    Notice that expanding in Taylor series one finds $\operatorname{Coeff}_q\FFZ_D(q) = 2P_D(0) + \frac{\mathrm d\mathsf P_D}{\mathrm dq}(0)$ and that $\operatorname{Coeff}_q\FFZ_D(q)$ counts how many skew Ferrers diagrams of $D$ boxes there are such that, once attached to the partition $\lambda = (1)$, give a partition of $D+1$ boxes. This is equivalent to the fact that $\operatorname{Coeff}_q\FFZ_D(q) = p_{D+1}$. Using \Cref{lemma:specialisation} concludes the proof of the first equation. The second and the third equations are computed with the same method using the previous ones.
\end{proof}
\end{lemma}

\subsection{The series \texorpdfstring{$\FFZ_{\bk}(q)$}{} for arbitrary \texorpdfstring{$\bk\in\BZ_{\geqslant 0}^s$}{}}
\label{subsec:g.f.tantinest}

In this section we prove \Cref{thm:intro:Z_k} in its general form. First we define a generalisation of reverse plane partition over a skew Ferrers diagram \textemdash\; their enumeration plays an important role in the generating series of the Euler characteristic of nested Hilbert schemes. See also \cite{double-nested-1,Mon_double_nested,CONSTELLATION} for recent works relating Hilbert schemes to the combinatorics of Young fillings.

Recall that we work with the partial order on $\BN^2$ given by
\begin{equation}
\label{eqn:PO}
(a,b) \geq (c,d) \iff a \geq c \mbox{ and }b\geq d.
\end{equation}

\begin{definition}
Fix a vector $\bk=(k_1,\ldots,k_s) \in \BZ_{\geq 0}^s$ and a skew Ferrers diagram $\lambda \in \mathrm Q$. A $\bk$-\emph{reverse plane partition} $\varphi$ \emph{of shape} $\lambda$ is a weakly decreasing function
    \[
    \begin{tikzcd}
\lambda\arrow[r,"\varphi"]&\Set{1,\ldots,s}
    \end{tikzcd}
    \]  
with respect to the partial order \eqref{eqn:PO}, such that
\[
\lvert \varphi^{-1} (i)\rvert =k_i,
\]
for $i=1,\ldots,s$. Similarly one defines the analogue up to translation for $\lambda\in\mathscr Q$.
\end{definition}

\begin{remark}
\label{rmk:madonna-spermatica}
Note that a $\bk$-reverse plane partition can be understood as a labelling of the boxes of a skew Ferrers diagram with the property that one can subsequently  remove the boxes labelled with the same number in decreasing order, while keeping the property of being a skew Ferrers diagram. Note also that the sum $\sum_{i=1}^sk_i$ agrees with the size $\lvert \lambda\rvert$, i.e.~with its cardinality. In this section we denote this number by $K$.
\end{remark}

\begin{definition}
Let $\lambda\in \mathscr Q $ be a skew Ferrers diagram, and let $\bk=(k_1,\ldots,k_s) \in \BZ_{\geq 0}^s$ be a vector of nonnegative integers. We denote by $\mathsf{RP}(\lambda;\bk)$ the number of $\bk$-reverse plane partition $\varphi$ of shape $\lambda$.
\end{definition}

 For a skew Ferrers diagram $\lambda\in\mathscr Q$, the number $\mathsf{RP}(\lambda;1, \dots,1)$  of $\bk$-reverse plane partitions of shape $\lambda$ is given by the classical hook length formula (see \cite{hook-length-formula} or \cite[Thm.~H, p.~60]{zbMATH03473265}), whereas $\mathsf{RP}(\lambda;K) = 1$ for any shape $\lambda$. See also \cite{Gansner_reversed_plane_partitions} for other results on the enumeration of reverse plane partitions.

\begin{figure}[H]
\ytableausetup{boxsize=1.5em}
\begin{ytableau}
1 &2&4&7&8\\
3&5&6&9 \\
10
\end{ytableau} 
\label{fig:Young_tableau}
\caption{A $(1,\ldots,1)$-reverse plane partition of shape $\lambda = (5,4,1)$. Removing simultaneously all boxes labelled $\{10, 9, \dots, i\}$, for any $i=10, \dots, 1$ gives again an underlying Ferrers diagram.}
\end{figure}

\begin{figure}[H]
\label{fig:tuamadre}
\begin{ytableau}
1 &3&4\\
2&3&4 \\
3
\end{ytableau} 
$\quad\qquad \qquad \qquad$
\begin{tikzpicture}
    \draw (0,0)--(1.5em,0)--(1.5em,1.5em)--(3em,1.5em)--(3em,3em)--(4.5em,3em)--(4.5em,6em)--(3em,6em)--(3em,4.5em)--(0em,4.5em)--cycle;
    \draw (0,1.5em)--(1.5em,1.5em)--(1.5em,4.5em) (0,3em)--(3em,3em)--(3em,4.5em) --(4.5em,4.5em) ;

    \draw  (0,0)--(0,4.5em)--(3em,4.5em)--(3em,6em)--(4.5em,6em);
    
    \node at (0.75em,0.75em) {$3$};
    \node at (0.75em,2.25em) {$2$};
    \node at (0.75em,3.75em) {$1$};
    \node at (2.25em,2.25em) {$3$};
    \node at (2.25em,3.75em) {$3$};
    \node at (3.75em,3.75em) {$3$};
    \node at (3.75em,5.25em) {$2$};
\end{tikzpicture} 
\caption{On the left, a $(1,1,3,2)$-reverse plane partition of shape  $\lambda=(3,3,1)$. 
On the right, a  $(1,2,4)$-reverse plane partition with shape a proper skew Ferrers diagram.}
\end{figure} 

\begin{prop}\label{prop:YTs}
Let $s>0$ be a positive integer. Given a sequence of nonnegative integers $\bk=(k_1,\ldots,k_s)\in \BZ^s_{\geqslant 0}$, the generating series $\FFZ_{\bk} (q)$ can be obtained as  
\begin{equation}
    \FFZ_{\bk}(q) = \!\!\!\! \sum_{\lambda\in \mathscr{Q}^{[K]}} \,\mathsf{RP}(\lambda; \bk)  \cdot T_{\lambda}.F^{[K]}(\vec{y};q) \Big{|}_1.
\end{equation}
\end{prop}

\begin{proof}
It follows from the definition of ${ \mathsf{RP}}(\lambda; \bk)$, of the operator $T_{\lambda}$ in \Cref{eq:defTlambda} and of $\chi^{[n,n+k_1,n+k_1+k_2,\ldots,n+\sum_{i=1}^sk_i]}$ and finally from \Cref{thm:skewferret}: each addition of the skew Ferrers diagram $\lambda$ performed by $T_{\lambda}$ is weighted by the number of ways it could be added into $s$ steps so that adding the first $k_1$ boxes of $\lambda$, followed by other $k_2$ boxes of $\lambda$, \dots, up until the last $k_s$ boxes of $\lambda$. At each intermediate step of adding $k_i$ boxes the resulting diagram is indeed the Ferrers diagram of a partition (i.e. the definition of partition is not broken by the newly added boxes at the $i$-th step for $i=1, \dots, s$). 
\end{proof}

\begin{corollary}\label{cor:nested:skew}
Let $s>0$ be a positive integer. Given a sequence of nonnegative integers $\bk=(k_1,\ldots,k_s)\in \BZ^s_{\geqslant 0}$, the generating series $\FFZ_{\bk}(q)$ can be obtained as
\begin{equation}
    \FFZ_{\bk}(q) = \!\!\!\! \sum_{\lambda\in \mathscr{Q}^{[K]}} \!\!\!\!\!\! { \mathsf{RP}}(\lambda; \bk) \frac{\mathsf P_{\lambda}(q)}{\prod_{i=1}^{L_{\NW}(\lambda) - 1} (1 - q^{i})} \ZZ(q),
\end{equation} 
where the polynomials $\mathsf P_{\lambda}$ are the same as in \Cref{prop:Zlambda} multiplied by $\prod_{i=1}^{\max(L,V)-1}(1-q^i)$. They are divisible by $q^B$ and of degree at most $\binom{L_{\NW}-1}{2} + \binom{\max(L,V)}{2} + B $, for $B = \sum_{i=1}^M \ell_i \left( \sum_{k=1}^{i-1} v_k \right)$.
\begin{proof}
It follows from \Cref{prop:YTs} and the proof of \Cref{thm:ZD}.
\end{proof}
\end{corollary}

\begin{theorem}\label{thm:ZDnested}
Let $s>0$ be a positive integer. Given a sequence of nonnegative integers $\bk=(k_1,\ldots,k_s)\in \BZ^s_{\geqslant 0}$ summing up to $K$, the generating series $\FFZ_{\bk}(q)$ satisfies the relation
\begin{equation} 
    \frac{\FFZ_{\bk}(q)}{\ZZ(q)} = \frac{\PP_{\bk}(q)}{\prod_{j=1}^{K} (1 - q^j)},
\end{equation}
where $\PP_{\bk}(q)$ is a polynomial in $q$ of degree bounded by $\frac{5}{4}K^2 - 2K + 1$. In particular, the ratio $\FFZ_{\bk}/\ZZ$ is a rational function in $q$ with only roots of unity as poles.
\begin{proof}
Rationality and the type of factors of the denominator \textemdash\; though not their multiplicities \textemdash\; follow from \Cref{cor:nested:skew} and the proof of \Cref{thm:ZD}.
\end{proof}
\end{theorem}

\begin{lemma} We have
\label{rmk:specialisation}
    \[
    \mathsf P_{\bk}(0) = \sum_{\lambda \in \mathrm{P}^{[K]}} { \mathsf{RP}}(\lambda;\bk).
    \]
\end{lemma}
    \begin{proof}
    By \Cref{cor:nested:skew} we have
    \[
    \mathsf P_{\bk}(0) = \sum_{\lambda\in \mathscr{Q}^{[K]}}{ \mathsf{RP}}(\lambda;\bk)\mathsf P_{\lambda}(0)\mathsf Q_{\lambda}(0),
    \]
    where $\mathsf Q_{\lambda}(q)$ is the product of factors of the form $(1 - q^i)$ which arise by taking common denominator. Hence $\mathsf Q_{\lambda}(0) = 1$. By \Cref{lem:Plambda0110} we have that $\mathsf P_{\lambda}(0) = \delta_{\lambda \in \mathrm{P}}$, concluding the proof.
    \end{proof}


\begin{remark}\label{rmk:denominator-simplifications}
    Notice that many further simplifications of the denominator in \Cref{thm:ZDnested} could occur. For example, if $(1+q+\dots+q^N)$ divides some polynomial $\mathsf P_{\bk}(q) = \widetilde{\mathsf P}_{\bk}(q)\cdot(1+q+\dots+q^N)$ for some $N<K$ then
    $$
    \frac{\mathsf P_{\bk}(q)}{\prod_{i=1}^K(1-q^{i})} = \frac{\widetilde{\mathsf P}_{\bk}(q)}{(1-q)^2\prod_{i\neq 1,N+1}^K(1-q^{i})}.
    $$
    Notice that these simplifications decrease the degree of the polynomial at the numerator but leave invariant the total degree as a rational function. This phenomenon simplifies $\mathsf P_{\bk}$ and reduces the number of distinct factors at the denominator, though increasing their multiplicity. In fact we see from the numerics (see \Cref{sec:numerics}) that in the nested case and small $k_i$ this phenomenon shows up to a certain extent. It would be interesting to have a geometric explanation of the phenomenon and establish more control on how and when it occurs.
\end{remark}

\section{The Euler characteristic of the Quot scheme}
\label{sec:euler-quot}

The goal of this section is to extend the formalism and main results of Sections \ref{sec:g.f.}--\ref{sec:proof-B} to the higher rank setting, and to prove Theorems \ref{thm:intro-Z_Dr} and \ref{thm:Z_Dr3-intro}.

\subsection{Punctual nested Quot schemes and generating functions}
Fix integers $\ell$, $r>0$ and a nondecreasing $\ell$-tuple $\bn = (n_1,\ldots,n_\ell)$ of integers $n_i \in \BZ_{\geq 0}$. Let $S$ be a smooth quasiprojective surface. Consider the \emph{nested Quot scheme}
\[
\Quot_r^{[\bn]}(S)
\]
parametrising isomorphism classes of quotients
\[
\begin{tikzcd}
\OO_{S}^{\oplus r} \arrow[two heads]{r} & F_\ell\arrow[two heads]{r} & F_{\ell-1}\arrow[two heads]{r} & \cdots \arrow[two heads]{r} & F_1,
\end{tikzcd}
\]
where $F_i$ is a 0-dimensional coherent sheaf on $S$ such that $\chi(F_i) = n_i$ for $i=1,\ldots,\ell$. As in the case $r=1$, there is a Quot-to-Chow morphism $\mathsf{qc}\colon \Quot_r^{[\bn]}(S) \to \Sym^{n_\ell}(S)$ and we define, for a closed point $p \in S$, the \emph{punctual nested Quot scheme}
\[
\Quot_r^{[\bn]}(S)_p \subset \Quot_r^{[\bn]}(S)
\]
as the fibre of $\mathsf{qc}$ over $n_\ell \cdot p \in \Sym^{n_\ell}(S)$. This scheme does not depend on $(S,p)$, i.e.~there is an isomorphism $\Quot_r^{[\bn]}(S)_p \cong \Quot_r^{[\bn]}(\BA^2)_0$. We set
\begin{equation}
\label{eqn:chi_r}
\chi_r^{[\bn]} = \chi (\Quot_r^{[\bn]}(\BA^2)_0) = \chi (\Quot_r^{[\bn]}(\BA^2)).
\end{equation}
Clearly we have $\chi_1^{[\bn]} = \chi^{[\bn]} = \chi(S_p^{[\bn]})$. The second equality in \eqref{eqn:chi_r} follows from the fact that the fully punctual locus is invariant under the $\BG_{\mathrm{m}}$-action scaling the fibres of $\OO_S^{\oplus r}$ and contains the fixed locus of the full nested Quot scheme.

\begin{definition}
\label{quotseries}
Form the generating functions
\begin{align*}
    \QQ_r(q)&=\sum_{n\ge 0}\chi_r^{[n]}q^n\in\BZ\llbracket q\rrbracket, & \text{for } r\in\BZ_{\ge0},\\
    \QQ(q,s)&=\sum_{r\ge 0}\QQ_r(q)s^r\in\BZ\llbracket q,s\rrbracket,\\
    \FFQ_{r,\bk}(q)&=\sum_{n\ge 0}\chi_r^{[n,n+k_1,n+k_1+k_2,\ldots,n+\sum_{i=1}^sk_i]}q^n \in\BZ\llbracket q\rrbracket,& \text{for } r\in\BZ_{\ge 0},  \bk=(k_1,\ldots,k_\ell) \in \BZ_{\geq 0}^\ell,\\
    \FFQ_{\bk}(q,s)&=\sum_{r\ge 0}\FFQ_{r,\bk}(q)  s^r \in\BZ\llbracket q,s\rrbracket, & \text{for }  \bk=(k_1,\ldots,k_\ell) \in \BZ_{\geq 0}^\ell,\\ 
    \FFQ(q,s,\boldsymbol{\vec v} ) &=\sum_{\ell\ge 1}\sum_{\bk\in\BZ^{\ell}_{\ge0}}\FFQ_{\bk}(q,s) \boldsymbol{\vec v}^{\bk}\in\BZ\llbracket q,s,\boldsymbol{\vec v} \rrbracket,
\end{align*} 
where, we put $\boldsymbol{\vec{v}}=(v_i)_{i\ge1}$, and $\boldsymbol{\vec v}^{\bk}=\prod_{i=1}^\ell v_i^{k_i}$.
\end{definition}

In the unnested setup we have the following relation.

\begin{lemma}
\label{lemma:sZ-geometric}
There is an identity
\begin{equation}
\label{eqn:Q(qs)}
\QQ(q,s) = \frac{1}{1-s\mathsf Z(q)}.
\end{equation}
\end{lemma}

\begin{proof}
By Bifet's work \cite{Bifet} we have, for each $d>0$, a scheme isomorphism
\[
\Quot_{r}^{[n]}(\BA^d)_0^{\BG_{m}^r} \cong \coprod_{n_1+\cdots +n_r = n}\prod_{i=1}^r (\BA^d)_0^{[n_i]},
\]
where $\BG_{m}^r$ is the $r$-dimensional torus acting by scaling the fibres of the trivial bundle $\OO^{\oplus r}$ on $\BA^d$. Taking Euler characteristics for the case $d=2$ gives
\begin{equation}
\label{eqn:chi-quot}
\chi_r^{[n]} = \sum_{n_1+\cdots +n_r = n}\prod_{i=1}^r \chi^{[n_i]}.
\end{equation}
Now, summing over $n\geq 0$, we can use \Cref{eqn:chi_r} with $\ell=1$ to obtain the relation
\begin{equation}
\label{rank-r-to-rank-1}
\mathsf Q_r(q) = \mathsf Q_1(q)^r = \mathsf Z(q)^r,
\end{equation}
from which the sought after identity follows immediately (after observing that $\mathsf Q_0(q)=1$).
\end{proof}

The nested version of \Cref{lemma:sZ-geometric} will be proved in \Cref{thm:Z_Dr}.

We can now introduce the combinatorial objects that the integers $\chi_r^{[\bn]}$ are enumerating.

\begin{definition}
\label{def:r-colours}
Fix integers $r,n\in\BZ$ such that $r-1,n\ge 0$. An \textit{$r$-coloured Ferrers diagram of size $n$} is an $r$-tuple of Ferrers diagrams $\boldsymbol\lambda=(\lambda^{(1)},\ldots,\lambda^{(r)})\in\mathrm{P}^{[n_1]} \times \cdots \times \mathrm{P}^{[n_r]}$ such that the number $\lvert\boldsymbol\lambda\rvert= \lvert\lambda^{(1)}\rvert+\cdots+\lvert\lambda^{(r)}\rvert = n_1+\cdots+n_r$ is equal to $n$.
\end{definition}

Let $\boldsymbol\lambda,\boldsymbol\mu$ be $r$-coloured Ferrers diagrams. We write $\boldsymbol\lambda \subset \boldsymbol\mu$ to mean that $\lambda^{(i)} \subset \mu^{(i)}$ for every $i = 1,\ldots,r$. This, in particular, implies that $\lvert\boldsymbol\lambda\rvert \leq \lvert\boldsymbol\mu\rvert$ and allows us to define the nested version of \Cref{def:r-colours} as follows.

\begin{definition}
\label{def:nested-r-colours}
Fix $\ell > 0$ and a nondecreasing sequence $\boldsymbol n=(n_1,\dots,n_\ell)\in\BZ_{\ge 0}^\ell$. An \textit{$r$-coloured $\boldsymbol n$-nested Ferrers diagram} is a nesting $\boldsymbol\lambda_1\subset\cdots\subset\boldsymbol\lambda_\ell$ of $r$-coloured  Ferrers diagrams of respective sizes $|\boldsymbol{\lambda}_i|=n_i$, for $i=1,\ldots,\ell$. We denote by $\mathrm P^{[\boldsymbol n]}_r$ the set of all $r$-coloured $\boldsymbol n$-nested Ferrers diagrams.
\end{definition}

\begin{remark}
It is clear that a nested version of the identity \eqref{eqn:chi-quot} holds. Indeed, there is an identification
    \begin{equation}\label{eq:nestedQuot-hilb}
            \mathrm P_r^{[\boldsymbol n]}\simeq\coprod_{\boldsymbol n^{(1)}+\cdots+\boldsymbol n^{(r)}=\boldsymbol n}\prod_{i=1}^r\mathrm P^{[\boldsymbol n^{(i)}]},
    \end{equation}
where $\boldsymbol n^{(i)} \in \BN^\ell$ is a nondecreasing $\ell$-tuple of nonnegative integers, for all $i=1,\dots,r$. In particular, there is an identity $\chi_r^{[\boldsymbol n]}=\lvert\mathrm P_r^{[\boldsymbol n]}\rvert$.
\end{remark}

\subsection{\texorpdfstring{The proof of Theorems \ref{thm:intro-Z_Dr} and \ref{thm:Z_Dr3-intro}}{}}
\label{sec:C-D}

The generating function of the Euler characteristics of nested Quot schemes can be neatly expressed in terms the generating function of Euler characteristics of nested Hilbert schemes. We define the generating function by
\begin{equation*}
    \FFZ(q,\boldsymbol{\vec v}) = \sum_{\ell\ge1}\sum_{\bk\in\BZ^\ell_{\ge0}} \FFZ_{\bk}(q)\boldsymbol{\vec v}^{\bk}\in\BZ\llbracket q,\boldsymbol{\vec v}\rrbracket . 
\end{equation*}
{  Whenever we consider only nestings of length one, we write $\FFZ(q,v)$ for the specialisation $\FFZ(q,v)=\FFZ(q,\boldsymbol{\vec v})\big|_{{v_1=v},\,v_{>1} = 0}$. We adopt the same notation for the series in \Cref{quotseries} involving nestings of $r$-coloured Ferrers diagrams.

\begin{theorem}
\label{thm:Z_Dr}
There is an identity of formal power series
    \[
    \mathsf{FQ}(q,s,\boldsymbol{\vec v}) =
\frac{1}{ 1 - \mathsf{FZ}(q,\boldsymbol{\vec v})s}\in\BZ\llbracket q,s,\boldsymbol{\vec v}\rrbracket.
    \]
\end{theorem}

\begin{proof}    
We exploit first the nestings of length one and then we explain how to deduce the general formula. By definition of an $r$-coloured Ferrers diagram, we have an equality 
    \[
    \FFQ_{r,D}(q) = \sum_{\substack{D_1,\dots,D_r \ge 0 \\D_1+\cdots+D_r=D}}\prod_{i=1}^r \FFZ_{D_i}(q).
    \]
    Since $\FFZ_{D_i}(q)=\operatorname{Coeff}_{v_1^{D_i}}\FFZ(q,\boldsymbol{\vec v})$ we have that
    \begin{equation}
    \label{eqn:FQrd}
    \FFQ_{r,D}(q) = \operatorname{Coeff}_{v_1^D}\left(\FFZ(q,\boldsymbol{\vec v})^r\right).
    \end{equation}
    Similarly, one finds in the same way  the coefficient of $s^rv^{k_1}_1 v^{k_2}_2 \ldots v^{k_\ell }_\ell$.  Precisely, we have
    \[
    \operatorname{Coeff}_{s^rv^{k_1}_1 v^{k_2}_2 \ldots v^{k_\ell }_\ell}\left(\frac{1}{ 1 - \mathsf{FZ}(q,\boldsymbol{\vec v})s}\right)=\operatorname{Coeff}_{v^{k_1}_1 v^{k_2}_2 \ldots v^{k_\ell }_\ell}\left(\FFZ(q,\boldsymbol{\vec v})^r\right)= \operatorname{Coeff}_{s^rv^{k_1}_1 v^{k_2}_2 \ldots v^{k_\ell }_\ell}( \mathsf{FQ}(q,s,\boldsymbol{\vec v})),
    \]
    and the second equality is a consequence of \eqref{eq:nestedQuot-hilb}.

    Then, by  summing over  $r\ge0$, $\ell\ge 1$, and $\bk\in \BZ^\ell_{\ge 0}$, one obtains
    \[
    \FFQ(q,s,\boldsymbol{\vec v}) = \sum_{r \ge 0}\sum_{\ell\ge1}\sum_{\bk \in\BZ^\ell_{\ge 0}} \FFQ_{r,\bk}(q) \boldsymbol{\vec v}^{\bk} s^r= \sum_{r \ge 0}\FFZ(q,\boldsymbol{\vec v})^rs^r = \frac{1}{1 - \FFZ(q,\boldsymbol{\vec v})s}.\qedhere
    \]
\end{proof}
 
     


A straightforward consequence of \Cref{thm:Z_Dr} is that, thanks to \Cref{thm:ZD}, $\FFQ_{r,\bk}(q)$  is factorised as the product of a rational function by $\ZZ(q)^r$.

\begin{corollary}
\label{cor:FQ/Z^r}
For any choice of $r,D\in\BZ_{\ge0}$, there is an identity
\begin{equation}
\label{eqn:polynomial_Z_rD}
\frac{\FFQ_{r,D}(q)}{\QQ(q)}=\frac{\FFQ_{r,D}(q)}{\ZZ(q)^r}=\frac{\mathsf P_{r,D}(q)}{\prod_{j=1}^D(1-q^j)^{\min\left(r, \floor{\frac{D}{j}}\right)}}, \end{equation}
where $\mathsf P_{r,D}(q)\in\BZ[q]$ is a polynomial.
\end{corollary}

\begin{proof}
    Let $\mathsf G(q,v)$ be the formal power series
    \[
    \mathsf G(q,v)=\frac{\FFZ(q,v)}{\mathsf Z(q)}=\sum_{m\ge 0}\frac{\mathsf P_m(q)}{\prod_{i=1}^m(1-q^i)}v^m,
    \]
    where the existence of the polynomials $\mathsf P_m(q)$ on the right-hand side is provided by \Cref{thm:ZD}. By \Cref{eqn:FQrd}, we have that
    \begin{align*}
       \frac{\FFQ_{r,D}(q)}{\mathsf Z(q)^r}
       &=
       \operatorname{Coeff}_{v^D}\left(\mathsf G(q,v)^r\right)
       \\
       &=\sum_{\substack{d_1,\dots,d_r \ge 0 \\d_1+\cdots+d_r=D}}\prod_{j=1}^r\frac{\mathsf P_{d_j}(q)}{\prod_{i=1}^{d_j}(1-q^i)}
       \\
       &=\sum_{\substack{d_1,\dots,d_r \ge 0 \\d_1+\cdots+d_r=D}}\frac{ \prod_{j=1}^r \mathsf P_{d_j}(q)}{\prod_{i=1}^D(1-q^i)^{N_i(\boldit{d})}}
       \\
       & = \frac{\mathsf P_{r,D}(q) }{\prod_{i=1}^D(1-q^i)^{\max_{\boldit{d}} N_i(\boldit{d})}},
    \end{align*}
    where for each $1\le i\le D$ we have set $N_i(\boldit{d})=|\set{j|d_j\ge i}|$. In particular, $\sum_i N_i=D$. Notice that $N_i \le r$, and since $\sum_j d_j = D$ the constraint $N_i \le \floor{\frac{D}{i}}$ must hold. Therefore $N_i \le \min\left( r , \floor{\frac{D}{i}}\right)$.  
\end{proof}

\begin{corollary}
\label{cor:FQ/Z^r:multinested} 
For any choice of $r\in\BZ_{\ge0}$ and vector ${\bk} = (k_1, \dots, k_\ell)\in\BZ^\ell_{\ge 0}$, we have an identity
\begin{equation}
\label{eqn:polynomial_Z_r:multi}
\frac{\FFQ_{r,\bk}(q)}{\ZZ(q)^r} 
= \frac{P_{\bk}(q)}{\prod_{j=1}^K (1 - q^j)^{u_j}}
\in\BZ\llbracket q \rrbracket,
\end{equation}
 where $K = \sum_{i=1}^\ell k_i$, {$u_j = u_j(\bk)$}  are nonnegative integers, and $P_{\bk}(q)$ is a polynomial that depends on the vector $\bk$. 

\begin{proof}
{  \Cref{thm:Z_Dr} expresses the relation between the series $\FFQ_{r,\bk}(q)$ and $\FFZ(q,\boldsymbol{\vec v})$. To conclude, one can argue as in the proof of \Cref{cor:FQ/Z^r}.}
\end{proof}
\end{corollary}

\begin{remark}
    In general, as it is the case in \Cref{rmk:denominator-simplifications}, many simplifications can occur in the right-hand side of \Cref{eqn:polynomial_Z_rD}, and the degree of the denominator can be quite lower.  For instance, one has
    \begin{equation*} 
\begin{split}
\FFQ_{r,1}(q)&= r\FFZ_1(q)\mathsf Z(q)^{r-1},\\
\FFQ_{r,2}(q)&= r\FFZ_2(q)\mathsf Z(q)^{r-1}+\binom{r}{2}\FFZ_1(q)^2\mathsf Z(q)^{r-2}.
\end{split}    
\end{equation*}
\end{remark}

Let us consider specialisation
\[
\FFQ(q,s,v) = \FFQ(q,s,\boldsymbol{\vec v})\big|_{{v_1=v},\,v_{>1} = 0}.
\]
We now confirm that $\FFQ(q,s,v)$ is entirely determined by the application of certain differential operators to the series $\QQ(q,s)$. This will prove \Cref{thm:Z_Dr3-intro} from the introduction.

\begin{theorem}
\label{thm:Z_Dr3} 
There is an identity of formal power series
    \begin{equation}\label{eq:thm:Z_Dr3}
        \FFQ(q,s,v ) =\exp{\left(\FFZ(q,v)-1\right)}\Big{|}_{v^k \mapsto v^ks^k\frac{\dd^k}{\dd s^k}} 
        .\QQ(q,s)
        \in\BZ\llbracket q,s,v\rrbracket,
    \end{equation}
where $\FFZ(q,v) = \FFZ(q,\boldsymbol{\vec v})\big|_{{v_1=v},\,v_{>1} = 0}$ and the substitution of variables is meant after the expansion of the exponential as a power series in $v$.
\end{theorem}

\begin{proof}
    The result in \Cref{thm:Z_Dr} allows us to compute $\FFQ_D(q,s)$ as the $D$-th coefficient in $v$ of $\FFQ(q,v,s)$, i.e.
    \[
    \FFQ_D(q,s)=\sum_{k\ge 0}s^k\operatorname{Coeff}_{v^D}\FFZ(q,v)^k.
    \]
    Thus we have
    \begin{align*}
        \FFQ_{D}(q,s)&=\sum_{k\ge 0}s^k\ZZ(q)^k\operatorname{Coeff}_{v^D}\left(\frac{\FFZ(q,v)}{\ZZ(q)}\right)^k\\
        &=\sum_{k\ge 0}s^k\ZZ(q)^k\sum_{\substack{(d_1,\dots,d_k)\in\BZ_{\ge 0}^k\\d_1+\cdots+d_k=D}}\prod_{i=1}^k\frac{\FFZ_{d_i}(q)}{\ZZ(q)}\\
        &=\sum_{\lambda\vdash D}s^{l(\lambda)}\ZZ(q)^{l(\lambda)}\prod_{i=1}^{l(\lambda)}\left(\frac{\FFZ_{\lambda_i}(q)}{\ZZ(q)}\right)\sum_{p\ge 0}\frac{(l(\lambda)+p)!}{p!\prod_{j=1}^D m_j(\lambda)!}s^p\ZZ(q)^p\\
        &=\sum_{\lambda\vdash D}\frac{l(\lambda)!}{\prod_{j=1}^D m_j(\lambda)!}\prod_{i=1}^{l(\lambda)}\left(\frac{\FFZ_{\lambda_i}(q)}{\ZZ(q)}\right)\frac{s^{l(\lambda)}\ZZ(q)^{l(\lambda)}}{\left(1-s\ZZ(q)\right)^{l(\lambda)+1}}\\
        &=\sum_{\lambda\vdash D}\frac{1}{\prod_{j=1}^D m_j(\lambda)!}\prod_{i=1}^{l(\lambda)}\left(\frac{\FFZ_{\lambda_i}(q)}{\ZZ(q)}\right)s^{l(\lambda)}\frac{\dd^{l(\lambda)}}{\dd s^{l(\lambda)}}\QQ(q,s)\\
        &=\left(\sum_{\lambda\vdash D}\mathscr D_\lambda\right)\QQ(q,s),
    \end{align*}
    where we used \Cref{lemma:sZ-geometric} in the second-to-last identity and in the last line we defined the differential operator $\mathscr D_\lambda$, acting on formal power series in $\BZ\llbracket q,s\rrbracket$, as
    \[
    \mathscr D_\lambda \defeq \frac{1}{\prod_{j\ge 1}m_j(\lambda)!}\prod_{i=1}^{l(\lambda)}\left(\frac{\FFZ_{\lambda_i}(q)}{\ZZ(q)}\right)s^{l(\lambda)}\frac{\dd^{l(\lambda)}}{\dd s^{l(\lambda)}}.
    \]
    Thus, summing over $D\ge 0$, we get
    \begin{align*}
        \FFQ(q,v,s)&=\sum_{D\ge 0}\FFQ_D(q,s)v^D\\
        &=\left(\sum_{D\ge 0}v^D\sum_{\lambda\vdash D}\mathscr D_\lambda\right)\QQ(q,s)\\
        &=\exp\left(\sum_{n\ge 1}\frac{\FFZ_{n}(q)}{\ZZ(q)}v^n\right)\Bigg{|}_{v^k \mapsto v^ks^k\frac{\dd^k}{\dd s^k}} \QQ(q,s),
    \end{align*}
    where the substitution of variable is meant to be made after the expansion of the exponential as a power series in $v$.
\end{proof}

\begin{remark}
In light of \Cref{thm:Z_Dr3}, the expansion of the first orders of \Cref{eq:thm:Z_Dr3} in the $v$-variable is:
\begin{align*}
\operatorname{Coeff}_{v^0}\FFQ(q,s,v) = \FFQ_0(q,s)&=1. \QQ(q,s)
\\
\operatorname{Coeff}_{v^1}\FFQ(q,s,v) = \FFQ_1(q,s)&=\left[\left(\frac{s}{1-q}\right)\frac{\mathrm d}{\mathrm ds}\right]. \QQ(q,s)
\\
\operatorname{Coeff}_{v^2}\FFQ(q,s,v) = \FFQ_2(q,s)&=\left[\frac{\mathsf P_2(q)}{(1-q)(1-q^2)}\left(s\frac{\mathrm d}{\mathrm d s}\right)+\frac{\mathsf P_1(q)^2}{(1-q)^2}\left(\frac{s^2}{2}\frac{\mathrm d^2}{\mathrm ds^2}\right)\right]. \QQ(q,s).
\end{align*}  
\end{remark}

\begin{remark}
\label{rmk:old-stuff}
Define
\[
\FFQ^{[\bullet]}(q,s,v)=\sum_{n\ge 0}\sum_{r\ge 0}\sum_{k\ge 0} \chi_r^{[n,k]}v^ks^rq^n.
\]
This series can also be expressed as
\[
\FFQ^{[\bullet]}(q,s,v)=  \FFQ(qv,s,v).
\]
For low values of $n$, the series 
\[
\FFQ^{[n]}(s,v)=\sum_{r\ge 0}\sum_{k\ge 0} \chi_r^{[n,k]}v^ks^r
\]
can be computed directly. For instance, one can show the following identity:
\begin{equation}
\label{madonna-luridissima-sempre}
\FFQ^{[2]}(q,s)=\left(s^2-\frac{2s}{1-q}+\left(2s(1-s+s^2)-\frac{2s^2}{1-q}\right)\frac{\dd}{\dd s}+\frac{s^2(1-s)^2}{2}\frac{\dd^2}{\dd s^2}\right)\QQ(q,s).
\end{equation}
\end{remark}

\appendix

\section{Some explicit examples}
In this subsection we present some explicit computations. For one nesting and two gaps, the base cases can be computed via \eqref{eq:gottsche} and  \Cref{thm:motivepunctual,motive-punctual-3n} and are
\begin{equation}
\label{eqn:1-step-basic}
 \begin{split}
    [S_p^{[0,2]}] &= \BP^1, \\
    [S_p^{[1,3]}] &= \BP^2, \\
    [S_p^{[2,4]}] &= 2\BL^3+3\BL^2+2\BL+1 ,\\
    [S_p^{[3,5]}] &= 2\BL^4 + 4\BL^3 + 4\BL^2 + 2\BL + 1.
\end{split}    
\end{equation}
For two nestings and one gap, the base cases are
\begin{equation}
\label{eqn:2-step-basic}
 \begin{split}
    [S_p^{[0,1,2]}] &=  [S_p^{[1,2]}]=[S_p^{[2]}]=[\BP^1] ,\\
    [S_p^{[1,2,3]}] &=  [S_p^{[2,3]}]=[\BP^1]^2 ,\\
    [S_p^{[2,3,4]}] &= [\BP^1](1+2\BL+2\BL^2)  .
\end{split}    
\end{equation}
The first two equalities are obtained using the formulas in \cite{Gottsche-motivic}. The last equality is obtained via the stratification
\[
\nested{2,3,4}_0= U_1\amalg U_2,
\]
where
\[
U_i=\Set{[Z_2,Z_3,Z_4]\in \nested{2,3,4}_0 | \bh_{Z_3}(1)=i}.
\]

One can also easily compute some of the global motives, e.g.~one finds 
\begin{equation}
\label{eqn:global-basic}
 \begin{split}
(\BA^2)^{[0,2]} &=  [\BP^1]\BL^3 ,\\
    (\BA^2)^{[1,3]} &=  (3\BL^4 + 2\BL^3 - \BL - 1 )\BL^2.
\end{split}    
\end{equation}
In order to compute the motive of 
$\nested{1,3}$ we first consider the composition
\[ 
\begin{tikzcd}
\nested{1,3} \arrow[swap]{dr}{\theta}\arrow{r}{\pr_2} & \nested{3} \arrow{d} \\
&\Sym^3(\BA^2),
\end{tikzcd}
\]
and we pullback along $\theta$ the stratification of the symmetric product by partitions $\Sym^3(\BA^2) = \coprod_{\lvert \lambda \rvert = 3}\Sym^3_\lambda(\BA^2)$ to obtain the stratification 
\[
\nested{1,3} =\coprod _{\lambda \in \mathrm{P}^{[3]}} H_{\lambda}=H_{\scalebox{0.2}{\yng(3)}}\amalg H_{\scalebox{0.2}{\yng(2,1)}}\amalg H_{\scalebox{0.2}{\yng(1,1,1)}},
\]
where, explicitly, we have
\[
H_\lambda=\Set{ [Z_1,Z_2]\in\nested{1,3} | \theta(Z_2) \in \Sym^3_\lambda(\BA^2)}.
\]
Now, we compute
\begin{align*}
    [H_{\scalebox{0.2}{\yng(3)}}]= & 3[\Sym^3(\BA^2)\smallsetminus \Delta]=3(\BL -  1)[\BP^1] \BL^4 \\
    [H_{\scalebox{0.2}{\yng(2,1)}}]=& 2[(\BA^2\times\BA^2\smallsetminus\BA^2)\times \BP^1]=2(\BL-1)[\BP^1]^2\BL^2 \\
    [H_{\scalebox{0.2}{\yng(1,1,1)}}] =& [(\BA^2)^{[1,3]}_0\times \BA^2]=\BL^2[\BP^2] ,  
\end{align*} 
and the sum gives the claimed formula displayed in \Cref{eqn:global-basic}. 

\section{Arbitrary surfaces}
\label{sec:power-structure}
Let $S$ be a smooth quasiprojective surface, and fix a length $\ell > 0$ for the nestings. Consider the motivic generating function
\[
\mathsf{Quot}_{S,r}(\boldit{q}) 
= \sum_{\bn} \,[\Quot^{[\bn]}_r(S)] \boldit{q}^{\bn} \,\in\,K_0(\Var_{\BC})\llbracket \boldit{q} \rrbracket,
\]
where we use the multi-index notation $\boldit{q}^{\bn} = q_1^{n_1}\cdots q_\ell^{n_\ell}$. Note that the coefficients are `global' (i.e.~nonpunctual) motives now. The power structure over $K_0(\Var_{\BC})$ gives the relation
\begin{equation}
\label{PW-quot}
\mathsf{Quot}_{S,r}(\boldit{q}) 
= \left(\sum_{\bn} \,[\Quot^{[\bn]}_r(\BA^2)_0] \boldit{q}^{\bn} \right)^{[S]}.
\end{equation}
See \cite{zbMATH05493515} for the case $r=1$ (in arbitrary dimension) and \cite{MR_nested_Quot, double-nested-1} for technical details on multivariable power structures. See also \cite[Thm.~C]{MOTIVES} for several explicit formulas for motives of punctual Quot schemes in all dimensions, all ranks and low number of points. As the Euler characteristic $\chi\colon K_0(\Var_{\BC}) \to \BZ$ is a homomorphism of rings with power structure, \Cref{PW-quot} induces an identity
\[
\chi\mathsf{Quot}_{S,r}(\boldit{q}) 
= \sum_{\bn}\chi(\Quot^{[\bn]}_r(S)) \boldit{q}^{\bn} = \left( \sum_{\bn} \chi_r^{[\bn]}\boldit{q}^{\bn}\right)^{\chi(S)}.
\]
{ Now fix $\ell=2$, and note that the series in big round brackets in the last identity is nothing but 
\[
\operatorname{Coeff}_{s^r}\FFQ(q,s,v), 
\]
and by our result this series is completely accessible. Furthermore, the motive of the Quot scheme only depends on the rank of the bundle we are taking quotients of, in the sense that replacing $\OO_S^{\oplus r}$ with any locally free sheaf of rank $r$ gives the same class in $K_0(\Var_{\BC})$, see \cite{ricolfi2019motive}. Therefore, to sum up, our results give complete access to the generating function of Euler characteristics
\[
\sum_{\bn} \chi(\Quot^{[\bn]}(E)) \boldit{q}^{\bn},
\]
where $E$ is an arbitrary locally free sheaf of rank $r$ on $S$. One can of course bootstrap from the single nesting ($\ell=2$)  to a nesting of  arbitrary length, just as we did in the rank 1 case, and the same power structure  machinery applies verbatim.}

\section{Some generating functions}
\label{sec:numerics}
In this section we list some of the polynomials $\mathsf P_D(q)$ and some of the ratios $\FFZ_{\bk}(q)/\ZZ(q)$ of \Cref{thm:ZD,thm:ZDnested}. We have computed the required initial data using the computer software sagemath \cite{sagemath}.
    \subsection{Some polynomials \texorpdfstring{$\PP_{D}(q)$}{}}
    \label{sec:some-polynomials}
    \begin{itemize}
        \item[$D=1$:] $1$
        \item[$D=2$:] $2-q$
        \item[$D=3$:] $3-q-q^2$
        \item[$D=4$:] $5-3q+q^2-2q^3-q^4+q^5$
        \item[$D=5$:] $7-3q-q^2+q^3-2q^4-5q^5+3q^6+q^7-q^8+2q^9-q^{10}$
        \item[$D=6$:]  $11-7q+q^2+q^3+q^4-11q^5+3q^6-2q^7+2q^8+q^{10}+4q^{11}-4q^{12}+2q^{13}-q^{14}$
        \item[$D=7$:]  $ 15-8q-q^2+4q^3-7q^5-3q^6-14q^7+12q^8+2q^9-9q^{10}+7q^{11}+5q^{12}+q^{13}+q^{14}-4q^{15}-q^{16}+5q^{17}-7 q^{18}+2q^{19}+2q^{20}-q^{21}$
        \item[$D=8$:]  
        $
        22 - 14q + 4q^3 + 11q^4 - 19q^5 + 6q^6 - 27q^7 + 7q^8 + 4q^9 - q^{10} - 13q^{11} + 15q^{12} + 4q^{13} + q^{14} + 13q^{15} - 8q^{16} - 3q^{17} + 6q^{18} + q^{19} - 15q^{20} + 5q^{21} + q^{22} + 4q^{23} - 7q^{24} + 3q^{25} + 3q^{26} - 2q^{27}
        $
        \item[$D=9$:]  
        $
        30 - 18q - 4q^2 + 13q^3 + 8q^4 - 16q^5 + 9q^6 - 33q^7 - q^8 - 6q^9 + q^{10} - 5q^{11} + 9q^{12} - 16q^{13} + 7q^{14} + 32q^{15} + 6q^{16} - 12q^{17} + 8q^{18} + 6q^{19} - 10q^{20} + 2q^{21} - 5q^{22} + q^{23} - 22q^{24} + 16q^{25} + 7q^{26} - 15q^{27} + 4q^{28} + 12q^{29} - 11q^{30} + 6q^{31} + q^{32} - 6q^{33} + 4q^{34} - t^{35}$
        \item[$D=10$:] 
        $
        42 - 28q - 2q^2 + 11q^3 + 23q^4 - 23q^5 + 24q^6 - 64q^7 + 25q^8 - 32q^9 - 7q^{10} - 6q^{11} + 38q^{12} - 76q^{13} + 23q^{14} + 31q^{15} + 13q^{16} + 8q^{17} + 23q^{18} - 7q^{19} + 16q^{20} - 8q^{21} - 15q^{22} + 47q^{23} - 47q^{24} - 26q^{25} + 15q^{26} + 12q^{27} - 33q^{28} + 24q^{29} - 19q^{30} + 19q^{31} - 5q^{32} - 3q^{33} + 25q^{34} - 7q^{35} - 28q^{36} + 20q^{37} + 9q^{38} - 9q^{39} - 6q^{40} - q^{41} + 8q^{42} - 3q^{43}.
        $
        \end{itemize} 
\subsection{Some ratios \texorpdfstring{$\FFZ_{\bk}(q)/\ZZ(q)$}{}}
\label{sec:some-ratios} 

        \begin{itemize}
        \item[$\bk=$]$(1,1)$:  $\displaystyle\frac{2}{(1-q)(1-q^2)}$
        \smallbreak
        \item[$\bk=$]$(1,1,1)$:  $\displaystyle\frac{4-2q}{(1-q)^2(1-q^2)}$
        \smallbreak
        \item[$\bk=$]$(1,1,1,1)$:  $\displaystyle\frac{10-4q-2q^2}{(1-q)^2(1-q^2)^2}$
        \smallbreak
        \item[$\bk=$]$(1,1,1,1,1)$:  $\displaystyle\frac{26-28q+6q^2}{(1-q)^3(1-q^2)^2}$
        \smallbreak
        \item[$\bk=$]$(1,1,1,1,1,1)$:  $\displaystyle\frac{76-72q-12q^2+16q^3}{(1-q)^3(1-q^2)^3}$
        \smallbreak
        \item[$\bk=$]$(1,2)$: $\displaystyle\frac{3+2q-q^2-q^3}{(1-q)(1-q^2)(1-q^3)}$ 
        \smallbreak
        \item[$\bk=$]$(2,1)$:  $\displaystyle\frac{ 4-q+2q^2-2q^3}{(1-q)(1-q^2)(1-q^3)}$
         \smallbreak
       \item[$\bk=$]$(2,2)$: $\displaystyle\frac{8-3q+8q^2-4q^3-2q^4-q^5}{(1-q)(1-q^2)(1-q^3)(1-q^4)}$.
    \end{itemize}


\bibliographystyle{amsplain-nodash}
\bibliography{The_Bible}

@Misc{Galkin-Shinder,
 Author = {Galkin, Sergey and Shinder, Evgeny},
 Title = {The {Fano} variety of lines and rationality problem for a cubic hypersurface},
 Year = {2014},
 HowPublished = {\href{https://arxiv.org/abs/1405.5154}{ArXiv:1405.5154}}
}

@Article{zbMATH02069674,
 Author = {Larsen, Michael and Lunts, Valery A.},
 Title = {Motivic measures and stable birational geometry},
 FJournal = {Moscow Mathematical Journal},
 Journal = {Mosc. Math. J.},
 ISSN = {1609-3321},
 Volume = {3},
 Number = {1},
 Pages = {85--95},
 Year = {2003}
}

@Article{zbMATH06837514,
 Author = {Borisov, Lev A.},
 Title = {The class of the affine line is a zero divisor in the {Grothendieck} ring},
 FJournal = {Journal of Algebraic Geometry},
 Journal = {J. Algebr. Geom.},
 ISSN = {1056-3911},
 Volume = {27},
 Number = {2},
 Pages = {203--209},
 Year = {2018}
}

@Article{Rasul-irr-nested,
 Author = {Gangopadhyay, Chandranandan and Rasul, Parvez and Sebastian, Ronnie},
 Title = {Irreducibility of some nested {Hilbert} schemes},
 FJournal = {Proceedings of the American Mathematical Society},
 Journal = {Proc. Am. Math. Soc.},
 ISSN = {0002-9939},
 Volume = {152},
 Number = {5},
 Pages = {1857--1870},
 Year = {2024}
}

@manual{sagemath,
  Key          = {SageMath},
  Author       = {{The Sage Developers}},
  Title        = {{S}ageMath, the {S}age {M}athematics {S}oftware {S}ystem ({V}ersion 9.3)},
  note         = {{\tt https://www.sagemath.org}},
  Year         = {2021},
}

@article{Gottsche-motivic,
	author = {G{\"{o}}ttsche, Lothar},
	journal = {{Math. Res. Lett.}},
	pages = {613--627},
	title = {{On the motive of the Hilbert scheme of points on a surface}},
	volume = {8},
	year = {2001}}

@article{8POINTS,
author = {Dustin Cartwright and Daniel Erman and Mauricio Velasco and Bianca Viray},
title = {{Hilbert schemes of 8 points}},
volume = {3},
journal = {Algebra \& Number Theory},
number = {7},
publisher = {MSP},
pages = {763 -- 795},
year = {2009},
doi = {10.2140/ant.2009.3.763},
URL = {https://doi.org/10.2140/ant.2009.3.763}
}

@misc{Fantechi-Ricolfi-structural, 
howpublished={In: Moduli, Motives and Bundles: New Trends in Algebraic Geometry. London Mathematical Society Lecture Note Series \textbf{499}, Cambridge University Press.}, 
title={On the stack of 0-Dimensional coherent sheaves: structural aspects}, 
author={Fantechi, Barbara and Ricolfi, Andrea T.}, 
year={2025}, 
volume={499},
pages={1--32}}

@article{MR-hyperquot,
 author = {Monavari, Sergej and Ricolfi, Andrea T.},
 title = {Hyperquot schemes on curves: virtual class and motivic invariants},
 fjournal = {Mathematische Annalen},
 journal = {Math. Ann.},
 issn = {0025-5831},
 volume = {392},
 number = {2},
 pages = {1665--1709},
 year = {2025}
}

@article{Fantechi-Ricolfi-motivic,
 author = {Fantechi, Barbara and Ricolfi, Andrea T.},
 title = {On the stack of 0-dimensional coherent sheaves: motivic aspects},
 fjournal = {Bulletin of the London Mathematical Society},
 journal = {Bull. Lond. Math. Soc.},
 issn = {0024-6093},
 volume = {57},
 number = {6},
 pages = {1607--1649},
 year = {2025},
 doi = {10.1112/blms.70096},
 zbMATH = {8059832}
}

@Article{Lissite-quot,
 Author = {Monavari, Sergej and Ricolfi, Andrea T.},
 Title = {{Sur la lissit\'e du sch\'ema Quot ponctuel emboît\'e}},
 FJournal = {Canadian Mathematical Bulletin},
 Journal = {Can. Math. Bull.},
 ISSN = {0008-4395},
 Volume = {66},
 Number = {1},
 Pages = {178--184},
 Year = {2023}
}

@article{MR_nested_Quot,
	author = {Monavari, Sergej and Ricolfi, Andrea T.},
	doi = {10.1016/j.jalgebra.2022.07.011},
	fjournal = {Journal of Algebra},
	issn = {0021-8693},
	journal = {J. Algebra},
	mrclass = {14C35}, 
	mrreviewer = {Remke Kloosterman},
	pages = {99--118},
	title = {On the motive of the nested {Q}uot scheme of points on a curve},
	url = {https://mathscinet.ams.org/mathscinet-getitem?mr=4463016},
	volume = {610},
	year = {2022}}

@misc{GMMR2,
	author = {Graffeo, Michele and Monavari, Sergej and Moschetti, Riccardo and Ricolfi, Andrea T.},
	howpublished = {{\href{https://arxiv.org/abs/2501.10267}{ArXiv:2501.10267}}},
	title = {{Enumeration of partitions via socle reduction}},
	year = {2025}}

@article{Mon_double_nested,
    AUTHOR = {Monavari, Sergej},
     TITLE = {Double nested {H}ilbert schemes and the local stable pairs
              theory of curves},
   JOURNAL = {Compos. Math.},
  FJOURNAL = {Compositio Mathematica},
    VOLUME = {158},
      YEAR = {2022},
    NUMBER = {9},
     PAGES = {1799--1849}
}

@incollection {Bertin,
    AUTHOR = {Bertin, Jos\'{e}},
     TITLE = {The punctual {H}ilbert scheme: an introduction},
 BOOKTITLE = {Geometric methods in representation theory. {I}},
    SERIES = {S\'{e}min. Congr.},
    VOLUME = {24-I},
     PAGES = {1--102},
 PUBLISHER = {Soc. Math. France, Paris},
      YEAR = {2012},
      ISBN = {978-2-85629-356-0},
   MRCLASS = {14C05 (13A50 14E15 14E16)}, 
MRREVIEWER = {Roberto\ Notari},
}

@unpublished{Kleppe,
    AUTHOR = {Kleppe, Jan O.},
    title = {{The Hilbert-flag scheme, its properties and its connection with the Hilbert scheme.
Applications to curves in 3-space}},
    note = {Unpublished, see \href{https://www.cs.hioa.no/~jank/papers.htm}{https://www.cs.hioa.no/~jank/papers.htm}}
}

@article{CR_framed_motivic,
	author = {Cazzaniga, Alberto and Ricolfi, Andrea T.},
	fjournal = {Mathematische Nachrichten},
	issn = {0025-584X},
	journal = {Math. Nachr.},
	mrclass = {14N35 (14C05)}, 
	number = {6},
	pages = {1096--1112},
	title = {Framed motivic {D}onaldson-{T}homas invariants of small crepant resolutions},
	url = {https://mathscinet.ams.org/mathscinet-getitem?mr=4453951},
	volume = {295},
	year = {2022}}

@article{CRR_higher_rank,
	author = {Cazzaniga, Alberto and Ralaivaosaona, Dimbinaina and Ricolfi, Andrea T.},
	doi = {10.1017/S0305004122000159},
	fjournal = {Mathematical Proceedings of the Cambridge Philosophical Society},
	issn = {0305-0041},
	journal = {Math. Proc. Cambridge Philos. Soc.},
	mrclass = {14N35}, 
	number = {1},
	pages = {97--122},
	title = {Higher rank motivic {D}onaldson-{T}homas invariants of {$\mathbb A^3$} via wall-crossing, and asymptotics},
	url = {https://mathscinet.ams.org/mathscinet-getitem?mr=4523112},
	volume = {174},
	year = {2023}}

@Article{fasola2023tetrahedron,
 Author = {Fasola, Nadir and Monavari, Sergej},
 Title = {Tetrahedron instantons in {Donaldson}-{Thomas} theory},
 FJournal = {Advances in Mathematics},
 Journal = {Adv. Math.},
 ISSN = {0001-8708},
 Volume = {462},
 Pages = {47},
 Note = {110099},
 Year = {2025},
 DOI = {10.1016/j.aim.2024.110099},
 zbMATH = {7972362}
}

@Article{FT_1,
	author = {Feyzbakhsh, Soheyla and Thomas, Richard P.},
Title = {Rank {{\(r\)}} {DT} theory from rank 1},
Journal = {J. Am. Math. Soc.},
 Volume = {36},
 Number = {3},
 Pages = {795--826},
 Year = {2023}
}

@article{FMR_higher_rank,
	author = {Nadir {Fasola} and Sergej {Monavari} and Andrea T. {Ricolfi}},
	fjournal = {{Forum of Mathematics, Sigma}},
	issn = {2050-5094/e},
	journal = {{Forum Math. Sigma}},
	msc2010 = {14N35 14C05},
	note = {e15},
	pages = {1--51},
	publisher = {Cambridge University Press, Cambridge},
	title = {{Higher rank \(K\)-theoretic Donaldson-Thomas theory of points}},
	volume = {9},
	year = {2021}}

@article{MR1616606,
	author = {Cheah, Jan},
	fjournal = {Pacific Journal of Mathematics},
	issn = {0030-8730},
	journal = {Pacific J. Math.},
	mrclass = {14C05 (14C15)}, 
	mrreviewer = {Raquel Mallavibarrena},
	number = {1},
	pages = {39--90},
	title = {Cellular decompositions for nested {H}ilbert schemes of points},
	volume = {183},
	year = {1998}}

@article{Fogarty_Hilb,
	author = {Fogarty, John},
	fjournal = {American Journal of Mathematics},
	issn = {0002-9327},
	journal = {Amer. J. Math.},
	mrclass = {14.20},
	mrreviewer = {S. Lubkin},
	pages = {511--521},
	title = {Algebraic families on an algebraic surface},
	volume = {90},
	year = {1968}}

@article{Gansner_reversed_plane_partitions,
	author = {Gansner, Emden R.},
	doi = {10.1016/0097-3165(81)90041-8},
	fjournal = {Journal of Combinatorial Theory. Series A},
	issn = {0097-3165},
	journal = {J. Combin. Theory Ser. A},
	mrclass = {05A17 (05A15)},
	mrnumber = {607040},
	mrreviewer = {Marshall L. Cates},
	number = {1},
	pages = {71--89},
	title = {The {H}illman-{G}rassl correspondence and the enumeration of reverse plane partitions},
	url = {https://mathscinet-ams-org.proxy.library.uu.nl/mathscinet-getitem?mr=607040},
	volume = {30},
	year = {1981}}

@book{Euler_partitions,
	author = {Euler, Leonhard},
	publisher = {V. 1, Lausanne},
	title = {Introductio in analysin infinitorum},
	year = {1748}}

@Misc{Macaulay,
 Author = {Macaulay, Francis S.},
 Title = {The algebraic theory of modular systems.},
 Year = {1916},
 HowPublished = {Cambridge: {University} press, {XIV} u. 112 {S}. {{\(8^{\circ}\)}}}
}

@Article{Atkin,
 Author = {Atkin, Arthur O. L. and Bratley, Paul and Macdonald, Ian G. and McKay, John K. S.},
 Title = {Some computations for $m$-dimensional partitions},
 FJournal = {Proceedings of the Cambridge Philosophical Society},
 Journal = {Proc. Camb. Philos. Soc.},
 ISSN = {0008-1981},
 Volume = {63},
 Pages = {1097--1100},
 Year = {1967}
}

@article {BFT_flags,
    AUTHOR = {Bonelli, Giulio and Fasola, Nadir and Tanzini, Alessandro},
     TITLE = {Flags of sheaves, quivers and symmetric polynomials},
   JOURNAL = {Forum Math. Sigma},
  FJOURNAL = {Forum of Mathematics. Sigma},
    VOLUME = {12},
      YEAR = {2024},
     PAGES = {Paper No. e74, 51},
      ISSN = {2050-5094},
   MRCLASS = {14C05 (05E05 14D20 14D21 14N35)}, 
       DOI = {10.1017/fms.2024.43},
       URL = {https://doi.org/10.1017/fms.2024.43},
}

@article{cazzaniga2020framed,
	author = {Alberto Cazzaniga and Andrea T. Ricolfi},
	journal = {Math. Z.},
	pages = {745--760},
	title = {{Framed sheaves on projective space and Quot schemes}},
	volume = {300},
	year = {2022}}

@Article{dCM_motives,
 Author = {de Cataldo, Mark A. and Migliorini, Luca},
 Title = {The {Chow} groups and the motive of the {Hilbert} scheme of points on a surface},
 FJournal = {Journal of Algebra},
 Journal = {J. Algebra},
 ISSN = {0021-8693},
 Volume = {251},
 Number = {2},
 Pages = {824--848},
 Year = {2002}
}

@misc{mozgovoy2019motivic,
	author = {Sergey Mozgovoy},
	title = {{Motivic classes of Quot-schemes on surfaces}},
    howpublished = {\href{https://arxiv.org/abs/1911.07561}{ArXiv:1911.07561}},
	year = {2019}}

@article{ricolfi2019motive,
	author = {Andrea T. Ricolfi},
	journal = {J. Math. Pures Appl.},
	pages = {50--68},
	title = {{On the motive of the Quot scheme of finite quotients of a locally free sheaf}},
	volume = {144},
	year = {2020}}

@article{Bifet,
	author = {Bifet, Emili},
	fjournal = {Comptes Rendus de l'Acad\'{e}mie des Sciences. S\'{e}rie I. Math\'{e}matique},
	issn = {0764-4442},
	journal = {C. R. Acad. Sci. Paris S\'{e}r. I Math.},
	mrclass = {14C05 (14L30)},
	mrreviewer = {V. Balaji},
	number = {9},
	pages = {609--612},
	title = {{Sur les points fixes du sch\'{e}ma {${\textrm{Quot}}_{{\mathscr O}^r_X/X/k}$} sous l'action du tore {${\mathbf G}^r_{m,k}$}}},
	url-dummy = {https://mathscinet.ams.org/mathscinet-getitem?mr=1053288},
	volume = {309},
	year = {1989}}

@incollection{Grothendieck_Quot,
	author = {Grothendieck, Alexander},
	booktitle = {S\'{e}minaire {B}ourbaki, Exp. No. 221, no.~6},
	mrclass = {14C05},
	pages = {249--276},
	publisher = {Soc. Math. France},
	title = {Techniques de construction et th\'{e}or{\`e}mes d'existence en g\'{e}om\'{e}trie alg\'{e}brique. {IV}. {L}es sch\'{e}mas de {H}ilbert},
	url-dummy = {https://mathscinet.ams.org/mathscinet-getitem?mr=1611822},
	year = {1961}}

@book{Nakajima,
	author = {Nakajima, Hiraku},
	doi-dummy = {10.1090/ulect/018},
	isbn = {0-8218-1956-9},
	mrclass = {14C05 (17B69)},
	mrreviewer = {Mark Andrea A. de Cataldo},
	pages = {xii+132},
	publisher = {American Mathematical Society, Providence, RI},
	series = {University Lecture Series},
	title = {Lectures on {H}ilbert schemes of points on surfaces},
	url-dummy = {https://mathscinet.ams.org/mathscinet-getitem?mr=1711344},
	volume = {18},
	year = {1999}}

@article {Iarropunctual,
    AUTHOR = {Iarrobino, Anthony},
     TITLE = {Punctual {H}ilbert schemes},
   JOURNAL = {Mem. Amer. Math. Soc.},
  FJOURNAL = {Memoirs of the American Mathematical Society},
    VOLUME = {10},
      YEAR = {1977},
    NUMBER = {188},
     PAGES = {viii+112},
      ISSN = {0065-9266,1947-6221},
   MRCLASS = {14C05}, 
MRREVIEWER = {Mary\ Schaps},
       DOI = {10.1090/memo/0188},
       URL = {https://doi.org/10.1090/memo/0188},
}

@incollection{Hilb_11,
	author = {Douvropoulos, Theodosios and Jelisiejew, Joachim and N{\o}dland, Bernt Ivar Utst{\o}l and Teitler, Zach},
	booktitle = {Combinatorial algebraic geometry},
	mrclass = {14C05 (14D15)},
	mrreviewer = {Paolo Lella},
	pages = {321--352},
	publisher = {Fields Inst. Res. Math. Sci., Toronto, ON},
	series = {Fields Inst. Commun.},
	title = {The {H}ilbert scheme of 11 points in {$\mathbb A^3$} is irreducible},
	url-dummy = {https://mathscinet.ams.org/mathscinet-getitem?mr=3752506},
	volume = {80},
	year = {2017}}

@book{Sernesi-Deformations,
	author = {Sernesi, Edoardo},
	isbn = {978-3-540-30608-5; 3-540-30608-0},
	mrclass = {14D15 (14B05 14B10 14B12 14B20 14B25 14D20)},
	mrreviewer = {Marko Roczen},
	pages = {xii+339},
	publisher = {Springer-Verlag, Berlin},
	series = {Grundlehren der Mathematischen Wissenschaften},
	title = {Deformations of algebraic schemes},
	url-dummy = {https://mathscinet.ams.org/mathscinet-getitem?mr=2247603},
	volume = {334},
	year = {2006}}

@Book{EISENBUD,
 Author = {Eisenbud, David},
 Title = {Commutative algebra. {With} a view toward algebraic geometry},
 FSeries = {Graduate Texts in Mathematics},
 Series = {Grad. Texts Math.},
 ISSN = {0072-5285},
 Volume = {150},
 ISBN = {3-540-94269-6; 3-540-94268-8},
 Year = {1995},
 Publisher = {Berlin: Springer-Verlag}
}

@Article{zbMATH06941785,
 Author = {Kuznetsov, Alexander and Shinder, Evgeny},
 Title = {Grothendieck ring of varieties, {D}- and {L}-equivalence, and families of quadrics},
 FJournal = {Selecta Mathematica. New Series},
 Journal = {Sel. Math., New Ser.},
 ISSN = {1022-1824},
 Volume = {24},
 Number = {4},
 Pages = {3475--3500},
 Year = {2018}
}

@article{MR18,
	author = {Riccardo Moschetti and Andrea T. Ricolfi},
	doi-dummy = {https://doi.org/10.1016/j.jalgebra.2018.09.028},
	issn = {0021-8693},
	journal = {{J. Algebra}},
	pages = {471--489},
	title = {On coherent sheaves of small length on the affine plane},
	url-dummy = {http://www.sciencedirect.com/science/article/pii/S0021869318305477},
	volume = {516},
	year = {2018}}

@article {ALESSIONESTED,
    AUTHOR = {Ramkumar, Ritvik and Sammartano, Alessio},
     TITLE = {Rational {S}ingularities of {N}ested {H}ilbert {S}chemes},
   JOURNAL = {Int. Math. Res. Not. IMRN},
  FJOURNAL = {International Mathematics Research Notices. IMRN},
      YEAR = {2024},
    NUMBER = {2},
     PAGES = {1061--1122},
      ISSN = {1073-7928,1687-0247},
   MRCLASS = {14 (13)}, 
       DOI = {10.1093/imrn/rnac365},
       URL = {https://doi.org/10.1093/imrn/rnac365},
}

@article{double-nested-1,
 author = {Graffeo, Michele and Lella, Paolo and Monavari, Sergej and Ricolfi, Andrea T. and Sammartano, Alessio},
 title = {The geometry of double nested {Hilbert} schemes of points on curves},
 fjournal = {Transactions of the American Mathematical Society},
 journal = {Trans. Am. Math. Soc.},
 issn = {0002-9947},
 volume = {378},
 number = {9},
 pages = {6013--6047},
 year = {2025},
 doi = {10.1090/tran/9247},
}

@Article{MacMahon-statistics,
 Author = {Amanov, Alimzhan and Yeliussizov, Damir},
 Title = {{MacMahon}'s statistics on higher-dimensional partitions},
 FJournal = {Forum of Mathematics, Sigma},
 Journal = {Forum Math. Sigma},
 ISSN = {2050-5094},
 Volume = {11},
 Pages = {23},
 Note = {Id/No e63},
 Year = {2023}
}

@article{Gott2,
	author = {G{\"{o}}ttsche, Lothar},
	journal = {Math. Ann.},
	pages = {193--207},
	title = {{The Betti numbers of the Hilbert scheme of points on a smooth projective surface}},
	volume = {286},
	year = {1990}}

@phdthesis{Rydh1,
	address = {Stockholm},
	author = {Rydh, David},
	school = {KTH},
	title = {{Families of cycles and the Chow scheme}},
	year = {2008}}

@article{Bria1,
	author = {Joel Brian{\c{c}}on},
	doi-dummy = {10.1007/BF01390164},
	fjournal = {{Invent. Math.}},
	issn = {0020-9910; 1432-1297/e},
	journal = {{Invent. Math.}},
	msc2010 = {14C05},
	pages = {45--89},
	publisher = {Springer, Berlin/Heidelberg},
	title = {{Description de $\textrm{Hilb}^n\mathbb C\{x,y\}$}},
	volume = {41},
	year = {1977},
	zbl = {0353.14004}}

@Article{zbMATH05493515,
 Author = {Sabir M. {Gusein-Zade} and Ignacio {Luengo} and Alejandro {Melle-Hern\'andez}},
 Title = {On the power structure over the {Grothendieck} ring of varieties and its applications},
 FJournal = {Proceedings of the Steklov Institute of Mathematics},
 Journal = {Proc. Steklov Inst. Math.},
 ISSN = {0081-5438},
 Volume = {258},
 Pages = {53--64},
 Year = {2007}
}

@incollection{LooijengaMM,
	author = {Looijenga, Eduard},
	booktitle = {{S\'eminaire Bourbaki. Volume 1999/2000. Expos\'es 865--879}},
	msc2010 = {14F42 14G10 14G05 14J10 14D05},
	pages = {267--297, ex},
	publisher = {Paris: Soci\'et\'e Math\'ematique de France},
	title = {{Motivic measures}},
	year = {2002},
	zbl = {0996.14011}}

@book {RPStanley-VolI,
    AUTHOR = {Stanley, Richard P.},
     TITLE = {Enumerative combinatorics. {V}olume 1},
    SERIES = {Cambridge Studies in Advanced Mathematics},
    VOLUME = {49},
   EDITION = {Second},
 PUBLISHER = {Cambridge University Press, Cambridge},
      YEAR = {2012},
     PAGES = {xiv+626},
      ISBN = {978-1-107-60262-5},
   MRCLASS = {05-02 (05A15 06-02)}, 
}

@article{GLMHilb,
	author = {Sabir M. {Gusein-Zade} and Ignacio {Luengo} and Alejandro {Melle-Hern\'andez}},
	doi-dummy = {10.1307/mmj/1156345599},
	fjournal = {{Michigan Mathematical Journal}},
	issn = {0026-2285},
	journal = {{Mich. Math. J.}},
	msc2010 = {14C05 18F30},
	number = {2},
	pages = {353--359},
	publisher = {University of Michigan, Department of Mathematics, Ann Arbor, MI},
	title = {{Power structure over the Grothendieck ring of varieties and generating series of Hilbert schemes of points}},
	volume = {54},
	year = {2006},
	zbl = {1122.14003}}

@article{ESHilb,
	author = {Geir {Ellingsrud} and Stein Arild {Str{\o}mme}},
	doi-dummy = {10.1007/BF01389419},
	fjournal = {{Invent. Math.}},
	issn = {0020-9910; 1432-1297/e},
	journal = {{Invent. Math.}},
	msc2010 = {14C05 14N05},
	pages = {343--352},
	publisher = {Springer, Berlin/Heidelberg},
	title = {{On the homology of the Hilbert scheme of points in the plane}},
	volume = {87},
	year = {1987},
	zbl = {0625.14002}}

@incollection{DenefLoeser1,
	author = {Denef, Jan and Loeser, Fran{\c c}ois},
	booktitle = {{3rd European congress of mathematics (ECM), Barcelona, Spain, July 10--14, 2000. Volume I}},
	isbn = {3-7643-6417-3/hbk; 3-7643-6419-X/set},
	msc2010 = {14B05 14B10 14A15 14F43},
	pages = {327--348},
	publisher = {Basel: Birkh\"auser},
	title = {{Geometry on arc spaces of algebraic varieties}},
	year = {2001},
	zbl = {1079.14003}}

@article {multigraded,
    AUTHOR = {Haiman, Mark and Sturmfels, Bernd},
     TITLE = {Multigraded {H}ilbert schemes},
   JOURNAL = {J. Algebraic Geom.},
  FJOURNAL = {Journal of Algebraic Geometry},
    VOLUME = {13},
      YEAR = {2004},
    NUMBER = {4},
     PAGES = {725--769},
      ISSN = {1056-3911,1534-7486},
   MRCLASS = {14C05 (13D40)}, 
MRREVIEWER = {Roy\ M.\ Skjelnes},
       DOI = {10.1090/S1056-3911-04-00373-X},
       URL = {https://doi.org/10.1090/S1056-3911-04-00373-X},
}

@article {ELEMENTARY,
    AUTHOR = {Jelisiejew, Joachim},
     TITLE = {Elementary components of {H}ilbert schemes of points},
   JOURNAL = {J. Lond. Math. Soc. (2)},
  FJOURNAL = {Journal of the London Mathematical Society. Second Series},
    VOLUME = {100},
      YEAR = {2019},
    NUMBER = {1},
     PAGES = {249--272},
      ISSN = {0024-6107,1469-7750},
   MRCLASS = {14C05 (13D10 14Q15)}, 
MRREVIEWER = {A.\ Prabhakar\ Rao},
       DOI = {10.1112/jlms.12212},
       URL = {https://doi.org/10.1112/jlms.12212},
}

@article {EMSALEM,
    AUTHOR = {Emsalem, Jacques},
     TITLE = {G\'{e}om\'{e}trie des points \'{e}pais},
   JOURNAL = {Bull. Soc. Math. France},
  FJOURNAL = {Bulletin de la Soci\'{e}t\'{e} Math\'{e}matique de France},
    VOLUME = {106},
      YEAR = {1978},
    NUMBER = {4},
     PAGES = {399--416},
      ISSN = {0037-9484},
   MRCLASS = {14D15 (12H05)}, 
MRREVIEWER = {K.\ Wolffhardt},
       URL = {http://www.numdam.org/item?id=BSMF_1978__106__399_0},
}

@incollection {Geramita,
    AUTHOR = {Geramita, Anthony V.},
     TITLE = {Inverse systems of fat points: {W}aring's problem, secant
              varieties of {V}eronese varieties and parameter spaces for
              {G}orenstein ideals},
 BOOKTITLE = {The {C}urves {S}eminar at {Q}ueen's, {V}ol. {X} ({K}ingston,
              {ON}, 1995)},
    SERIES = {Queen's Papers in Pure and Appl. Math.},
    VOLUME = {102},
     PAGES = {2--114},
 PUBLISHER = {Queen's Univ., Kingston, ON},
      YEAR = {1996},
      ISBN = {0-88911-717-9},
   MRCLASS = {13D40 (13C13 14M05)}, 
MRREVIEWER = {Martin\ Kreuzer},
}

@book {Iarrobook,
    AUTHOR = {Iarrobino, Anthony and Kanev, Vassil},
     TITLE = {Power sums, {G}orenstein algebras, and determinantal loci},
    SERIES = {Lecture Notes in Mathematics},
    VOLUME = {1721},
      NOTE = {Appendix C by Iarrobino and Steven L. Kleiman},
 PUBLISHER = {Springer-Verlag, Berlin},
      YEAR = {1999},
     PAGES = {xxxii+345},
      ISBN = {3-540-66766-0},
   MRCLASS = {14M12 (13C40 13H10 13N10 14C05)}, 
MRREVIEWER = {Juan\ C.\ Migliore},
       DOI = {10.1007/BFb0093426},
       URL = {https://doi.org/10.1007/BFb0093426},
}

@InCollection{CARLINI,
 Author = {Carlini, Enrico},
 Title = {Reducing the number of variables of a polynomial},
 BookTitle = {Algebraic geometry and geometric modeling. Based on the workshop, Nice-Sophia Antipolis, France, September 27--29, 2004},
 ISBN = {3-540-33274-X},
 Pages = {237--247},
 Year = {2006},
 Publisher = {Berlin: Springer}
}

@article{CONSTELLATION,
author = {Graffeo, Michele},
title = {{Moduli spaces of $\mathbb Z/k\mathbb Z$-constellations over $\mathbb A^2$}},
journal = {Communications in Contemporary Mathematics},
volume = {27},
number = {03},
pages = {2450019},
year = {2025},
doi = {10.1142/S0219199724500196},
URL = { 
    
        https://doi.org/10.1142/S0219199724500196
    
},
eprint = { 
    
        https://doi.org/10.1142/S0219199724500196
 
}

}

@article {UPDATES,
    AUTHOR = {Franco Giovenzana and Luca Giovenzana and Michele Graffeo and Paolo Lella},
     TITLE = {{U}nexpected but recurrent phenomena for {Q}uot and {H}ilbert schemes of points},
   JOURNAL = {Rend. Semin. Mat. Univ. Politec. Torino},
  FJOURNAL = {Rendiconti del Seminario Matematico. Universit\`a{} e
              Politecnico Torino},
    VOLUME = {82},
      YEAR = {2024},
    NUMBER = {1},
     PAGES = {145--170}
}

@misc{MOTIVES,
      title={{The motive of the Hilbert scheme of points in all dimensions}}, 
      author={Graffeo, Michele and Monavari, Sergej and Moschetti, Riccardo and Ricolfi, Andrea T.},
      howpublished = {To appear in Proc. Lond. Math. Soc.},
      year={2026}
}

@article{BULOIS,
    AUTHOR = {Bulois, Micha\"el and Evain, Laurent},
     TITLE = {Nested punctual {H}ilbert schemes and commuting varieties of
              parabolic subalgebras},
   JOURNAL = {J. Lie Theory},
  FJOURNAL = {Journal of Lie Theory},
    VOLUME = {26},
      YEAR = {2016},
    NUMBER = {2},
     PAGES = {497--533}
}

@Misc{zbMATH03473265,
 Author = {Knuth, Donald E.},
 Title = {The art of computer programming. {Vol}. 3: {Sorting} and searching},
 Year = {1973},
 HowPublished = {Addison-{Wesley} {Series} in {Computer} {Science} and {Information} {Processing}. {Reading}, {Mass}. etc.: {Addison}-{Wesley} {Publishing} {Company}. {XI}, 722 pages}
}

@Article{hook-length-formula,
 Author = {Frame, J. S. and Robinson, Gilbert de B. and Thrall, R. M.},
 Title = {The hook graphs of the symmetric group},
 FJournal = {Canadian Journal of Mathematics},
 Journal = {Can. J. Math.},
 ISSN = {0008-414X},
 Volume = {6},
 Pages = {316--324},
 Year = {1954}
}

@book {HARRIS,
    AUTHOR = {Harris, Joe},
     TITLE = {Algebraic geometry},
    SERIES = {Graduate Texts in Mathematics},
    VOLUME = {133},
      NOTE = {A first course},
 PUBLISHER = {Springer-Verlag, New York},
      YEAR = {1992},
     PAGES = {xx+328},
      ISBN = {0-387-97716-3},
   MRCLASS = {14-01}, 
MRREVIEWER = {Liam\ O'Carroll},
       DOI = {10.1007/978-1-4757-2189-8},
       URL = {https://doi.org/10.1007/978-1-4757-2189-8},
}

@article {Negut_flags,
    AUTHOR = {{Negu\c{t}}, Andrei},
     TITLE = {Moduli of flags of sheaves and their {$K$}-theory},
   JOURNAL = {Algebr. Geom.},
  FJOURNAL = {Algebraic Geometry},
    VOLUME = {2},
      YEAR = {2015},
    NUMBER = {1},
     PAGES = {19--43},
      ISSN = {2313-1691,2214-2584},
   MRCLASS = {14J60 (14D21 14F05 19E08)}, 
MRREVIEWER = {P.\ E.\ Newstead},
       DOI = {10.14231/AG-2015-002},
       URL = {https://doi.org/10.14231/AG-2015-002},
}

@misc{2step,
      title={{New components of Hilbert schemes of points and 2-step ideals}}, 
      author={Franco Giovenzana and Luca Giovenzana and Michele Graffeo and Paolo Lella},
      year={2025},
      eprint={2507.02789},
      archivePrefix={arXiv},
      primaryClass={math.AG},
      url={https://arxiv.org/abs/2507.02789}, 
}

\bigskip

\bigskip
\noindent
{\small{Nadir Fasola \\
\address{SISSA, Via Bonomea 265, 34136, Trieste (Italy)} \\
\href{mailto:nfasola@sissa.it}{\texttt{nfasola@sissa.it}}
}}

\bigskip
\noindent
{\small{Michele Graffeo \\
\address{Dipartimento di Matematica, Politecnico di Milano, Piazza Leonardo da Vinci 32, 20133, Milan (Italy)} \\
\href{mailto:michele.graffeo@polimi.it}{\texttt{michele.graffeo@polimi.it}}
}}
 
\bigskip
\noindent
{\small{Danilo Lewa\'nski \\
\address{Università di Trieste, Dipartimento MIGe, Via Valerio 12/1, 34127, Trieste (Italy) and Istituto Nazionale di Fisica Nucleare (INFN), sezione di Trieste (Italy)} \\
\href{mailto:danilo.lewanski@units.it}{\texttt{danilo.lewanski@units.it}}
}}

\bigskip
\noindent
{\small Andrea T. Ricolfi (corresponding author) \\
\address{SISSA, Via Bonomea 265, 34136, Trieste (Italy)} \\
\href{mailto:aricolfi@sissa.it}{\texttt{aricolfi@sissa.it}}}

 \end{document}